\documentclass[a4paper, 12pt]{article}

\usepackage[a4paper, margin=2.5cm]{geometry}
\usepackage[utf8]{inputenc}
\usepackage[english]{babel}
\usepackage{graphicx}
\usepackage[colorlinks=true, citecolor=blue]{hyperref}
\usepackage{mathtools}
\usepackage{amsfonts}
\usepackage{amsthm}
\usepackage{parskip}
\usepackage{bbm}
\usepackage{natbib}
\usepackage{authblk}

\everymath{\displaystyle}

\makeatletter
\def\thm@space@setup{%
  \thm@preskip=\parskip \thm@postskip=0pt
}
\makeatother

\newtheorem{theorem}{Theorem}
\newtheorem{lemma}[theorem]{Lemma}
\newtheorem{definition}[theorem]{Definition}
\newtheorem{corollary}[theorem]{Corollary}

\title{Neural Field Models with Transmission Delays and Diffusion}
\author[1,*]{Len Spek}
\author[2,1]{Yuri A. Kuznetsov}
\author[1,2]{Stephan A. van Gils}
\affil[1]{Department of Applied Mathematics, University of Twente, Enschede, The Netherlands}
\affil[2]{Department of Mathematics, Utrecht University, Utrecht, The Netherlands}
\affil[*]{Corresponding author: l.spek@utwente.nl}
\date{}
\begin{document}
\maketitle

\begin{abstract}
\noindent A neural field models the large scale behaviour of large groups of neurons. We extend results of \citet{van_gils_local_2013} and \citet{dijkstra_pitchforkhopf_2015} by including a diffusion term into the neural field, which models direct, electrical connections. We extend known and prove new sun-star calculus results for delay equations to be able to include diffusion and explicitly characterise the essential spectrum. For a certain class of connectivity functions in the neural field model, we are able to compute its spectral properties and the first Lyapunov coefficient of a Hopf bifurcation. By examining a numerical example, we find that the addition of diffusion suppresses non-synchronised steady-states, while favouring synchronised oscillatory modes.
\end{abstract}

\section{Introduction}
In the study of neurological disease, non-invasive imaging techniques are often used to get an understanding of the structure and functioning of the brain on intermediate scales. As they give a course-grained view of the neuronal activity, mean-field models are a natural fit to describe the observed dynamics \citep{jirsa_field_1996,jirsa_spatiotemporal_2002}. In this paper we will use a neural field model with gap-junctions, electrical connections between neurons, which are thought to be related to observed synchronisation of neural tissue in Parkinson's disease \citep{schwab_pallidal_2014,schwab_synchronization_2014}. We will study the effect of gap junctions on the dynamics of the model. We will mainly focus on the stability of steady-states, periodic oscillations and the bifurcations which lead to a qualitative change in behaviour. 

To properly address the difference in time-scales between gap-junctions and synaptic connections, we use a neural field with transmission delays for the synaptic connections. This leads to a complicated model which is infinite-dimensional and has spatially-distributed delays. The dynamical theory for such models is not readily available. In this paper, we address the analytic problems which arise from these abstract delay differential equations. 

We use the sun-star calculus as the basic functional analytic tool to cast the equation in the variation-of-constants form. We exploit the results by  \citet{janssens_class_2019, janssens_class_2020} that allow the linear part of the equation, without the delays, to be unbounded, as is the case for the diffusion operator.

\subsection{Background}
Neural field models try to bridge the gap between single neurons models \citep{hodgkin_quantitative_1952}, and whole brain models \citep{sanz_leon_virtual_2013}, by modelling the qualitative behaviour of large groups of neurons. In the seminal work of \citet{wilson_excitatory_1972,wilson_mathematical_1973}, they modelled two populations of excitatory and inhibitory neurons and analysed the dynamical properties of the resulting model. A neural field uses spatial and temporal averaging of the membrane voltage of a population of neurons. The synaptic connections are modelled by a convolution of a connectivity kernel and a nonlinear activation function. This leads to a set of two a integro-differential equations with delays.

These models have been simplified by \citet{amari_dynamics_1977} by combining the excitatory and inhibitory populations into a single population and made more realistic by \citet{nunez_brain_1974} by including transmission delays. These delays arise from the finite propagation speed of action potentials across an axon and the delay due to dendritic integration. There has been considerable interest in the role of these delays in the spatiotemporal dynamics. \citep{jirsa_time_2007,liley_spatially_2002,hutt_pattern_2003,hutt_analysis_2005,hutt_spontaneous_2007,hutt_local_2008,coombes_waves_2005,ermentrout_large_1980,ermentrout_mathematical_2010}. Further modelling work by Coombes, Venkov and collaborators show the usefulness of these neural fields for understanding neural activity \citep{coombes_stephen_delays_2009,coombes_large-scale_2010,coombes_tutorial_2014,venkov_dynamic_2007}.

Roxin and collaborators first did a bifurcation analysis for neural fields with a single fixed delay \citep{roxin_role_2005,roxin_rate_2006,roxin_how_2011}. Faugeras and collaborators investigated the stability properties of stationary solutions of these neural fields with distance dependent delays \citep{faye_theoretical_2010,veltz_localglobal_2010,veltz_stability_2011,veltz_center_2013} using a functional analytic approach based on formal projectors. In \citet{van_gils_local_2013} it was shown that the neural fields can be studied as abstract delay differential equations to which the sun-star framework can be applied. They used this to compute normal form coefficients for bifurcations of equilibria. \citet{dijkstra_pitchforkhopf_2015} expanded their analysis to Pitchfork-Hopf bifurcations and \citet{visser_standing_2017} analysed a neural field with delays on a spherical domain. We will build on the work of \citet{van_gils_local_2013} and \citet{dijkstra_pitchforkhopf_2015} by introducing gap-junctions into the neural field model and studying the resulting bifurcations and dynamics. 

Gap-junctions are electrical connections between neurons, which directly exchange ions through a connexin-protein. This is in contrast to synaptic connections, where a potential is induced across the synapse by neurotransmitters. These gap-junctions are thought to be related to Parkinson's disease by synchronising neurons in the {\it globus pallidus} \citep{schwab_pallidal_2014,schwab_synchronization_2014}. Gap-junctions can be modelled as a simple diffusion process \citep{coombes_tutorial_2014}. There have been some attempts to incorporate gap-junctions into networks of coupled neurons \citep{amitai_spatial_2002,laing_exact_2015,ostojic_synchronization_2009}, but to our knowledge not yet within a proper neural field model. 

\subsection{Theoretical framework}
As mentioned before, we use the sun-star calculus for delay differential equations to formally analyse these neural field models with transmission delays. This mathematical theory for delay differential equations was constructed in the book by \citet{diekmann_delay_1995} and the references therein. This theory uses the space $X^{\odot}$, pronounced $X$-sun, which is the largest subspace of strong continuity of the adjoint semigroup. It allows us to employ the classical Fredholm alternative, which plays a key role in the computation of the normal form coefficients. As a result, many of the mathematical techniques developed for the analysis of ODE's, such as the center manifold reduction and the Hopf-bifurcation theorem, can be generalised for these abstract delay differential equations. 

Recently, \citet{janssens_class_2019,janssens_class_2020} has begun expanding the sun-star calculus to the case where the linear part, which contains no delays, is an unbounded operator. This allows us to study both the neural field with and without diffusion in the same framework. This unifying theory then allows us to fill in the gap in the proofs of \citet{van_gils_local_2013}, while obtaining the same results for a neural field with diffusion.

There are also other theoretical frameworks possible. The first approach to develop a geometric theory for delay equations along the lines of ODEs was proposed by \citet{hale_theory_1971} who used formal adjoint operators. Formal adjoint operators were also used by \citet{faria_normal_1995,faria_normal_1995-1} and \citet{arino_normal_2006} to study Hopf and Bogdanov-Takens bifurcations.  \citet{wu_theory_2012} used the formal adjoint method to study reaction-diffusion systems with delays and prove the necessary theorems for bifurcation analysis.

There is a difference whether to take as a starting point an abstract integral equation, like we do, or an abstract ODE like in the integrated semigroup approach of \citet{magal_center_2009,magal_semilinear_2009} and \citet{liu_normal_2014}. Integrated semigroups have been used to deal with classical delay differential equations as abstract ODEs with non-dense domains. By classical we here mean that the state space is $\mathbbm{R}^n $. In the case of the neural field equations we consider, the state space is an abstract Banach space. It might very well be possible that the formalism of integrated semigroups is general enough to cover this as well, but as far as we know, this has not been done as yet. We prefer the sun-star formalism, as it allows us to work with the variation-of-constants formula in the state space $X$, albeit after an excursion in the bigger space $X^{\odot\ast}$. In addition, the projectors are based on duality pairing and the classical Fredholm alternative, while in the integrated semigroup formalism the projectors are based on a formal inner product \citep{liu_projectors_2008}.

There are also two approaches to compute normal form coefficients. In the first approach, the abstract ODE is split into a finite dimensional and an infinite dimensional one. By decoupling these step by step, the center manifold is rectified and the equation on it is normalised \citep{magal_center_2009,magal_semilinear_2009,liu_normal_2014}. In the second approach, which we follow, we parametrise the center manifold and assume that the finite dimensional ODE on it is in normal form. As the delay differential equation has an abstract state space, this ODE is also an abstract ODE. The Taylor coefficients of the center manifold are obtained in a step by step procedure that simultaneously gives us the coefficients of the normal form \citep{coullet_amplitude_1983,elphick_simple_1987}. In this way, the sun-star calculus approach leads to explicit, compact and easy to evaluate expression for the normal form coefficients \citep{janssens_normalization_2010}. These coefficients are obtained using the true duality pairing, for which the classical Fredholm alternative holds. Of course, the resulting formulas are equivalent but the approach we adopted is more straightforward.

In the sun-star calculus we choose to model the neural field as a continuous function in space. \citet{faye_theoretical_2010} and \citet{veltz_stability_2011} instead choose to use the $L^2$-functions, based on the work by \citet{webb_functional_1976} and \citet{batkai_semigroups_2001,batkai_semigroups_2005}. This leads to some mathematical complications dealing with the smoothness of the non-linearity, as laid out previously in section 2.4 of \citep{van_gils_local_2012}. This was later rectified by \citet{veltz_erratum_2015}. Moreover, from a physiological point of view, it is not clear why the potential of the neural field should be merely square integrable, instead of continuous.

Finally, we want to briefly comment on the need of a theoretical framework to study these neural fields. Software packages, such as DDE-BIFTOOL \citep{engelborghs_numerical_2002}, can perform numerical bifurcation analysis of delay equations. However, they can not directly be applied to these delayed integro-differential equations. While a discretised model can be studied with these software packages, there is no guarantee that the dynamical properties converge to those of the full neural field. In this work, the formulas of the normal form coefficients are exact and can be evaluated to arbitrary precision

In this paper we will build on the work of  \citet{janssens_class_2019,janssens_class_2020} and prove the necessary theorems to use the sun-star calculus to study our neural field model with diffusion and without diffusion. We will then derive the spectrum and resolvent of a neural field with delays, diffusion and a connectivity kernel of a sum of exponentials. Finally we will compute the first Lyapunov coefficient of a Hopf-bifurcation and verify our results by simulating the full neural field numerically.

\subsection{Modelling}\label{sec:modelling}
In this section we will derive the neural field model with transmission delays and gap junctions. This is largely based on a derivation by  \citet{ermentrout_large_1980}. 

We start with a collection of neurons $i=1,2,3,\cdots$ and denote the (somatic) potential of neuron $i$ at time $t$ by $u_i(t)$ and its firing rate by $f_i(t)$. We assume there is a nonlinear dependence of $f_i$ on $u_i$ given by
\[f_i(t) = S_i(u_i(t))\]
We define $\Phi_{i,j}(t)$ to be the postsynaptic potential appearing on postsynaptic cell $i$ due to a single spike from presynaptic cell $j$. We assume a linear summation of the postsynaptic potentials, so the total potential received at the soma due to the synaptic connection between cell $i$ and $j$ can be modelled as
\[G_{i,j}(t)= \int_{-\infty}^t \Phi_{i,j}(t-s)f_j(s-\tau_{i,j})\,ds\]
where $\tau_{i,j}$ is the delay due to the finite propagation  speed of action potentials along an axon and other factors such as dendritic integration. We define $\Psi_{i}(t)$ to be the potential appearing in neuron $i$ due to a gap-junction current $I_{i,gap}(t)$. The resulting model for $u_i$ becomes
\begin{equation}\label{eq:model1}
u_i(t) = \Psi_i(t) + \sum_j \int_{-\infty}^t \Phi_{i,j}(t-s)S_j(u_j(s-\tau_{i,j}))\,ds
\end{equation}
We can reduce this integral equation if we have a model for $\Phi$ and $\Psi$. For cell $i$, let us consider a passive membrane with a time constant $1/\alpha_i$, a resistance $R_i$ and an injected postsynaptic current $I_{i,j,syn}(t)$
\[\frac{1}{\alpha_i} \frac{d \Phi_{i,j}}{dt} + \Phi_{i,j} = R_i I_{i,j,syn}(t)\]
and similarly when a gap-junction current is injected
\[\frac{1}{\alpha_i} \frac{d \Psi_i}{dt} + \Psi_i = R_i I_{i,gap}(t)\]
If we now apply the Laplace transform $\mathcal{L}$ to equation \eqref{eq:model1}, we get
\[\left(\frac{s}{\alpha_i}+1\right)\mathcal{L}(u_i)(s) = R_i \mathcal{L}(I_{i,gap})(s) + R_i \sum_j \mathcal{L}(I_{i,j,syn})(s) \mathcal{L}(S_j(u_j(\cdot -\tau_{i,j})))(s)\]
We assume that the synaptic dynamics are dominated by the time-scale of the membrane. This means we can reduce $I_{i,j,syn}(t)$ to $w_{i,j} \delta(t)$, where $\delta$ is the Dirac-delta distribution and $w_{i,j}$ represents the strength of the synaptic connection, where a negative value corresponds to inhibition. Taking the inverse Laplace transform results in a system of differential equations
\begin{equation}\label{eq:model2}
\left(\frac{1}{\alpha_i} \frac{d}{dt}+1\right) u_i(t) = R_i I_{i,gap}(t) + R_i \sum_j w_{i,j} S_j(u_j(t-\tau_{i,j}))
\end{equation}

We want to model this network of cells by a neural field. Suppose we have a sequence of similar neurons $i=1,2,\cdots, M$ on the interval $\Omega=[-1,1]$ and we model the gap-junctions as a simple resistor between adjacent neurons we arrive at the formula
\begin{equation}\label{eq:model3}
\left(\frac{1}{\alpha} \frac{d}{dt}+1\right) u_i(t) = R g (u_{i-1}(t) - 2 u_i(t) + u_{i+1}(t)) + R \sum_j w_{i,j} S(u_j(t-\tau_{i,j}))
\end{equation}
We will now take the limit as $M \rightarrow \infty$ while scaling $g$ by $M^2$ and $w_{i,j}$ by $1/M$, to find our neural field model
\begin{equation}\label{eq:model4}
\frac{\partial u}{\partial t}(t,x)= d \frac{\partial^2 u}{\partial x^2}(t,x) - \alpha u(t,x)+ \alpha \int_{\Omega}J(x,x')S(u(t-\tau(x,x'),x'))\,dx'
\end{equation}
We haven't specified yet what happens with the gap-junctions at the boundary of our domain. It is natural to assume that no current leaks away at the boundaries, which corresponds to Neumann boundary conditions in the neural field:
\[\frac{\partial u}{\partial x}(t,\pm 1) =0\]

\subsection{Overview}
This paper is divided into three parts, each of which can mostly be read independently. 

In section \ref{sec:duality}, we construct the sun-star calculus for abstract delay differential equations and derive the variation-of-constants formula. In particular we prove a novel characterisation for sun-reflexivity. Furthermore we consider linearisation, the corresponding spectrum and a normal form derivation for Hopf bifurcation of the nonlinear equations. In appendix \ref{sec:diffusion} we elaborate on the case when the unbounded linear operator is the diffusion operator. We expect the reader to be familiar with the basics of the sun-star framework in the book by \citet{diekmann_delay_1995}.

In section \ref{sec:single_pop} we derive formulas for the eigenvalues and eigenvectors for a neural field with a connectivity defined by a sum of exponentials. We also explicitly construct the solution to the resolvent problem for this class of neural field models. 

The section \ref{sec:num} we do a numerical study for a neural field model with specific parameter values. We compute the first Lyapunov coefficient for the Hopf bifurcation and investigate how it is influenced by the diffusion term. We will also investigate the emergence of periodic behaviour using numerical simulations of the neural field.

\section{Abstract Delay Differential Equations in the Sun-Star Framework}
\label{sec:duality}
In this section we will first develop the sun-star calculus for a large class of abstract delay differential equations, \eqref{ADDE}. This leads to a variation-of-constants formulation of the \eqref{ADDE}. Next we study the linearisation and obtain results on the spectrum. Finally we construct a method for computing the first Lyapunov coefficient for a Hopf bifurcation of nonlinear equations. We build on the theory developed by \citet{janssens_class_2019}, who considers a class of abstract delay differential equations with a possibly unbounded linear part.

Consider two Banach spaces $Y$ and $X = C([-h,0];Y)$ over $\mathbb{R}$ or $\mathbb{C}$. Let $S$ be a strongly continuous semigroup on $Y$ with its generator $B$ and let $G: X\rightarrow Y$ be a (nonlinear) globally Lipschitz-continuous operator. Note that the assumption that the semigroup $S$ is compact is not necessary, in contrast to what is assumed by \citet{wu_theory_2012}. 

We introduce now our main object of study
\begin{equation}
\begin{cases}
\dot{u}(t)=B u(t)+G(u_t)\\
u_0=\varphi\in X
\end{cases}
\label{ADDE} \tag{ADDE}
\end{equation}
Here $u_t \in X$, where $u_t(\theta)= u(t+\theta)$ for $t\geq 0$ and $\theta \in [-h,0]$.

In the remaining sections we are mainly interested in the case where $B$ is a diffusion operator acting in the space of continuous functions $Y=C([-a,a];\mathbb{R})$. We have summarised the relevant properties of the diffusion operator in Appendix \ref{sec:diffusion}. However the theorems which are proven in this section hold for any operator $B$ that generates a strongly continuous semigroup $S$ on $Y$. This fills in some technical details missing in \citep{van_gils_local_2013}, where $B=-\alpha I$, which does not generate a compact semigroup.

On $X$ we consider the strongly continuous semigroup $T_0$ defined by 
\begin{equation}
\label{eq:T0}
(T_0(t)\varphi)(\theta):=\begin{cases}
\varphi(t+\theta) &t+\theta\in [-h,0]\\
S(t+\theta)\varphi(0) &t+\theta>0
\end{cases}
\end{equation}
Here $\varphi\in X, t\geq 0$ and $\theta \in [-h,0]$. This semigroup is related to the problem for $G\equiv 0$, i.e.
\begin{align}
\begin{cases}
\dot{v}(t)=B v(t) &\text{for } t>0\\
v_0=\varphi &\text{for } t \in [-h,0]
\end{cases}
\label{eq:gzero}
\end{align}
The solution of problem \eqref{eq:gzero} is then given by $v_t:=T_0(t)\varphi$. 

\begin{lemma} \citep[Theorem VI.6.1]{engel_one-parameter_1999}
The generator $A_0$ of the semigroup \sloppy{${T_0}$} is given by
\begin{equation}\label{eq:A0}
A_0 \varphi = \dot{\varphi}, \quad D(A_0)=\{\varphi \in C^1([-h,0];Y)| \varphi(0)\in D(B) \text{ and }  \dot{\varphi}(0)=B\varphi(0)\}
\end{equation}
\end{lemma}

We will interpret the \eqref{ADDE} as problem \eqref{eq:gzero} with some nonlinear perturbation $G:X\rightarrow Y$ and use a variation-of-constants formula in $X$ to obtain results about the perturbed problem, such as normal form coefficients for local bifurcations. As $G$ maps $X$ into $Y$, we would like to embed $Y$ in a natural way into $X$. A naive approach would be to use a delta-function as an embedding. However, this embedding is not bounded, so the domain of $A_0$ would not be preserved under perturbation. This is indeed the case, as the rule for extending a function beyond its original domain, i.e. $\dot{\varphi}(0)=B\varphi(0)$, is incorporated in $D(A_0)$. Hence adding a perturbation to the rule for extension changes the domain of the generator. A way out is to embed this problem into a larger space. A natural choice would be $Y \times X$, where we have a continuous embedding $\ell:Y \rightarrow Y\times \{0\}$ and we can separate the extension and translation part of $A_0$ into $Y \times \{0\}$ and $\{0\} \times X$ respectively.

\begin{figure}[ht!]
\includegraphics[width=0.95\linewidth]{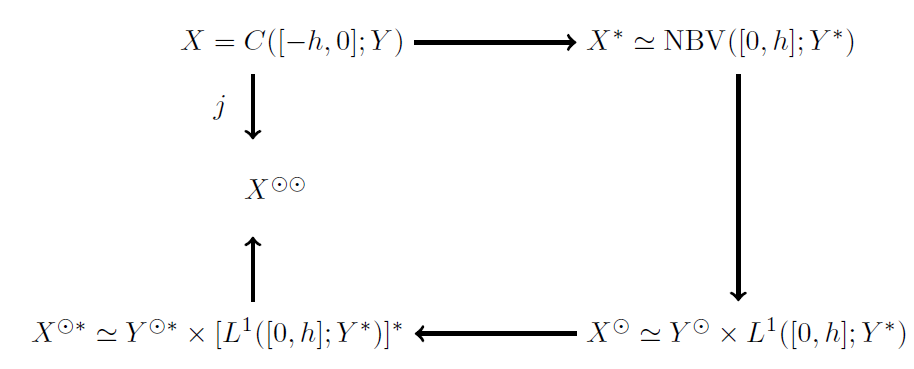}
\caption{A schematic representation of the various Banach spaces in sun-star calculus \cite{janssens_class_2019}}\label{Figure1}
\end{figure}

More formally we use the sun-star calculus as developed in the book by \citet{diekmann_delay_1995} to construct the space $X^{\odot *}$, which contains the space $Y\times X$. We will first restrict the dual space $X^*$ to the sun space $X^\odot$, on which $T_0^*$ is strongly continuous. Then taking the dual we obtain the dual space $X^{\odot *}$. It is convenient to present the relationship of the various spaces schematically in the following `duality' diagram, see figure \ref{Figure1}. 

\section*{Figures}

\subsection{Characterisation of the sun-dual}
Using a generalisation of the Riesz Representation Theorem, we can find a representation of $X^*$, the dual space of $X$ \citep{gowurin_uber_1936}. It can be represented as $NBV([0,h];Y^*)$, the space of functions $f:[0,h]\rightarrow Y^*$ of bounded variation on $[0,h]$, normalised such that $f(0)=0$ and $f$ is right continuous on $(0,h)$. The (complex valued) duality pairing between $X$ and $X^*$ is given by the Riemann-Stieltjes integral, for $\varphi \in X$ and $f \in X^*$
\[\langle f, \varphi \rangle := \int_0^h \varphi(-\theta) \,df(\theta) \]
Results of scalar functions of bounded variation and the corresponding Riemann-Stieltjes integral can be extended to $Y$-valued functions \citep{gowurin_uber_1936}.

It is possible to find an explicit representation of the adjoint operator $A_0^*$ and its corresponding domain $D(A_0^*)$. The adjoint operator exists and is unique as the domain $D(A_0)$ is dense.
\begin{theorem}
The domain of $A_0^*$ is given by
\begin{equation}
\begin{split}
D(A_0^*) := \{ f\in NBV([0,h];Y^*)|  \text{ there exists } y^*\in D(B^*) \text{ and } g\in NBV([0,h];Y^*)&\\ 
\text{ with } g(h)=0 \text{ such that } f(t) = y^* \chi_0(t) + \int _0^t g(\theta)\,d\theta \}&
\end{split}
\label{eq:da0}
\end{equation}
and the action of $A_0^*$ is given by $A_0^* f = B^*y^*\chi_0+g$, where $\chi_0 = \mathbbm{1}_{(0,h]}$, i.e the characteristic function of $(0,h]$.
\label{thm:a0*}
\end{theorem}
\textit{Proof.} We first prove the inclusion $\subseteq$ for the domain $D(A_0^*)$.
Let $f\in D(A_0^*)$ and $\varphi \in D(A_0)$. Without loss of generality we can write $A_0^* f = c\chi_0 + g$, where $c \in Y^*$ and $g \in NBV([0,h];Y^*)$ and $g(h)=0$. Using the integration by parts formulas for Riemann-Stieltjes integrals \citep[Appendix H]{bartle_modern_2001} we obtain
\begin{equation}
\begin{split}
\int_0^h \dot{\varphi}(-\theta) \,df(\theta) &= \langle f, A_0 \varphi \rangle = \langle A_0^* f, \varphi \rangle = \langle c\chi_0 + g, \varphi \rangle \\
&= \int_0^h \varphi(-\theta) \,d(c\chi_0(\theta)) + \int_0^h \varphi(-\theta) \,dg(\theta)\\
&= \langle c, \varphi(0) \rangle + \langle g(\theta), \varphi(-\theta) \rangle |_0^h - \int_0^h g(\theta) \,d\varphi(-\theta)\\
&= \langle c, \varphi(0) \rangle + \int_0^h \langle g(\theta), \dot{\varphi}(-\theta) \rangle \,d\theta
\end{split} 
\label{eq:da0*_2}
\end{equation}
We will now want to use some limiting argument. However, the Riemann-Stieltjes integral lacks good convergence properties. In the scalar case, we could interpret this integral as a Lebesque-Stieltjes integral, which has better convergence properties. For a general Banach space $Y$ and continuous integrands, the equivalent would be the Bartle integral \citep{singer_linear_1957,bartle_general_1956}. The Bartle integral has an equivalent theorem to the Lebesque Dominated Convergence Theorem. For uniformly bounded, pointwise converging sequences we can interchange the limit and the integral. \citep[Theorem 6]{bartle_general_1956}. 

For some $0<s<t\leq h$ and $y\in Y$, we may choose $(\dot{\varphi}_n)_{n\in \mathbb{N}}$ as a uniformly bounded sequence in $X$ such that $\dot{\varphi_n}(0)=\varphi_n(0)=0$ and it converges pointwise to $y \mathbbm{1}_{[-t,-s]}$, i.e. the characteristic function of $[-t,-s]$. We then substitute $\varphi$ for $\varphi_n$ in \eqref{eq:da0*_2}
\[\int_0^h \dot{\varphi}_n(-\theta) \,df(\theta) = \int_0^h \langle g(\theta), \dot{\varphi}_n(-\theta) \rangle \,d\theta\]
Taking the limit as $n \rightarrow \infty$, we get using the dominated convergence of the Bartle integral that
\begin{align*}
\int_0^h y \mathbbm{1}_{[-t,-s]}(-\theta) \,df(\theta) &= \int_0^h \langle g(\theta), y \mathbbm{1}_{[-t,-s]}(-\theta) \rangle \,d\theta\\
\langle f(t)-f(s), y \rangle &= \int_s^t \langle g(\theta) , y \rangle \,d\theta
\end{align*}
Since $y$ was arbitrary, we infer that
\[f(t)= f(s) + \int_s^t g(\theta) \,d\theta\]
Letting $s\downarrow 0$, we obtain for $t\in [0,h]$
\[f(t)= y^*\chi_0(t) + \int_0^t g(\theta) \,d\theta\]
where $y^* = \lim_{s\downarrow 0} f(s)$. Now we substitute this formula for $f$ into $\langle f, A_0 \varphi \rangle$ and use integration by parts and the fact that $\dot{\varphi}(0)=B\varphi(0)$ to find that
\begin{align*}
\langle f, A_0 \varphi \rangle &= \langle y^* , \dot{\varphi}(0) \rangle + \int_0^h \langle g(\theta), \dot{\varphi}(-\theta) \rangle \,d\theta\\
&= \langle y^* , B\varphi(0) \rangle + \int_0^h \langle g(\theta), \dot{\varphi}(-\theta) \rangle \,d\theta
\end{align*}
We compare this to equation \eqref{eq:da0*_2}
\[\langle f, A_0 \varphi \rangle = \langle c, \varphi(0) \rangle + \int_0^h \langle g(\theta), \dot{\varphi}(-\theta) \rangle \,d\theta\]
Since $\varphi(0)$ can be chosen arbitrary, $\langle y^* , B\varphi(0) \rangle = \langle c, \varphi(0) \rangle$ implies that $c \in D(B^*)$ and $c=B^* y^*$. 

Finally we prove the other inclusion $\supseteq $ for the domain $D(A_0^*)$ and simultaneously obtain the formula for the action of $A_0^*$. Let $f$ be of the form in $\eqref{eq:da0}$, then by the above computations we find that
\[\langle f, A_0 \varphi \rangle= \langle B^*y^* , \varphi(0) \rangle + \int_0^h \varphi(-\theta) \,dg(\theta) = \langle B^*y^*\chi_0+g, \varphi \rangle\]
 \qed

We can characterise the sun-dual $X^{\odot}$ as the subspace of $X^*$ where $T_0^*$ is strongly continuous or equivalently $X^{\odot} = \overline{D(A_0^*)}$, where the closure is with respect to the norm on $X^*$. Similarly we can characterise the sun-dual $Y^{\odot}$ as the subspace of $Y^*$ where $B^*$ is strongly continuous or equivalently $Y^{\odot} = \overline{D(B^*)}$, where the closure is with respect to the norm on $Y^*$. In case $B$ is the diffusion operator, see \ref{sec:diffusion} for an explicit characterisation of $Y^\odot$.

The following theorem can be proved by showing that $T_0^*$ is strongly continuous on some set $E$ given by \eqref{eq:E}, that $D(A_0^*)\subseteq E$, and that $E$ is closed.

\begin{theorem}\citep[Theorem 1 and Remark 4]{janssens_class_2019}\label{thm:xsun}
The space $X^\odot$, the sun-dual of $X$ with respect to $T_0$, is given by the set
\begin{equation}\label{eq:E}
\begin{split}
\{ f:[0,h]\rightarrow Y^*| & \text{ there exists } y^\odot\in Y^\odot \text{ and } g\in L^1([0,h];Y^*)\\
& \text{ such that } f(t) = y^\odot \chi_0(t) + \int _0^t g(\theta)\,d\theta \}
\end{split}
\end{equation}
Furthermore, the map $\iota: Y^\odot \times L^1([0,h];Y^*)\rightarrow X^\odot$ defined by
\begin{equation}\label{eq:isometry}
\iota (y^\odot, g)(t) := y^\odot \chi_0(t) +\int_0^t g(\theta)\,d\theta \quad
\forall t\in [0,h]
\end{equation}
is an isometric isomorphism.
\end{theorem}
From now on we will identify $X^\odot$ with $Y^\odot \times L^1 ([0,h];Y^*)$. The corresponding duality pairing between $X$ and $X^\odot$ is then given by 
\begin{equation}
\langle \varphi^\odot, \varphi \rangle := \langle y^\odot, \varphi(0) \rangle + \int_0^h \langle g(\theta), \varphi(-\theta) \rangle \,d\theta
\end{equation}

Now we can describe the action of $T_0^\odot$ and $A_0^\odot$, the restrictions of the operators $T_0^*$ and $A_0^*$ to the subspace $X^\odot$.
\begin{definition}
The strongly continuous semigroup $T_1$ on $L^1([0,h];Y^*)$ is defined as
\begin{equation}
(T_1(t)g)(\theta) := \begin{cases}
g(t+\theta) &t+\theta\in [0,h]\\
0 &t+\theta > h \end{cases}
\end{equation}
\end{definition}
\begin{theorem}\citep[Theorem 1]{janssens_class_2019}
\label{thm:T0sun}
For the action of $T_0^\odot$ on $X^\odot$ we have 
\begin{equation}
T_0^\odot(t)(y^\odot,g) := \left(S^\odot (t) y^\odot + \int_0^{\text{min}(t,h)} S^*(t-\theta)g(\theta)\,d\theta\,,\, T_1(t)g\right)
\label{eq:2.6}
\end{equation}
where the integral is the weak$^*$ Lebesque integral with values in $Y^\odot$.
\end{theorem} 

\begin{theorem}\label{thm:a0sun}
For the sun-dual of $A_0$ on $X^\odot$ we have that
\begin{equation}\label{eq:da0sun}
\begin{split}
D(A_0^\odot)= \{ (y^\odot,g) | g \in AC([0,h];Y^*) \text{ with } g(h)=0, y^\odot\in D(B^*)&\\  \text{ and } B^* y^\odot + g(0)\in Y^\odot \}&
\end{split}
\end{equation}
and $A_0^\odot (y^\odot,g) = (B^* y^\odot + g(0),\dot{g})$, with $\dot{g}$ a function in $L^1([0,h];Y^*)$ such that 
\begin{equation}
g(t) = g(0) + \int_0^t \dot{g}(\theta)\,d\theta
\end{equation}
for $t \in [0,h]$.
\end{theorem}
\textit{Proof.} By definition
\[D(A_0^\odot) := \{ \varphi^\odot \in X^\odot | \iota\varphi^\odot\in D(A_0^*), A_0^* \iota\varphi^\odot \in \iota (X^\odot) \}\]
and $\iota A_0^\odot \varphi^\odot = A_0^* \iota \varphi^\odot$. We first prove the equivalence of the definition and \eqref{eq:da0sun}.

Let $\varphi^\odot=(y^\odot,g) \in X^\odot$ such that $\iota\varphi^\odot\in D(A_0^*)$ and $A_0^* \iota\varphi^\odot \in \iota (X^\odot)$. Recall that the embedding $\iota$ is given by \eqref{eq:isometry}
\[\iota \varphi^\odot(t) = y^\odot \chi_0(t) + \int_0^t g(\theta)\,d\theta\]
From Theorem \ref{thm:a0*}, we can conclude that $\iota\varphi^\odot\in D(A_0^*)$ implies that $y^\odot \in D(B^*)$ and $g\in NBV([0,h];Y^*)$ with $g(h)=0$. As $A_0^* \iota \varphi^\odot = B^* y^\odot \chi_0+g \in \iota(X^\odot)$, Theorem \ref{thm:xsun} implies that $B^*y^\odot+g(0+)\in Y^\odot$ and we can write $g$ as
\[g(t) = g(0+)\chi_0 + \int_0^t \dot{g}(\theta) \,d\theta\]
where $g(0+)= \lim_{t \downarrow 0} g(t)$ and $\dot{g}$ some function in $L^1([0,h];Y^\odot)$. Hence $g$ is absolutely continuous on $(0,h]$. As $g$ is an $L^1$-function (class), we may redefine $g(0):=g(0+)$ to get a absolutely continuous function on $[0,h]$.

Conversely, let $\varphi^\odot=(y^\odot,g) \in X^\odot$ such that it is in the right hand side of \eqref{eq:da0sun}. From Theorem \ref{thm:a0*} and the fact that $g-g(0)\in NBV([0,h];Y^*)$, we conclude that $\iota\varphi^\odot\in D(A_0^*)$ and that $A_0^* \iota\varphi^\odot= (B^* y^\odot + g(0))\chi_0 + g$. As $g$ is absolutely continuous and $B^* y^\odot + g(0)\in Y^\odot$, this implies that $A_0^* \iota\varphi^\odot = \iota(B^* y^\odot + g(0),\dot{g})\in \iota (X^\odot )$. Hence, $A_0^\odot \varphi^\odot = (B^* y^\odot + g(0),\dot{g})$. \qed 

\subsection{Characterisation of the sun-star space}
We can represent $X^{\odot *}$, the dual of $X^\odot$, as $Y^{\odot *} \times (L^1([0,h];Y^*)^*$, where $Y^{\odot *}$ is the dual of $Y^{\odot}$. In case  $B$ is the diffusion operator, $Y^{\odot *}$ is explicitly characterised in Appendix \ref{sec:diffusion}.

In general, $(L^1([0,h];Y^*)^*$ cannot be identified with $L^\infty([-h,0];Y^{**})$. However, the latter space can be embedded into the former. 
\begin{theorem}\label{thm:Xsunstar}
\citep[Remark 1.4.18, Theorem 1.4.19]{cazenave_introduction_1998}
There exists an isometric embedding of $L^\infty([-h,0];Y^{**})$ into $(L^1([0,h];Y^*)^*$ with the duality pairing 
\[\langle \varphi,g \rangle = \int_{0}^h \langle \varphi(-\theta),g(\theta) \rangle \,d\theta \]
for $g\in L^1([0,h];Y^*)$ and $\varphi \in L^\infty([-h,0];Y^{**})$.

Moreover, $(L^1([0,h];Y^*)^*$ can be identified with $L^\infty([-h,0];Y^{**})$ if and only if $Y^{**}$ has the Radon-Nikodym property.
\end{theorem}
\begin{lemma}\label{lem:reflexive}
(Dunford-Pettis) If $Y$ is reflexive then it has the Radon-Nikodym property.
\end{lemma}

We can embed both $Y$ and $X$ into $Y \times X$ which is a subspace of $Y^{\odot *} \times L^\infty([-h,0];Y^{**})$. The canonical embedding $j: X \rightarrow X^{\odot *}$ is defined as $\langle j \varphi,\varphi^\odot \rangle = \langle \varphi^\odot, \varphi \rangle$. The continuous embedding $\ell: Y \rightarrow X^{\odot *}$ is defined as $\ell = (j_Y y, 0)$, where $j_Y$ is the canonical embedding of $Y$ into $Y^{\odot *}$ \citep{janssens_class_2019}. It is possible to find an explicit representation of $j$. 

\begin{lemma}\label{lem:embedj}
For $\varphi \in X$, $j\varphi = (j_Y \varphi(0),\varphi)$. Moreover, $j$ is a continuous embedding and $j^{-1}: j(X) \rightarrow X$ is bounded. $T_0^{\odot *}(t)j = j T_0(t)$, consequently $j(X)$ is contained in $X^{\odot \odot}$, which is the subspace of $X^{\odot *}$ on which $T_0^{\odot *}$ is strongly continuous.
\end{lemma}
\textit{Proof.} Let $\varphi \in X$ and $\varphi^\odot=(y^\odot,g) \in X^\odot$, then
\begin{align*}
\langle j \varphi,\varphi^\odot \rangle &= \langle \varphi^\odot, \varphi \rangle \\
&= \langle y^\odot, \varphi(0) \rangle + \int_0^h \langle g(\theta), \varphi(-\theta) \rangle \,d\theta \\
&= \langle j_Y \varphi(0), y^\odot \rangle + \int_0^h \langle \varphi(-\theta), g(\theta) \rangle \,d\theta \\
&= \langle (j_Y\varphi(0),\varphi), \varphi^\odot \rangle 
\end{align*}
Hence $j\varphi = (j_Y \varphi(0),\varphi)$. The other statements are generally known to hold for the canonical embedding of $X$ into $X^{\odot *}$, \citep[Appendix II, Cor. 3.16, Prop. 3.17]{diekmann_delay_1995}. \qed

As we don't have an explicit norm or measure on $(L^1([0,h];Y^*)^*$ we cannot say anything in general about $A_0^{\odot *}$. However, it is possible to find a representation of $A_0^{\odot *}$ restricted to the space \sloppy{${Y^{\odot *} \times L^\infty([-h,0];Y^{**})}$}.
\begin{theorem}\label{thm:a0sunstar}
For $(y^{\odot *},\varphi) \in X^{\odot *}$ the following statements are equivalent
\begin{enumerate}
\item $(y^{\odot *},\varphi) \in D(A_0^{\odot *})$ and $A_0^{\odot *}(y^{\odot *},\varphi)\in Y^{\odot *} \times L^\infty([-h,0];Y^{**})$
\item $\varphi$ has an a.e. derivative $\dot{\varphi} \in L^\infty([-h,0];Y^{**})$ for which
\[\varphi(t) = y^{\odot *} - \int_t^0 \dot{\varphi}(\theta) \,d\theta \]
and $\varphi(0)=y^{\odot *}\in D(B^{\odot *})$. 
\end{enumerate}
In this case the action of $A_0^{\odot *}$ is given by $A_0^{\odot *}(y^{\odot *},\varphi) = (B^{\odot *} y^{\odot *}, \dot{\varphi})$
\end{theorem}
\textit{Proof.} Let $(y^{\odot *},\varphi) \in D(A_0^{\odot *})$ such that $A_0^{\odot *}(y^{\odot *},\varphi)= (\gamma, \psi)\in Y^{\odot *} \times L^\infty([-h,0];Y^{**})$ and let $(y^\odot,g) \in D(A_0^\odot)$. We have that
\begin{equation}\label{eq:a0sunstar}
\begin{split}
\langle y^{\odot *}, B^* y^\odot +g(0) \rangle + \langle \varphi, \dot{g} \rangle 
&= \langle (y^{\odot *},\varphi), A_0^{\odot}(y^\odot,g) \rangle \\
&= \langle A_0^{\odot *}(y^{\odot *},\varphi),(y^\odot,g) \rangle \\
&= \langle \gamma, y^\odot \rangle + \int_{0}^h \langle \psi(-\theta),g(\theta) \rangle \,d\theta 
\end{split}
\end{equation}
Let $\Phi \in L^\infty([-h,0];Y^{**})$ such that
\[\Phi(t) = \Phi(0) - \int_t^0 \psi(\theta)\, d\theta\]
Then by Lemma \ref{lem:int_by_parts} and Theorem \ref{thm:a0sun}, i.e. $g(h)=0$, we can rewrite \eqref{eq:a0sunstar} as
\begin{equation}\label{eq:a0sunstar2}
\langle y^{\odot *}, B^* y^\odot +g(0) \rangle = \langle \gamma, y^\odot \rangle + \langle \Phi(0), g(0) \rangle + \langle \Phi - \varphi, \dot{g} \rangle
\end{equation}

Taking $g\equiv 0$ we get that $\langle y^{\odot *}, B^* y^\odot \rangle = \langle \gamma, y^\odot \rangle$ for all $y^\odot \in Y^\odot$ such that $B^* y^\odot \in Y^\odot$ by Theorem \ref{thm:a0sun}. Hence $y^\odot \in D(B^\odot)$, which implies that $y^{\odot *} \in D(B^{\odot *})\subseteq Y^{\odot \odot}$ and $\gamma = B^{\odot *}y^{\odot *}$. As $Y^{\odot \odot}$ can be embedded in $Y^{**}$ \citep[Corollary 4.2]{clement_perturbation_1986}, we find that $y^{\odot *} \in Y^{**}$. Furthermore, we have for all $y^{\odot *} \in D(B^{\odot *})$ and $y^\odot \in D(B^*)$ \citep[Theorem 4.3]{clement_perturbation_1986}.
\[\langle B^{\odot *} y^{\odot *}, y^\odot \rangle = \langle y^{\odot *}, B^* y^\odot \rangle\]
Alternatively, we take $\Phi(0) = y^{\odot *}$, $g(0) = -\int_0^h \dot{g}(\theta)d\theta$ and $y^\odot \in D(B^*)$ such that $B^{*}y^\odot + g(0) \in Y^\odot$. Then \eqref{eq:a0sunstar2} reduces to $\langle \Phi -\varphi, \dot{g} \rangle =0$ for all $\dot{g}\in L^1([0,h];Y^{*})$, hence $\Phi \equiv \varphi$.

Conversely, let $(y^{\odot *},\varphi) \in Y^{\odot *} \times L^\infty([-h,0];Y^{**})$, where $\varphi(0)=y^{\odot *}\in D(B^{\odot *})$ and $\varphi$ has an a.e. derivative $\dot{\varphi} \in L^\infty([-h,0];Y^{**})$ for which
\[\varphi(t) = y^{\odot *} - \int_t^0 \dot{\varphi}(\theta) \,d\theta \]
Then again using Lemma \ref{lem:int_by_parts} we get that for any $(y^\odot,g) \in D(A_0^\odot)$
\begin{align*}
\langle (y^{\odot *},\varphi), A_0^{\odot}(y^\odot,g) \rangle 
&= \langle y^{\odot *}, B^* y^\odot +g(0) \rangle + \int_{0}^h \langle \varphi(-\theta), \dot{g}(\theta) \rangle \,d\theta\\
&= \langle y^{\odot *}, B^* y^\odot \rangle + \int_{0}^h \langle \dot{\varphi}(-\theta),g(\theta) \rangle \,d\theta \\
&= \langle B^{\odot *}y^{\odot *}, y^\odot \rangle + \int_{0}^h \langle \dot{\varphi}(-\theta),g(\theta) \rangle \,d\theta \\
&= \langle (B^{\odot *} y^{\odot *}, \dot{\varphi}), (y^\odot, g) \rangle
\end{align*}
Hence $A_0^{\odot *}(y^{\odot *},\varphi) = (B^{\odot *} y^{\odot *}, \dot{\varphi}) \in Y^{\odot *} \times L^\infty([-h,0];Y^{**})$. \qed

\begin{corollary}\label{cor:a0sunstar}
For $\varphi \in X$ the following statements are equivalent
\begin{enumerate}
\item $j\varphi \in D(A_0^{\odot *})$ and $A_0^{\odot *}j\varphi \in Y^{\odot *} \times L^\infty([-h,0];Y^{**})$
\item $j_Y\varphi(0) \in D(B^{\odot *})$ and $\varphi$ has an a.e. derivative $\dot{\varphi} \in L^\infty([-h,0];Y)$
\end{enumerate}
In this case, the action of $A_0^{\odot *}$ is given by $A_0^{\odot *}j\varphi = (B^{\odot *} j_Y \varphi(0), \dot{\varphi})$. 
\end{corollary}
\textit{Proof.} This follows immediately from Theorem \ref{thm:a0sunstar} and Lemma \ref{lem:embedj}. \qed

Note that for $A_0^{\odot *}$ the rule for extension, $\dot{\varphi}(0) = B\varphi(0)$, is no longer included in the domain of $A_0^{\odot *}$, but is represented in the action of $A_0^{\odot *}$, which resolves the problem with $A_0$ stated at the beginning of this section.

The previous theorem allows us to formulate an equivalence between the sun-reflexivity of $X$, i.e. $X^{\odot \odot} = j(X)$ and the ordinary reflexivity of $Y$, i.e. $Y^{**} = j_Y(Y)$
\begin{theorem}\label{thm:sun-reflexivity}
$X$ is sun-reflexive with respect to $T_0$ if and only if $Y$ is reflexive.
\end{theorem}
\textit{Proof.}
Suppose that $Y$ is reflexive. Then by Theorem \ref{thm:Xsunstar} and Lemma \ref{lem:reflexive}, $X^{\odot *}$ can be represented as $Y^{\odot *} \times L^\infty([-h,0];Y)$ and hence the full domain of $A_0^{\odot *}$ is given by Theorem \ref{thm:a0sunstar}
\[D(A_0^{\odot *})= \{(y^{\odot *},\varphi)\in X^{\odot *} | \varphi(0)=y^{\odot *}\in D(B^{\odot *}), \varphi \text{ has an a.e. derivative} \}\]
We use that $X^{\odot \odot}$ is the closure of $D(A_0^{\odot *})$ with respect to the norm on $X^{\odot *}$. First the closure of $D(B^{\odot *})$ with respect to the $Y^{\odot *}$-norm results in the space $Y^{\odot \odot}$. As reflexivity implies sun-reflexivity \citep[Corollary 2.5]{van_neerven_reflexivity_1990}, we have that $Y^{\odot \odot}= j_Y(Y)$. Next we note that $C^1$ functions are dense in the continuous functions and $C^0$ is closed with respect to the $L^\infty$-norm. Hence we conclude that
\[X^{\odot \odot} = \{(y^{\odot \odot},\varphi)\in j_Y(Y) \times C([-h,0];Y)|\varphi(0)=y^{\odot \odot}\}= j(X)\]

Conversely, suppose that $Y$ is not reflexive. From Theorem \ref{thm:Xsunstar}, $Y^{\odot *} \times L^\infty([-h,0];Y)$ is a subset of $X^{\odot *}$ and hence
\[\{(y^{\odot *},\varphi)\in X^{\odot *} | \varphi(0)=y^{\odot *}\in D(B^{\odot *}), \varphi \text{ has an a.e. derivative}\} \subseteq D(A_0^{\odot *}) \]
Taking the norm closure of both sides, we conclude that 
\[\{(y^{\odot \odot},\varphi)\in Y^{\odot \odot} \times C([-h,0];Y^{**})|\varphi(0)=y^{\odot \odot}\}\subseteq X^{\odot \odot}\]
As $Y$ is not reflexive, $C([-h,0];Y)$ is a proper subset of $C([-h,0];Y^{**})$. Hence $j(X)$ is a proper subset of $X^{\odot \odot}$, so $X$ is not sun-reflexive. \qed

In case $B$ is the diffusion operator, we use that $Y$ is the space of continuous functions. As this is a non-reflexive Banach space, $X$ in this case is not sun-reflective.

\subsection{Variation-of-constants formulation}
As the space $X^{\odot *}$ solves the problems mentioned in the beginning of this section, we can formulate a variation-of-constants formula for the \eqref{ADDE} as an abstract integral equation 
\begin{equation}
u_t=T_0(t)\varphi+j^{-1}\int_0^t T_0^{\odot *}(t-\tau)\ell G(u_\tau)\,d\tau
\label{AIE} \tag{AIE}
\end{equation}
Here the embeddings $j$ and $\ell$ are as defined in the previous section. As the integrand of \eqref{AIE} takes values in $X^{\odot *}$, the integral is taken to be a weak$^*$ integral. It is possible to show that the integral maps to the range of $j(X)$ and hence the \eqref{AIE} is well-defined.

\begin{lemma} \citep[Proposition 8]{janssens_class_2019}
\label{lem:2.8}
Let $u \in C(\mathbb{R}^+,Y)$ be given, then 
\begin{equation}
\label{eq:A05}
\int_0^t T_0^{\odot\ast} (t-\tau) \ell u(\tau) \,d\tau = j\psi \quad \forall t\geq 0
\end{equation}
where 
\begin{equation}
\psi(\theta):=\int_0^{\max\{(t+\theta), 0\}}S(t-\tau +\theta)u(\tau)\,d\tau \quad \forall \theta\in [-h,0]
\end{equation}
Moreover, 
\begin{equation}\label{eq:3.7}
\| \psi \| \leq M e^{\omega h} \frac{e^{\omega t}-1}{\omega}\sup_{0\leq \tau \leq t}\|u(\tau)\| \quad \forall t\geq 0
\end{equation}
where $M,\omega>0$ are such that $\|S(t)\| \leq M e^{\omega t}$ for all $t\geq 0$.
\end{lemma}

The Banach Fixed Point Theorem in combination with the bound in \eqref{eq:3.7} gives the existence of a unique global solution of \eqref{AIE}. 
\begin{corollary}\citep[Corollary 9]{janssens_class_2019}
\label{cor:2.9}
Let $G:X\rightarrow Y$ be globally Lipschitz continuous. For every initial condition 
$\varphi \in X$ there exists a unique solution $v\in C(\mathbb{R}_+,X)$ such that $u_t=v(t)$ satisfies \eqref{AIE} for all $t \geq 0$.
\end{corollary}

We would like to show that this unique solution of the \eqref{AIE} can be translated over to a (classical) solution of the \eqref{ADDE}. However, this is in general not the case when $B$ is unbounded. Therefore we recall a weaker solution concept from \citet{wu_theory_2012}. 

\begin{definition} \label{def:classical_solution}
A function $u\in C([-h,\infty);Y)$ is called a \textbf{classical solution} of \eqref{ADDE} if $u$ is continuously differentiable on $\mathbb{R}_+$, $u(t)\in D(B)$ for all $t\geq 0$ and $u$ satisfies the \eqref{ADDE}.
\end{definition}

\begin{definition} \label{def:mild solution}
A function $u\in C([-h,\infty);Y)$ is called a \textbf{mild solution} of \eqref{ADDE} if $u_0 =\varphi$ and $u$ satisfies
\begin{equation}\label{eq:A08}
u(t)=S(t)\varphi(0) +\int_0^t S(t-\tau) G(u_\tau)\,d\tau\quad \forall t\geq 0
\end{equation}
\end{definition}

Note that Definition \ref{def:classical_solution} is quite restrictive as only specific initial conditions $\varphi\in X$ are admissible. There is the following correspondence between classical and mild solutions of \eqref{ADDE}
\begin{lemma}\citep[Theorem 2.1.4]{wu_theory_2012}
A classical solution of \eqref{ADDE} is also a mild solution of \eqref{ADDE}

Conversely when $G$ has a globally Lipschitz continuous Fréchet derivative and $\varphi \in C^1([-h,0];Y)$, $\varphi(0) \in D(B)$ and $\dot{\varphi}(0) = B\varphi(0) + G(\varphi)$ then a mild solution of \eqref{ADDE} is also a classical solution of \eqref{ADDE}.
\end{lemma}

Note that Theorem \ref{thm:generatorA} below implies that the conditions in the second statement, starting with conversely, are equivalent to the condition that $\varphi \in D(A)$. 

It is possible to construct a one-to-one correspondence between solutions of \eqref{AIE} and mild solutions of \eqref{ADDE}.
\begin{theorem}\citep[Theorem 16]{janssens_class_2019}
\label{thm:2.15}
Let $\varphi \in X$ be an initial condition. The following two statements hold.
\begin{enumerate}
\item
Suppose that $u$ is a mild solution of \eqref{ADDE}. Define $v:\mathbb{R}_+\rightarrow X$ by
\[ v(t):=u_t\quad \forall t\geq 0 \]
Then $v$ is a solution of \eqref{AIE}.
\item
Suppose that $v$ is a solution of \eqref{AIE}. Define $u:[-h,\infty)\rightarrow Y$ by
\[ u(t):=\begin{cases}
\varphi(t)\quad &-h\leq t\leq 0\\
v(t)(0) \quad & t\geq 0
\end{cases}\]
Then $u$ is a mild solution of \eqref{ADDE}. 
\end{enumerate}
\end{theorem}

\begin{corollary}
Suppose $G$ is a globally Lipschitz operator and it has a globally Lipschitz Fréchet derivative then for all $\varphi \in C^1([-h,0];Y)$ with $\varphi(0) \in D(B)$ and $\dot{\varphi}(0) = B\varphi(0) + G(\varphi)$, there exists a unique classical solution of the \eqref{ADDE}. 
\end{corollary}

\subsection{Linearisation}
We want to investigate the behaviour near a fixed point. We will show that for the linearised problem, we can perturb the semigroup $T_0$ with generator $A_0$ to a semigroup $T$ with generator $A$. In the next section we will investigate the spectral properties of $A$.

Linearising equation \eqref{ADDE} near a fixed point $u$, which we take without loss of generality to be $u\equiv 0$, results in the linear problem $\eqref{LP}$.
\begin{equation}
\begin{cases}
\dot{u}(t)=Bu(t)+DG(0)u_t\\
u_0=\varphi\in X
\end{cases}
\label{LP} \tag{LINP}
\end{equation}

As with the general nonlinear problem we can define an abstract integral equation.
\begin{equation}\tag{AIE}
u_t = T_0(t) \varphi + j^{-1}\int_0^t T_0^{\odot *}(t-s)L u_t\,ds
\end{equation}
where $L := \ell DG(0)$. Then due to Lemma \ref{lem:2.8} and Corollary \ref{cor:2.9} we can define the strongly continuous semigroup $T(t)\varphi:=u_t$ when $DG(0)$ is globally Lipschitz.

\begin{lemma}\citep[Theorem 19]{janssens_class_2019}
\label{thm:3.6}
Let $DG(0)$ be globally Lipschitz continuous, then there exists a unique strongly continuous semigroup $T$ on $X$ such that 
\begin{equation}
\label{eq:3.10}
T(t)\varphi = T_0(t) \varphi + j^{-1} \int_0^t T_0^{\odot\ast} LT(\tau) \varphi \,d\tau
\end{equation}
for all $\varphi \in X$ and for all $t\geq 0$.
\end{lemma}
The strongly continuous semigroup $T$ has a generator $A$. We want to establish how the perturbed generator $A$ relates to the original generator $A_0$, which can be done using the sun-star framework. A technical detail which we need to check is that the sun dual space $X^\odot$ is the same with respect to $T$ and $T_0$.

\begin{lemma}\citep[Proposition 20]{janssens_class_2019}
$X^{\odot}$ is also the maximal subspace of strong continuity of the adjoint semigroup $T^*$ on $X^*$. The adjoint generator $A^*$ is given by 
\begin{equation}
A^*=A_0^* + L^* \text{ with } D(A^*)=D(A_0^*)
\end{equation}
and the generator $A^{\odot}$ of the $T^\odot$ is given by
\begin{equation}
A^\odot=A_0^\odot + L^\odot \text{ with } D(A^\odot)=D(A_0^\odot)
\end{equation}
Finally $X^{\odot \odot}$ is also the maximal subspace of strong continuity of the sun-star semigroup $T^{\odot \odot}$.
\label{lem:A*}
\end{lemma}
One could think that we could extend this argument and show that $D(A^{\odot\ast})=D(A_0^{\odot\ast})$ and
$A^{\odot\ast}=A_0^{\odot\ast}+L j^{-1}$. However, this is not the case when we lack sun-reflexivity, i.e. $X^{\odot \odot}\neq j(X)$. We can circumvent these problems by restricting the domain to $j(X)$

\begin{lemma}\citep[Proposition 22]{janssens_class_2019}
\label{lem:3.12} It holds that
\begin{equation}
\label{eq:3.15}
D(A^{\odot\ast})\cap j (X) = D(A_0^{\odot\ast}) \cap j (X)
\end{equation}
and $A^{\odot\ast}=A_0^{\odot\ast} + L j^{-1}$ on this subspace.
\end{lemma}

We can extend the Corollary \ref{cor:a0sunstar} for $A_0^{\odot *}$ to $A^{\odot *}$, which will be needed for the computation of normal form coefficients.
\begin{corollary}\label{cor:asunstar}
For $\varphi \in X$ the following statements are equivalent
\begin{enumerate}
\item $j\varphi \in D(A^{\odot *})$ and $A^{\odot *}j\varphi \in Y^{\odot *} \times L^\infty([-h,0];Y^{**})$
\item $j_Y\varphi(0) \in D(B^{\odot *})$ and $\varphi$ has an a.e. derivative $\dot{\varphi} \in L^\infty([-h,0];Y)$
\end{enumerate}
In this case, the action of $A^{\odot *}$ is given by $A^{\odot *}j\varphi = (B^{\odot *} j_Y \varphi(0) + j_Y DG(0)\varphi, \dot{\varphi})$. 
\end{corollary}
\textit{Proof.} The statement on the domain follows immediately from Lemma \ref{lem:3.12} and Corollary \ref{cor:a0sunstar}. Furthermore we have that 
\[A^{\odot*}j \varphi = A_0^{\odot *} j \varphi + \ell DG(0) \varphi = (B^{\odot *} j_Y \varphi(0), \dot{\varphi}) + (j_Y DG(0)\varphi,0)\]\qed

We are now able to state the result which relates $A$ to $A_0$. 
\begin{theorem}\citep[Corollary 23]{janssens_class_2019}
For the generator $A$ of the semigroup $T$ we have that
\begin{equation}\label{eq:generatorA}
\begin{split}
D(A)&=\{\varphi \in X | j\varphi \in D(A_0^{\odot *}), A_0^{\odot *}j\varphi + L \varphi \in j(X)\}\\
A &= j^{-1}(A_0^{\odot *} j + L)
\end{split}
\end{equation}
\end{theorem}

We can cast \eqref{eq:generatorA} in a form which can also be found in the book by \citet[Theorem VI.6.1]{engel_one-parameter_1999} by using Corollary \ref{cor:a0sunstar}.
\begin{theorem}\label{thm:generatorA}
For the generator $A$ of the semigroup $T$ we have that
\begin{equation}\label{eq:generatorA2}
\begin{split}
D(A)&=\{\varphi \in C^1([-h,0];Y) | \varphi(0)\in D(B), \dot{\varphi}(0) = B\varphi(0) + DG(0)\varphi\}\\
A\varphi &= \dot{\varphi}
\end{split}
\end{equation}
\end{theorem}
\textit{Proof.} Let $j\varphi\in D(A_0^{\odot *})$ and $A_0^{\odot *}j\varphi + L \varphi \in j(X)$. As $L\varphi \in j_Y(Y)\times \{0\}$, we have that $A_0^{\odot *}j\varphi \in Y^{\odot *} \times L^{\infty}([-h,0];Y^{**})$. By Corollary \ref{cor:a0sunstar}, $j_Y\varphi(0) \in D(B^{\odot *})$ and $\varphi$ has an a.e. derivative $\dot{\varphi} \in L^\infty([-h,0];Y)$. Furthermore we have that 
\[A_0^{\odot *}j\varphi + L \varphi = (B^{\odot *}j_Y\varphi(0) + j_Y DG(0)\varphi, \dot{\varphi}) \in j(X)\]
By Lemma \ref{lem:embedj} this implies that $B^{\odot *}j_Y\varphi(0) + j_Y DG(0)\varphi \in j_Y(Y), \dot{\varphi} \in C([-h,0];Y)$ and $\dot{\varphi}(0) = B\varphi(0) + DG(0)\varphi$. Hence $\varphi\in C^1([-h,0];Y)$ and $B^{\odot *}j_Y\varphi(0) \in j_Y(Y)$. 

Let $B^{\odot *}j_Y\varphi(0)=j_Y y$ with $y\in Y$. As $B^{\odot *}j_Y\varphi(0) \in Y^{\odot \odot}$, $j_Y\varphi(0)\in D(B^{\odot \odot})$. Let $S^{\odot \odot}$ be the strongly continuous semigroup generated by $B^{\odot \odot}$. This implies that 
\[j_Y\frac{1}{t}(S(t)\varphi(0)-\varphi(0))=\frac{1}{t}(S^{\odot \odot}(t)j_Y\varphi(0)-j_Y\varphi(0))\]
for all $t>0$ \citep[Appendix II Proposition 3.17]{diekmann_delay_1995} . By continuity of $j_Y^{-1}$, this converges in norm as $t\downarrow 0$ to $j_Y B\varphi(0) = B^{\odot \odot} j_Y\varphi(0)$ with $\varphi(0)\in D(B)$.

Conversely, let $\varphi \in C^1([-h,0];Y),  \varphi(0)\in D(B)$ and $\dot{\varphi}(0) = B\varphi(0) + DG(0)\varphi$. Furthermore, let $y^\odot \in D(B^\odot)$, then
\[\langle j_Y B\varphi(0), y^\odot \rangle = \langle y^\odot, B\varphi(0) \rangle = \langle B^\odot y^\odot, \varphi(0) \rangle = \langle j_Y \varphi(0), B^\odot y^\odot \rangle \]
Hence $j_Y \varphi(0)\in D(B^{\odot *})$ and, by Corollary \ref{cor:a0sunstar}, $j\varphi \in D(A_0^{\odot *})$. Furthermore 
\begin{align*}
A_0^{\odot *}j\varphi + L \varphi &=(B^{\odot *}j_Y\varphi(0) + j_Y DG(0)\varphi, \dot{\varphi})\\
&= (j_Y B\varphi(0) + j_Y DG(0)\varphi, \dot{\varphi}) = j \dot{\varphi}\in j(X)
\end{align*}

Finally, for the action of $A$ we derive 
\[A\varphi = j^{-1}(A_0^{\odot *} j + L)\varphi = j^{-1}(j_Y B\varphi(0) + j_Y DG(0)\varphi, \dot{\varphi}) =\dot{\varphi}\] \qed

\subsection{Spectral properties}
In this section we state some results on the spectrum of the operator $A$, notably its essential spectrum and a method for computing its eigenvalues.

For an operator $A$ on $X$ the resolvent set $\rho(A)$ is the set of all $z\in \mathbb{C}$ such that the operator $z-A$ has a bounded inverse. The resolvent operator $R(z,A): X\rightarrow D(A)$ is then defined as $R(z,A)= (z-A)^{-1}$ for $z \in \rho(A)$. The spectrum of $A$, $\sigma(A)=\mathbb{C}\setminus \rho(A)$, can be decomposed into the point spectrum $\sigma_p(A)$ and the essential spectrum $\sigma_{ess}(A)$. We use Weyl's definition of the essential spectrum, i.e. $\sigma_{ess}(A):= \{\lambda\in \mathbb{C} | \lambda -A$ is not a Fredholm operator$\}$ \citep{kato_perturbation_1995}. Then $\sigma_P(A)=\sigma(A) \setminus\sigma_{ess}(A)$ is the discrete spectrum, i.e. isolated eigenvalues with a finite dimensional eigenspace.
\begin{lemma}\label{thm:spectrum_a0}
For the respective spectra we have $\sigma(A_0)=\sigma(A_0^*)=\sigma(A_0^\odot)=\sigma(A_0^{\odot *})=\sigma(B)$. Furthermore $\sigma_{ess}(A_0)=\sigma_{ess}(B)$.
\end{lemma}
\textit{Proof.}
We have that $\sigma(A_0)=\sigma(A_0^*)=\sigma(A_0^\odot)=\sigma(A_0^{\odot *})$ \citep[Proposition IV.2.18]{engel_one-parameter_1999}. 

Next we consider the eigenvalues of $A_0$. For some $\lambda \in \sigma(A_0)$, we need to find a $\varphi\in D(A_0)$ such that $\dot{\varphi}=\lambda \varphi$. Clearly, this is the case if and only if $\varphi(\theta)= q e^{\lambda \theta}$ for $\theta \in [-h,0]$, with $q\in D(B)$ and $Bq = B \varphi(0) = \dot{\varphi}(0) = \lambda q$. Therefore $\lambda \in \sigma_p(A_0)$ if and only if $\lambda \in \sigma_p(B)$ as the corresponding eigenspaces have the same dimension. 

Finally we show that $\rho(A_0)=\rho(B)$, which completes the proof. If $z\in \rho(B)$ then we can find the resolvent of $A_0$ explicitly as for all $\varphi \in X$ and $\theta \in [-h,0]$,  \citep[Proposition VI.6.7]{engel_one-parameter_1999}
\begin{equation}
[R(z,A_0))\varphi](\theta)= e^{z \theta}R(z,B)\varphi(0) + \int_{\theta}^0 e^{z(\theta-s)}\varphi(s)\,ds
\end{equation}
Hence $z \in \rho(A_0)$.

Conversely, suppose that $z \in \rho(A_0)$ en let $y\in Y$. Then the constant function $\psi(\theta) :=y$ for $\theta \in [-h,0]$ is in X and hence $\varphi := R(z,A_0)\psi \in D(A_0)$. This implies that $\varphi(0) \in D(B)$ and $(z-B)\varphi(0) =  z \varphi(0) - \dot{\varphi}(0)  = ((z-A_0)\varphi)(0) = \psi(0) = y$. Hence $z-B$ is surjective. As $z$ is not an eigenvalue of $A_0$, by the above reasoning it is not an eigenvalue of $B$ and hence $z-B$ is injective.

So we conclude that $\sigma(A)= \sigma(B)$ and $\sigma_{ess}(A_0)=\sigma_{ess}(B)$. \qed

If $DG(0)$ is compact then we can make inferences on the essential spectrum of $A$ from the spectrum of $A_0$.
\begin{theorem}\label{thm:ess_spectrum}
If $DG(0)$ is compact then $\sigma_{ess}(A) = \sigma_{ess}(B)$.
\end{theorem}
\textit{Proof.} We will proof this by working in the dual space. This is possible as $\sigma_{ess}(A) = \sigma_{ess}(A^*)$, which is a consequence of the properties of Fredholm operators. \citep[Theorem IV.5.14]{kato_perturbation_1995}

On $X^*$, $A^*=A_0^*+L^*$ due to Lemma \ref{lem:A*}. As $\ell$ is bounded, $L= \ell DG(0)$ is compact and so is its adjoint $L^*$ due to Schauder's theorem, \citep[Theorem III.4.10]{kato_perturbation_1995}. Hence $A^*$ is a compact perturbation of $A_0^*$. One of the defining properties of Weyl's essential spectrum is that it is invariant under compact perturbations \citep[Theorem IV.5.35]{kato_perturbation_1995}.

So we conclude that 
\[\sigma_{ess}(A) = \sigma_{ess}(A^*)=\sigma_{ess}(A_0^*)=\sigma_{ess}(A_0)=\sigma_{ess}(B) \] \qed 

In case $B$ is the diffusion operator, its essential spectrum is empty, see Lemma \ref{lem:spectrumB}. This means that also the essential spectrum of $A$ is empty, when $DG(0)$ is compact.

For computation the eigenvalues we follow \citet{engel_one-parameter_1999}. We introduce the family of operators $K^z:Y\rightarrow Y$, $H^z:X \rightarrow X$ and $W^z: X \rightarrow Y$ parametrized by $z\in \mathbb{C}$, defined as
\begin{equation}
\begin{split}
K^z y &:= DG(0)(y e^{z\theta})\\
(H^z\varphi)(\theta) &:= \int_\theta^0 e^{z(\theta-s)}\varphi(\theta)\,ds\\
W^z\varphi &:= \varphi(0)+DG(0)H^z\varphi
\end{split}
\end{equation}
for $y \in Y, \varphi \in X$ and $\theta \in [-h,0]$. Using these we can define the characteristic operator $\Delta(z)$
\begin{equation}
\Delta(z) = z-B-K^z
\end{equation}

Now we formulate the main theorem of this section which allows us to reduce the computation of the eigenvalues and eigenvectors in $X$ to a computation on $Y$.
\begin{theorem}
\citep[Proposition VI.6.7]{engel_one-parameter_1999} For every $z \in \mathbb{C}$, $\varphi \in \mathcal{R}(z-A)$ if and only if 
\[\Delta(z)q = W^z \varphi\]
has a solution $q\in D(B)$. Moreover $z\in \rho(A)$ if and only if this $q$ is unique. In that case the resolvent is given by
\[(R(z,A)\psi)(\theta) = e^{z\theta}\Delta^{-1}(z)W^z\varphi+(H^z\psi)(\theta)\]
where $\theta\in [-h,0]$ and $\psi \in X$. Finally, $\psi \in D(A)$ is an eigenvector corresponding to $\lambda \in \sigma_p(A)$ if and only if $\psi(\theta)=e^{\lambda \theta}q$, where $q\in D(B)$ is non-trivial and satisfies \[\Delta(\lambda)q=0\]
\label{thm:spectraltheorem}
\end{theorem}

\subsection{Hopf bifurcation}
We are interested in the nonlinear behaviour of \eqref{ADDE}. In this section we develop techniques to compute the first Lyapunov coefficient for (Andronov-)Hopf bifurcations. These techniques can be extended to other local bifurcations, but we will not address those here. In this section, we will follow the methods from \citet{van_gils_local_2013}.

Suppose that $\sigma(A)$ contains a pair of simple purely imaginary eigenvalues $\lambda = \pm i\omega$ with $\omega>0$ and no other eigenvalues on the imaginary axis. Let $\psi \in X$ be the corresponding eigenvector of $A$ and $\psi^{\odot} \in X^\odot$ be the corresponding eigenvector of $A^\odot$ respectively,
\begin{equation}
A\psi= i \omega \psi,\qquad A^\odot\psi^\odot= i \omega \psi^\odot
\end{equation}
We normalise these vectors such that 
\begin{equation}
\langle\psi^\odot,\psi\rangle=1
\end{equation}
The center subspace $X_0$ is spanned by the basis $\Psi=\{\psi,\bar{\psi}\}$ of eigenvectors corresponding to the critical eigenvalues of $A$. Here $\bar{\psi}$ denotes the complex conjugate of $\psi$. 

In order to extend this to the non-linear setting, we need a (locally) invariant critical center manifold $W^c_{loc}$, which is tangent to $X_0$ at the equilibrium at the origin. From \citet{janssens_class_2020}, we get a general result on the existence of this center manifold. 
\begin{theorem}\citep[Theorem 41]{janssens_class_2020}
If the strongly continuous semi-group $S$ generated by $B$ is immediately norm continuous,  $X_0$ is finite-dimensional,  $\sigma(A)$ is the pairwise disjoint union of the sets
\begin{align*}
    \sigma_- &:= \{\lambda \in \sigma(A)| \mathrm{Re }\, \lambda <0\}\\
    \sigma_0 &:= \{\lambda \in \sigma(A)| \mathrm{Re }\, \lambda =0\}\\
    \sigma_+ &:= \{\lambda \in \sigma(A)| \mathrm{Re }\, \lambda >0\}
\end{align*}
where $\sigma_-$ is closed and both $\sigma_0,\sigma_+$ are compact, and if 
\[\sup_{\lambda\in \sigma_-} \mathrm{Re }\, \lambda < 0 < \inf_{\lambda\in\sigma_+} \mathrm{Re }\, \lambda\] 
then there exist a $C^k$-smooth mapping $\mathcal{C}:X_0 \rightarrow X$ and an open neighbourhood $U$ of the origin in $X_0$ such that $\mathcal{C}(0)=0, D\mathcal{C}(0)= I_{X_0 \rightarrow X}$, the identity mapping, and $\mathcal{W}^c_{loc}=\mathcal{C}(U)$ is locally positively invariant for \eqref{ADDE} and contains every solution of \eqref{AIE} that exists on $\mathbb{R}$ and remains sufficiently small for all time. 
\end{theorem}
The conditions on $\sigma(A)$ can be easily satisfied when $\sigma_0$ and $\sigma_{+}$ are composed of finitely many eigenvalues of finite multiplicity. In case $B$ is the diffusion operator, it is immediately norm continuous by Lemma \ref{lem:prop_b} and the essential spectrum $\sigma_{ess}(A)=\sigma_{ess}(B)=\emptyset$ by Theorem \ref{thm:ess_spectrum} and Lemma \ref{lem:spectrumB}. Also when $B=-\alpha I$, $\alpha>0$, we get that the conditions are likewise satisfied.

If $\zeta \in X_0$ then we can write $\zeta = z \psi + \bar{z}\bar{\psi}$ for some $z\in \mathbb{C}$. Using this we can recast $\mathcal{C}(U)$ into the formal expansion $\mathcal{H}:\mathbb{C} \rightarrow W^c_{loc}$
\begin{equation}\label{eq:centermanifold}
\mathcal{H}(z,\bar{z})= z \psi + \bar{z}\bar{\psi} + \sum_{j+k \geq 2} \frac{1}{j! k!} h_{jk}z^j\bar{z}^k
\end{equation}
Due to Theorem \ref{thm:2.15} the \eqref{ADDE} and \eqref{AIE} formulations are equivalent. By weak$^*$ differentiation of \eqref{AIE} and exploiting the finite dimensionality of $\mathcal{W}^c_{loc}$, one can show that a solution $v\in C(\mathbb{R}^+;X)$, $v(t)=u_t$, of \eqref{AIE} satisfies the abstract ODE
\begin{equation}\label{eq:abstractode}
\dot{v}(t) =j^{-1}(A^{\odot *} j v(t) + \ell R(v(t)))
\end{equation}
Where the nonlinearity $R: X\rightarrow Y$ is given by 
\begin{equation}\label{eq:nonlinearity}
R(\varphi):=G(\varphi)-DG(0)(\varphi)=\frac{1}{2}D^2G(0)(\varphi,\varphi)+\frac{1}{6}D^3G(0)(\varphi,\varphi,\varphi)+\mathcal{O}(\|\varphi\|^4)
\end{equation}

Let $\zeta(t) = z(t) \psi + \bar{z}(t)\bar{\psi}$ be the projection of $v(t)$ onto the center subspace $X_0$. The function $z(t)$ satisfies a complex ODE which is smoothly equivalent to the Poincaré normal form 
\begin{equation}\label{eq:normal_form}
\dot{z}= i \omega z + c_1 z|z|^2 + \mathcal{O}(|z|^4)
\end{equation}
where $z,c_1 \in \mathbb{C}$. In polar coordinates, $z=r e^{i \theta}$, this is orbitally equivalent to
\begin{equation}\label{eq:normal_form_real}
\begin{cases}
\dot{r}&= l_1 r^3 + \mathcal{O}(|r|^4)\\
\dot{\theta} &= 1 + \mathcal{O}(|r|^2)
\end{cases}
\end{equation}
where $l_1$ is the {\it first Lyapunov coefficient} determined by the formula
\begin{equation}
l_1=\frac{1}{\omega}\mathrm{Re}(c_1)
\end{equation}
It is well known \citep{kuznetsov_elements_2004} that in generic unfoldings of \eqref{eq:normal_form_real}, $l_1<0$ implies that the bifurcation is supercritical and that a stable limit cycle exists near one of the branches. On the other hand, $l_1>0$ implies that the bifurcation is subcritical and that an unstable limit cycle exists near one of the branches.

The critical center manifold $\mathcal{W}^c_{loc}$ has the expansion \eqref{eq:centermanifold} and due to the time-invariance of $\mathcal{W}^c_{loc}$ we have
\begin{equation}
v(t)=\mathcal{H}(z(t),\bar{z}(t))
\end{equation}
If we differentiate both sides with respect to time and use the abstract ODE \eqref{eq:abstractode} for the left-hand side, we obtain the {\it homological equation}
\begin{equation}\label{eq:homological_equation}
A^{\odot *} j \mathcal{H}(z,\bar{z}) + \ell R(\mathcal{H}(z,\bar{z})) = j \mathcal{H}_z(z,\bar{z})\dot{z} + j \mathcal{H}_{\bar{z}}(z,\bar{z})\dot{\bar{z}}
\end{equation}
We can substitute the expansion of the nonlinearity \eqref{eq:nonlinearity}, the normal form \eqref{eq:normal_form}, and the expansion of the critical center manifold \eqref{eq:centermanifold} into the homological equation \eqref{eq:homological_equation} to derive the normal form coefficients. If we equate coefficients of the corresponding powers of $z$ and $\bar{z}$ we obtain the following equations
\begin{equation}\label{eq:Hopf_cond}
\begin{split}
(2 i \omega - A^{\odot *})j h_{20}=&\; \ell D^2G(0)(\psi,\psi)\\
-A^{\odot *}j h_{11}=&\; \ell D^2G(0)(\psi,\bar{\psi})\\
(i \omega - A^{\odot *})j h_{21}=&\; \ell D^3G(0)(\psi,\psi,\bar{\psi})+\ell D^2G(0)(h_{20},\bar{\psi})\\
&+ 2\ell D^2G(0)(\psi,h_{11})-2c_1j\psi
\end{split}
\end{equation}
They all have the form
\begin{equation}\label{eq:asunstarresolvent}
(z - A^{\odot *})\varphi^{\odot *}=\psi^{\odot *}
\end{equation}
Here $z \in \mathbb{C}$ and $\psi^{\odot *}\in X^{\odot *}$ are given. When $z \in \rho(A)$ then \eqref{eq:asunstarresolvent} has a unique solution. However, if $z \in \sigma(A)$, then a solution $\varphi^{\odot *}$ doesn't necessarily exist for all $\psi^{\odot *}$. The following lemma, which is equivalent to \citep[Lemma 33]{van_gils_local_2013}, provides a condition for solvability.\\
\begin{lemma}[Fredholm solvability]\label{lem:Fredholmsolvability}
Let $z \notin \sigma_{ess}(A)$. Then $z - A^\odot:D(A^\odot)\rightarrow X^\odot$ has closed range. In particular $(z -A^{\odot *})\varphi^{\odot *}=\psi^{\odot *}$ is solvable for $\varphi^{\odot *} \in D(A^{\odot *})$ given $\psi\in X^{\odot *}$ if and only if $\langle \psi^{\odot *},\psi^{\odot}\rangle=0$ for all $\psi^\odot \in \mathcal{N}(z-A^\odot)$. 
\end{lemma}
\textit{Proof.} From the definition of the essential spectrum, $\mathcal{R}(z -A)$ is closed \citep[Section IV.5.1]{kato_perturbation_1995}, and $\mathcal{R}(z-A^*)$ is also closed by Banach's Closed Range Theorem \citep[Theorem IV.5.13]{kato_perturbation_1995}. Let $(\psi_n^\odot)_{n\in \mathbb{N}}$ be a sequence in $\mathcal{R}(z-A^\odot)$ such that $\psi_n^\odot \rightarrow \psi^\odot \in X^\odot$. Then there is a sequence $(\varphi_n^\odot)_{n\in \mathbb{N}}$ in $D(A^\odot)$ such that
\[\psi_n^\odot = (z-A^\odot)\varphi_n^\odot = (z-A^*)\varphi_n^\odot \qquad \forall n \in \mathbb{N}\]
Hence $\psi_n^\odot \in \mathcal{R}(z-A^*)$ for all $n\in \mathbb{N}$, so there exists $\varphi^\odot \in D(A^*)$ such that $(z-A^*)\varphi^\odot = \psi^\odot$ and
\[A^*\varphi^\odot = z \varphi^\odot -(z-A^*)\varphi^\odot   = z \varphi^\odot - \psi^\odot \in X^\odot \]
Hence $\varphi^\odot \in D(A^\odot)$, $(z-A^\odot)\varphi^\odot = \psi^\odot$ and $\psi^\odot \in \mathcal{R}(z-A^\odot)$.

Due to Banach's Closed Range Theorem, $\varphi^{\odot *}$ is a solution of
\[(z-A^{\odot *})\varphi^{\odot *}=\psi^{\odot *}\]
given $\psi^{\odot *}$ if and only if 
\[\langle \psi^{\odot *}, \psi^\odot \rangle =0 \qquad \forall \psi^\odot \in \mathcal{N}(z-A^\odot)\]\qed

We now return to equations \eqref{eq:Hopf_cond}. As $\{0,2i \omega\} \subset \rho(A)=\rho(A^\odot)$ we can use the resolvent of $A^{\odot *}$ to solve the first two equations. However, $i \omega \in \sigma(A)$ so for the last equation of \eqref{eq:Hopf_cond} we need to use the theorem above. The corresponding eigenspace $\mathcal{N}(A^*-\lambda)$ is spanned by $\psi^\odot$, so we can compute for the normal form coefficient by
\begin{equation}
\begin{split}
j h_{20}&= R(2i \omega,A^{\odot *})\ell D^2G(0)(\psi,\psi)\\
j h_{11}&= R(0,A^{\odot *})\ell D^2G(0)(\psi,\bar{\psi})\\
c_1&=\frac{1}{2}\langle \ell D^3G(0)(\psi,\psi,\bar{\psi})+\ell D^2G(0)(h_{20},\bar{\psi})+ 2\ell D^2G(0)(\psi,h_{11}), \psi^\odot\rangle
\end{split}
\label{eq:Hopf_nfc}
\end{equation}
We are not yet able to compute the normal form coefficient explicitly as we don't have an explicit representation of $\psi^\odot$ or a representation of the resolvent of $A^{\odot *}$. However, we resolve this by using spectral projections.

Let $P^\odot$ and $P^{\odot *}$ be the spectral projections on $X^\odot$ and $X^{\odot *}$ corresponding to some eigenvalue $\lambda$, respectively. Then $P^{\odot *}\varphi^{\odot *} = \nu j \psi$ for some $\nu \in \mathbb{C}$ and
\[\langle \varphi^{\odot *}, \psi^\odot\rangle=\langle \varphi^{\odot *}, P^\odot \psi^\odot\rangle=\langle P^{\odot *}\varphi^{\odot *}, \psi^\odot\rangle=\nu \langle j\psi, \psi^\odot\rangle=\nu\]
Hence we seek to determine $\nu$. From the Dunford integral representation it follows that
\begin{equation}\label{eq:dunford_integral}
P^{\odot *}\varphi^{\odot *} = \frac{1}{2\pi i} \oint_{\partial C_\lambda} R(z,A^{\odot *})\varphi^{\odot *}\,dz= \nu j \psi
\end{equation}
where $C_\lambda$ is a sufficiently small open disk centered at $\lambda$ and $\partial C_\lambda$ its boundary. 
The element on the left in the pairing \eqref{eq:Hopf_nfc} is of the form $\varphi^{\odot *} = \ell y$, $y\in Y$. In this case we can reduce $R(z,A^{\odot *})\varphi^{\odot *}$ to $\Delta^{-1}(z)y$ by virtue the following theorem.
\begin{theorem}\label{thm:resolvent_Asunstar}
Suppose that $z \in \rho(A)$. For each $y\in Y$ the function $\varphi\in X$, defined as $\varphi(\theta) :
= e^{z \theta} \Delta^{-1}(z) y$ for $\theta \in [-h,0]$, is the unique solution in $\{\varphi \in C^1([-h,0];Y)| \varphi(0)\in D(B)\}$ of the system
\begin{equation}\label{eq:36}
\begin{cases}
(z-B)\varphi(0) - DG(0)\varphi &= y\\
z \varphi -\dot{\varphi} &= 0
\end{cases}
\end{equation}
Moreover, $\varphi^{\odot *} = j\varphi$ is the unique solution in $D(A^{\odot *})$ of $(z-A^{\odot *})\varphi^{\odot *} = \ell y $.
\end{theorem}
\textit{Proof.} Since $z \in \rho(A)$, by Theorem \ref{thm:spectraltheorem} it follows that $\Delta^{-1}(z)$ exists. We start by showing that $\varphi$ as defined above solves \eqref{eq:36}. Clearly $\varphi \in C^1([-h,0];Y)$ and $\varphi(0)=\Delta^{-1}(z)y\in D(B)$. Recall from the definition of $K^z$ that for $q\in Y$, $K^z q = DG(0)q e^{z \theta}$. Therefore, 
\[(z-B)\varphi(0) - DG(0)\varphi = (z-B)\Delta^{-1}(z)y - K^z \Delta^{-1}(z) y = y\]
Finally, by differentiating $\varphi$ we see that it satisfies the second equation in \eqref{eq:36}. 

When $\varphi(0)\in D(B)$ then $j_Y \varphi(0) \in D(B^{\odot *})$, because for all $y^\odot \in D(B^\odot)$
\[\langle j_Y B\varphi(0), y^\odot \rangle = \langle y^\odot, B\varphi(0) \rangle = \langle B^\odot y^\odot, \varphi(0) \rangle = \langle j_Y \varphi(0), B^\odot y^\odot \rangle \]
 
Then corollary \ref{cor:asunstar} implies that $j \varphi \in D(A^{\odot *})$. 
\[(z-A^{\odot *})\varphi^{\odot *} = (j_Y(z-B)\varphi(0) - j_Y DG(0)\varphi,z \varphi-\dot{\varphi})=(j_Y y,0)= \ell y\]
But by Theorem \ref{thm:spectraltheorem} $\rho(A^{\odot *}) = \rho(A)$, so $\varphi^{\odot *} = j \varphi$ is the unique solution of $(z-A^{\odot *})\varphi^{\odot *} = \ell y$. Consequently, $\varphi$ itself is the unique solution in $\{\varphi \in C^1([-h,0];Y)| \varphi(0)\in D(B)\}$. \qed

Now given that we can compute the resolvent $\Delta^{-1}(z)$ and the Fréchet derivatives of $G$, we have a method to compute the center manifold coefficients $h_{20}$ and $h_{11}$, and the first Lyapunov coefficient $l_1 = \tfrac{1}{\omega}\mathrm{Re}\,c_1$:
\small
\begin{equation}\label{eq:normal_form_computation}
\begin{split}
h_{20}(\theta)&= e^{2i \omega \theta}\Delta^{-1}(2i \omega) D^2G(0)(\psi,\psi)\\
h_{11}(\theta)&= \Delta^{-1}(0) D^2G(0)(\psi,\bar{\psi})\\ 
c_1 \psi(\theta) &= \frac{1}{4\pi i} \oint_{\partial C_\lambda} \!\!\!\!\! e^{z\theta}\Delta^{-1}(z)(D^3G(0)(\psi,\psi,\bar{\psi}) + D^2G(0)(h_{20},\bar{\psi}) + 2 D^2G(0)(\psi,h_{11}))\,dz
\end{split}
\end{equation}
\normalsize

\section{Characterisation of the Spectrum}
\label{sec:single_pop}
In this section we will return to the Neural Field as derived in section \ref{sec:modelling}. For certain choices we can derive some explicit conditions for the spectrum and find an explicit expression for the resolvent. 

We take $Y=C(\Omega)$ with $\Omega = [-1,1]$ and use the \eqref{ADDE1} formulation of the section \ref{sec:duality}
\begin{equation}
\begin{cases}
\dot{u}(t)=Bu(t)+G(u_t)\\
u_0=\varphi\in X
\end{cases}
\label{ADDE1} \tag{ADDE}
\end{equation}
Where $B: D(B)\rightarrow Y$ and $G:X\rightarrow Y$ are defined as
\begin{align*}
B q &:= d q'' - \alpha q\\
D(B) &:= \{ q\in Y | q\in C^2(\Omega), q'(\partial \Omega)=0 \}\\
G(\varphi) &:= \alpha \int_{\Omega} J(x,x')S(\varphi(t-\tau(x,x'),x'))\,dx'
\end{align*}
Here we assume that $d \geq 0, \alpha >0$, $J$ and $\tau$ are continuous functions and $S \in C^{\infty}(\mathbb{R})$, with $S(0)=0$ and $S'(0)\neq 0$. The assumption $S(0)=0$ makes sure we have an equilibrium at $u\equiv 0$. We interpret $u$ as the deviation from this physiological resting state. This interpretation then makes for cleaner notation.

We have the following properties for $G$ and its derivatives.
\begin{lemma}\label{lem:Gderivative}
\citep[Lemma 3, Proposition 11]{van_gils_local_2013} $G$ is compact, globally Lipschitz continuous and $k$ times Fréchet differentiable for any $k\in \mathbb{N}$. Furthermore the $k$th Fréchet derivative of $G$ at $\psi\in X$, $D^kG(\psi): X^k \rightarrow Y$, is compact and given by
\begin{align*}
(D^kG(\psi)(\varphi_1,\cdots, \varphi_k))(x)= \alpha \int_{\Omega} \Big[ J(x,x')S^{(k)}(\psi(-\tau(x,x'),x')) &\\
\prod_{m=1}^k(\varphi_m(-\tau(x,x'),x'))\Big] & \,dx'
\end{align*}
\end{lemma}
As $DG(0)$ is compact we can find, due to Theorem \ref{thm:ess_spectrum} and Lemma \ref{lem:spectrumB}, that the essential spectrum of the linearisation $A$ is given by
\begin{equation}
\sigma_{ess}(A) = \begin{cases}
\emptyset & d>0\\
\{-\alpha \} & d=0 
\end{cases}
\end{equation}

We want to be able to compute the eigenvalues, eigenvectors and resolvent for specific choices of $J$ and $\tau$. We take $J$ as a sum of exponentials and $\tau$ as a constant delay plus a finite propagation speed, which we can normalise to $1$ by scaling time.
\begin{align*}
J(x,x') &:= \sum_{j=1}^N \eta_j e^{-\mu_j |x-x'|}\\
\tau(x,x') &:= \tau^0 + |x-x'|
\end{align*}
Where we take $\tau^0 \geq 0$ and $\eta_j \neq 0$ for $j \in \{1,\cdots,N\}$. 

Due to Theorem \ref{thm:spectraltheorem} we have that $\lambda$ is an eigenvalue and $\psi$ an eigenvector if and only if $\psi(\theta) = q e^{\lambda \theta}$ and $q \in D(B)$ satisfies characteristic equation \eqref{CE}.
\begin{equation}\tag{CE}
\Delta(\lambda)q=(\lambda - B - K^\lambda)q=0
\end{equation}
Where in this case $K^z: Y \rightarrow Y$ is a parametrized family of operators for $z\in \mathbb{C}$ defined as 
\begin{equation}
\begin{split}
K^z &:= \sum_{j=1}^N K^z_j\\
K_j^z\, y(x) &:= c_j(z) \int_{-1}^1 e^{-k_j(z) |x-x'|} y(x')\,dx'
\end{split}
\end{equation}
where $c_j(z) := S'(0) \alpha \eta_j e^{-\tau^0 z}\neq 0$ and $ k_j(z) := \mu_j + z$.

The case without diffusion, i.e. $d=0$, has already been extensively studied \citep{van_gils_local_2013, dijkstra_pitchforkhopf_2015}, so in this section we will develop formula's for the eigenvalues, eigenvectors and resolvent with nontrivial diffusion, i.e. $d>0$.

For the following section we adopt the notational convention that bold-faced variables correspond to vectors $\mathbf{a}=(a_1\; \cdots \; a_n)^T$ where its length is clear from the context.

\subsection{Eigenvalues}
So we are looking for non-trivial solutions $q \in D(B)$ of 
\begin{equation}\label{CE}\tag{CE}
(z-B- K^z)q=0
\end{equation}
As this is a mixed differential-integral equation, it is in general hard to solve. We will use the method of \citet{dijkstra_pitchforkhopf_2015} to convert \eqref{CE} into a differential equation \eqref{ODE}, which we can solve. Then substituting the general solution of \eqref{ODE} back into \eqref{CE} yields appropriate conditions on $q$. This is possible due to the following observations.

\begin{lemma}\label{thm:cinf}
All solutions of \eqref{CE} are $C^{\infty}(\Omega)$. 
\end{lemma}
\textit{Proof.} As $q\in C^2(\Omega)$ and the range of $K^z$ is contained in $C^3(\Omega)$ we have that $B q \in C^2(\Omega)$, which means that $q\in C^4(\Omega)$. By induction, we conclude that $q \in C^{\infty}(\Omega)$. \qed

Differentiating the kernel functions in the \eqref{CE} in the distributional sense yields for $j\in \{1,\cdots, N\}$
\[\frac{\partial^2}{\partial x^2}e^{-k_j(z)|x-x'|}=\left[ k_j^2(z)-2k_j(z)\delta(x-x')\right]e^{-k_j(z)|x-x'|}\]
So we define the differential operator $L_j^z$ for $j \in \{1,\cdots, N\}$. 
\[L_j^z:= k_j^2(z)-\partial_x^2\]
For this operator $L_j$ we have that for $j\in \{1,\cdots, N\}$
\[L_j^z K_j^z q = 2 c_j(z) k_j(z) q\]
Hence by applying the operator $L^z = \prod_{p=1}^N L^z_p$ to \eqref{CE} we end up with an ordinary differential equation \eqref{ODE}
\begin{equation}\label{ODE} \tag{ODE}
L^z\Delta(z)q=(z-B)\prod_{p=1}^N L_p^z q - 2 \sum_{j=1}^N c_j(z) k_j(z) \prod_{\substack{p=1\\ p\neq j}}^N L_p^z q =0
\end{equation}
This differential equation has a characteristic polynomial corresponding to exponential solutions $e^{\rho x}$
\begin{equation}\label{eq:char_poly}
P^z(\rho):=(\alpha+z-d\rho^2)\prod_{p=1}^N (k_{p}(z)^2-\rho^2)-2 \sum_{j=1}^N c_j(z)k_j(z) \prod_{\substack{p=1\\ p\neq j}}^N (k_{p}(z)^2-\rho^2)
\end{equation}
$P^z$ is an even polynomial of order $2(N+1)$. Assuming that $z$ is such that $P^z$ has exactly $2(N+1)$ distinct roots $\pm \rho_1(z),\cdots,\pm \rho_{N+1}(z)$, the general solution $q$ of \eqref{ODE} is a linear combination of exponentials $e^{\pm \rho_j x}$. 
\begin{equation}\label{eq:generalsol}
q(x):=\sum_{m=1}^{N+1}\left[a_m \cosh(\rho_m(z) x)+b_m \sinh(\rho_m(z) x)\right]) 
\end{equation}
Writing $q$ as a linear combination of cosine hyperbolic and sine hyperbolic leads to cleaner notation below. 

Before we substitute \eqref{eq:generalsol} back into \eqref{CE}, we first prove two lemmas.
\begin{lemma}\label{lem:nonzero}
If the characteristic polynomial $P^z(\rho)$ has $2(N+1)$ distinct roots then $\rho_m(z)\neq 0$ for all $m\in \{1,\cdots, N+1\}$ and $k_j(z)\neq 0$ for all $j \in \{1,\cdots, N\}$.
\end{lemma}
\textit{Proof.} If $P^z(\rho)$ has $2(N+1)$ distinct roots $\pm \rho_1(z), \cdots, \pm \rho_{N+1}(z)$, then $\rho_m(z)$ is distinct from $- \rho_m(z)$ and hence $\rho_m(z) \neq 0$ for $m\in \{1,\cdots, N+1\}$. 

Let without loss of generality $k_1(z)=0$. In that case the characteristic polynomial becomes
\[P^z(\rho)=\rho^2(\alpha+z-d\rho^2) \prod_{p=2}^N (k_{p}(z)^2-\rho^2)-2 \rho^2 \sum_{j=2}^N c_j(z)k_j(z) \prod_{\substack{p=2\\ p\neq j}}^N (k_{p}(z)^2-\rho^2)\]
So $\rho=0$ is a root of $P^z$. Hence by contradiction we conclude by contradiction that $k_j(z)\neq 0$ for all $j \in \{1,\cdots, N\}$. \qed

Define the set $\mathcal{L}$ as follows
\begin{equation}\label{eq:exceptionL}
\mathcal{L}:=\{z\in\mathbb{C}| \exists j \in \{1,\cdots N\}, m\in \{1,\cdots, N+1\} \text{ such that } k_j(z)= \pm \rho_m(z) \}
\end{equation}
\begin{lemma}\label{lem:nonequal}
If characteristic polynomial $P^z$ has $2(N+1)$ distinct roots then 
\[\mathcal{L}=\{z\in\mathbb{C}| \exists j,p \in \{1,\cdots N\}, j\neq p\text{ such that } k_j^2(z)= k_p^2(z) \}\]
\end{lemma}
\textit{Proof.}
We have that $z\in \mathcal{L}$ if and only if $P^z(k_j(z)) = 0$ for some $j\in \{1,\cdots, N\}$. 
\[P^z(k_j(z))=-2c_j(z)k_j(z)\prod_{\substack{p=1\\ p\neq j}}^N (k_{p}^2(z)-k_j^2(z))\]
Hence $P^z(k_j(z)) = 0$ if and only if $k_j^2(z)= k_p^2(z)$ for some $p \in \{1,\cdots N\}, j\neq p$ \qed

For $z\notin \mathcal{L}$ we can rewrite $P^z(\rho_m)$ as 
\[P^z(\rho_m) = \left[\alpha + z - d \rho_m^2 - \sum_{j=1}^N\frac{2c_j(z)k_j(z)}{k_j^2(z)-\rho_m^2(z)}\right]\prod_{p=1}^N(k_p^2(z)-\rho_m^2(z))=0\]
We can divide out the product to conclude that for $m\in\{1,\cdots,N+1\}$ and $j\in \{1,\cdots, N\}$
\begin{equation}\label{eq:char_eq}
 \alpha + z - d \rho_m^2 - \sum_{j=1}^N\frac{2c_j(z)k_j(z)}{k_j^2(z)-\rho_m^2(z)}=0
\end{equation}

Next we find formula's for $K_j^z\cosh(\rho_m(z)x)$ and $K_j^z\sinh(\rho_m(z)x)$. To compute these integrals we split the interval $[-1,1]$ into the intervals $[-1,x]$ and $[x,1]$. On these intervals $e^{-k|x-x'|}$ is an $C^1$ function in $x'$ so we can compute the following anti-derivatives for these smooth branches.
\small
\begin{equation} \label{eq:antid}
\begin{split}
\int^{x'} e^{-k|x-s|} \cosh(\rho s) \,ds &= \begin{cases}
e^{-k|x-x'|}\frac{(\;\; k \cosh(\rho x') - \rho \sinh(\rho x'))}{k^2-\rho^2} + const. & -1\leq x' < x \leq 1\\
e^{-k|x-x'|}\frac{(-k \cosh(\rho x') - \rho \sinh(\rho x'))}{k^2-\rho^2} + const. & -1\leq x < x' \leq 1\\
\end{cases}\\
\int^{x'} e^{-k|x-s|} \sinh(\rho s) \,ds &= \begin{cases}
e^{-k|x-x'|}\frac{(\;\; k \sinh(\rho x') - \rho \cosh(\rho x'))}{k^2-\rho^2} + const. & -1\leq x' < x \leq 1\\
e^{-k|x-x'|}\frac{(-k \sinh(\rho x') - \rho \cosh(\rho x'))}{k^2-\rho^2} + const. & -1\leq x < x' \leq 1\\
\end{cases}
\end{split}
\end{equation}
\normalsize
Using these anti-derivatives, we can evaluate the integrals $K_j^z\cosh(\rho_m(z)x)$ and $K_j^z\sinh(\rho_m(z)x)$. For clarity we omit the dependence on $z$ in the remainder of this section.
\small
\begin{align*}
K_j\cosh(\rho_m x)&=\frac{2 c_j k_j\cosh(\rho_m x)-2 c_j e^{-k_j}\cosh(k_j x)(k_j\cosh(\rho_m)+\rho_m\sinh(\rho_m))}{k_j^2-\rho_m^2}\\
K_j\sinh(\rho_m x)&=\frac{2 c_j k_j\sinh(\rho_m x)-2 c_j e^{-k_j}\sinh(k_j x)(\rho_m\cosh(\rho_m)+k_j\sinh(\rho_m))}{k_j^2-\rho_m^2}
\end{align*}
\normalsize

Now we are ready to substitute the general solution $q$ of \eqref{ODE}, \eqref{eq:generalsol}, back into \eqref{CE}.

\begin{equation}\label{eq:subs}
\begin{split}
&\sum_{m=1}^N\left[a_m\cosh(\rho_m x)+b_m\sinh(\rho_m x)\right]\left[(\alpha+z+d\rho_m^2)+\sum_{j=1}^N\frac{2 c_j k_j}{k_j^2-\rho_m^2}\right]+\\
&\sum_{j=1}^N c_j e^{-k_j} \left[-\cosh(k_j x)\sum_{m=1}^{N+1} a_m \frac{k_j\cosh(\rho_m)+\rho_m\sinh(\rho_m)}{k_j^2-\rho_m^2}\right. \\
&\qquad\qquad\;\; \left. -\sinh(k_j x)\sum_{m=1}^{N+1} b_m \frac{\rho_m\cosh(\rho_m)+k_j\sinh(\rho_m)}{k_j^2-\rho_m^2}\right]=0
\end{split}
\end{equation}

Due to the characteristic equation \eqref{eq:char_eq} the first line in equation \eqref{eq:subs} vanishes. When $z\notin \mathcal{L}$, $\cosh(k_j x)$ and $\sinh(k_j x)$ for $j \in \{1,\cdots, N\}$ are linearly independent. Hence the second line vanishes if and only if $S^{z,even}\mathbf{a}=S^{z,odd}\mathbf{b}=\mathbf{0}$, where matrices $S^{z,even}$ and $S^{z,odd}$ are defined as
\begin{equation}
\begin{split}
S^{z,even}_{j,m}&:=\frac{k_j\cosh(\rho_m)+\rho_m\sinh(\rho_m)}{k_j^2-\rho_m^2}\\
S^{z,odd}_{j,m}&:=\frac{\rho_m\cosh(\rho_m)+k_j\sinh(\rho_m)}{k_j^2-\rho_m^2}
\end{split}
\end{equation}
for $j\in \{1,\cdots, N\}$ and $m \in \{1, \cdots, N+1\}$. 

As $q\in D(B)$, we also need to take the boundary conditions into account as 
\begin{equation}\label{eq:bound_eq}
q'(\pm 1) = \sum_{m=1}^{N+1}\left[b_m \rho_m \cosh(\rho_m)\pm a_m \rho_m \sinh(\rho_m)\right]=0
\end{equation}
To satisfy the boundary conditions, we augment the matrices $S^{z,even}$ and $S^{z,odd}$ as follows:
\begin{equation}
\begin{split}
S^{z,even}_{N+1,m}&:=\rho_m \sinh(\rho_m)\\
S^{z,odd}_{N+1,m}&:=\rho_m \cosh(\rho_m)
\end{split}
\end{equation}
Now we have square matrices $S^{z,even}, S^{z,odd} \in \mathbb{C}^{(N+1)\times (N+1)}$. There exists a non-trivial solution $q\in D(B)$ of \eqref{CE} if and only if $\det(S^{z,even})=0$ or $\det(S^{z,odd})=0$.

\begin{theorem} Suppose $\det(P^\lambda(\rho))$ has $2(N+1)$ distinct roots and $\lambda \notin \mathcal{L}$ for some $\lambda \in \mathbb{C}$ then we have that $\lambda\in \sigma_p(A)$ if and only if $\det(S^{\lambda,even})\det(S^{\lambda,odd})=0$.

When $\det(S^{\lambda,even})=0$, the corresponding eigenvector $\psi \in X$ is given by
\begin{equation}
\psi(\theta)(x):=e^{\lambda\theta}\sum_{m=1}^{N+1}a_m \cosh(\rho_m(\lambda) x)
\end{equation}
Where $\mathbf{a}$ is a vector in the nullspace of $S^{\lambda,even}$. 

When $\det(S^{\lambda,odd})=0$, the corresponding eigenvector $\psi \in X$ is given by
\begin{equation}
\psi(\theta)(x):=e^{\lambda\theta}\sum_{m=1}^{N(N+1)}b_m \sinh(\rho_m(\lambda) x)
\end{equation}
Where $\mathbf{b}$ is a vector in the nullspace of $S^{\lambda,odd}$. 
\label{thm:eigenvalues}
\end{theorem}
\textit{Proof.}
Let $q\in D(B)$ be a solution of \eqref{CE} for some $\lambda \in \mathbb{C}$. Then by Theorem \ref{thm:cinf} $q\in C^{\infty}$ so it is also a solution of \eqref{ODE}. 

Conversely, let $q$ be a solution of \eqref{ODE}. As $\det(P^\lambda(\rho))$ has $2(N+1)$ distinct roots, $q$ is of the form \eqref{eq:generalsol}. Due to \eqref{eq:subs} and \eqref{eq:bound_eq} it is a solution of \eqref{CE} if and only if $\det(S^{\lambda,even})\det(S^{\lambda,odd})=0$. \qed

We will call a eigenvalue 'even', respectively 'odd', when $\det(S^{\lambda,even})=0$, respectively $\det(S^{\lambda,odd})=0$.

\subsection{Resolvent}
Due to Theorem \ref{thm:resolvent_Asunstar}, to compute the normal form coefficients we need a representation of $\Delta^{-1}(z)y$. It is defined for $z\in \rho(A)$ as the unique solution $q \in D(B)$ of the resolvent equation \eqref{RE} 
\begin{equation}\label{RE}\tag{RE}
\Delta(z)q= (z-B-K^z)q =y
\end{equation}
We can find an explicit form for this resolvent using a variation-of-constants ansatz when $z\notin \mathcal{S}$, which is defined as
\begin{equation}\label{eq:expceptionS}
\mathcal{S}:=\sigma(B)\cup \mathcal{L}\cup\{z \in \mathbb{C}|P^z(\rho)\text{ has less than } 2(N+1)\text{ distinct zeros}\}
\end{equation}
With $\mathcal{L}$ as in \eqref{eq:exceptionL}. 

\begin{theorem}\label{thm:resolvent}
For $z\in \rho(A)$ with $z\notin \mathcal{S}$ the unique solution $q\in D(B)$ of \eqref{RE} is given by
\begin{equation}\label{eq:res_ansatz}
q(x):=R(z,B)y(x)+\sum_{m=1}^{N+1}\left[a_m(x) \cosh(\rho_m(z) x)+b_m(x) \sinh(\rho_m(z) x)\right]
\end{equation}
Where $R(z,B)$ is the resolvent operator of $B$ as in \eqref{eq:res_b} and $\mathbf{a}(x)$ and $\mathbf{b}(x)$ as in \eqref{eq:res_sol}
\end{theorem}
\textit{Proof.}
Our variation-of-constants Ansatz $q$ needs to satisfy 3 conditions. It must solve \eqref{RE}, $\Delta(z)q=y$, it must satisfy the boundary conditions $(q)'(\pm 1)=0$ and the regularity condition $q \in C^2(\Omega)$. When we a found some $a_m(x),b_m(x)$ such that $q$ satisfies these conditions, we have found the resolvent as it is unique due to Theorem \ref{thm:spectraltheorem}. As $R(z,B)$ maps into $D(B)$, the regularity condition is satisfied when $\mathbf{a}(x), \mathbf{b}(x) \in C^2(\Omega)$. For this proof we suppress the dependencies on $z$.

To aid in the calculation of $\Delta(z)q$, we first compute some integrals up front. We can integrate by parts by splitting the interval $[-1,1]$ into $[-1,x)$ and $(x,1]$ and using the anti-derivatives in \eqref{eq:antid} to end up with
\small
\begin{equation}\label{eq:antid2}
\begin{split}
K_j a_m(x)\cosh &(\rho_m x)= a_m(x)\cosh(\rho_m x)\frac{2 c_j k_j }{k_j ^2-\rho_m^2}\\
&+c_j e^{-k_j (1+x)}a_m(-1)S^{z,even}_{j,m}+c_j e^{-k_j (1-x)}a_m(1)S^{z,even}_{j,m}\\
&-c_j \int_{-1}^1 \frac{a_m'(x')}{k_j ^2-\rho_m^2}e^{-k_j |x-x'|}\left(\mathrm{sgn}(x-x')k_j \cosh(\rho_m x') - \rho_m \sinh(\rho_m x')\right)\,dx'\\
K_j b_m(x)\sinh &(\rho_m x)= b_m(x)\sinh(\rho_m x)\frac{2 c_j k_j }{k_j ^2-\rho_m^2}\\
&-c_j e^{-k_j (1+x)}b_m(-1)S^{z,odd}_{j,m}+c_j e^{-k_j (1-x)}b_m(1)S^{z,odd}_{j,m}\\
&-c_j \int_{-1}^1 \frac{b_m'(x')}{k_j ^2-\rho_m^2}e^{-k_j |x-x'|}\left(\mathrm{sgn}(x-x')k_j \sinh(\rho_m x') - \rho_m \cosh(\rho_m x')\right)\,dx'
\end{split}
\end{equation}
\normalsize
Now we substitute ansatz \eqref{eq:res_ansatz} into \eqref{RE} and collect the terms. Using the above calculations and the fact that $(z-B)R(z,B)y=y$, we have that
\small
\begin{subequations}
\begin{align}
0=&\sum_{m=1}^{N+1}\left[a_m(x)\cosh(\rho_m x)+b_m(x)\sinh(\rho_m x)\right]\left[(\alpha+z-d \rho_m^2(z))-\sum_{j=1}^N\frac{2 c_j k_j }{k_j ^2-\rho_m^2}\right]\label{eq:res2}\\
-&\sum_{m=1}^{N+1}d \left[(a_m''(x)+2\rho_m b_m'(x))\cosh(\rho_m x)+(b_m''(x)+2\rho_m a_m'(x))\sinh(\rho_m x)\right]\label{eq:res1} \\
-&\sum_{j=1}^N c_j e^{-k_j (1+x)}\left[\sum_{m=1}^{N+1}a_m(-1)S^{z,even}_{j,m}-\sum_{m=1}^{N+1}b_m(-1)S^{z,odd}_{j,m}\right]\label{eq:res3}\\
-&\sum_{j=1}^N c_j e^{-k_j (1-x)}\left[\sum_{m=1}^{N+1}a_m(1)S^{z,even}_{j,m}+\sum_{m=1}^{N+1}b_m(1)S^{z,odd}_{j,m}\right]\label{eq:res4}\\
-&\sum_{j=1}^Nc_j \int_{-1}^1e^{-k_j |x-x'|}\left. \bigg[R(z,B)y(x')\right.\nonumber\\
&\qquad\qquad\qquad-\sum_{m=1}^{N+1}\frac{a_m'(x')}{k_j ^2-\rho_m^2}\left(\mathrm{sgn}(x-x')k_j \cosh(\rho_m x') - \rho_m \sinh(\rho_m x')\right)\nonumber\\
&\qquad\qquad\qquad\left.-\sum_{m=1}^{N+1}\frac{b_m'(x')}{k_j ^2-\rho_m^2}\left(\mathrm{sgn}(x-x')k_j \sinh(\rho_m x') - \rho_m \cosh(\rho_m x')\right)\right]\,dx' \label{eq:res5}
\end{align}
\end{subequations}
\normalsize

We have that the above equation vanishes when all the terms within square brackets vanish. The term \eqref{eq:res2} vanishes naturally due to characteristic equation in \eqref{eq:char_eq} as $z\notin \mathcal{L}$. 

As $R(z,B)$ maps into $D(B)$, the boundary conditions $q'(\pm 1)=0$ reduces to
\begin{equation}\label{eq:res6}
\sum_{m=1}^{N+1}\left[(a_m'(\pm 1)+\rho_m b_m(\pm 1)) \cosh(\rho_m)\pm (b_m'(\pm 1)+\rho_m a_m(\pm 1)) \sinh(\rho_m)\right]=0
\end{equation}
We can split equation \eqref{eq:res6} into 3 sufficient equations
\begin{subequations}
\begin{align}
&\sum_{m=1}^{N+1}\left[a_m'(\pm 1)\cosh(\rho_m)\pm b_m'(\pm 1)\sinh(\rho_m)\right]=0 \label{eq:res7}\\
&\sum_{m=1}^{N+1}\left[\rho_m b_m(1) \cosh(\rho_m)+ \rho_m a_m(1) \sinh(\rho_m)\right]=0 \label{eq:res8}\\
&\sum_{m=1}^{N+1}\left[\rho_m b_m(-1) \cosh(\rho_m)- \rho_m a_m(-1) \sinh(\rho_m)\right]=0 \label{eq:res9}
\end{align}
\end{subequations}
Note that the equations \eqref{eq:res8} and \eqref{eq:res9} are equivalent to
\begin{equation}
\begin{split}
\sum_{m=1}^{N+1}a_m(-1)S^{z,even}_{N+1,m}&-\sum_{m=1}^{N+1}b_m(-1)S^{z,odd}_{N+1,m}=0\\
\sum_{m=1}^{N+1}a_m(1)S^{z,even}_{N+1,m}&+\sum_{m=1}^{N+1}b_m(1)S^{z,odd}_{N+1,m}=0
\end{split}
\label{eq:res10}
\end{equation}
If we combine the equations \eqref{eq:res10} with terms in square brackets in \eqref{eq:res3} and \eqref{eq:res4} we get the matrix equations.
\begin{equation}
\begin{split}
S^{z,even}\mathbf{a}(-1)-S^{z,odd}\mathbf{b}(-1)&=\mathbf{0}\\
S^{z,even}\mathbf{a}(1)+S^{z,odd}\mathbf{b}(1)&=\mathbf{0} 
\end{split}
\label{eq:res11}
\end{equation}

The term in square brackets in \eqref{eq:res1} vanishes if the following two equations vanish.
\begin{subequations}
\begin{align}
\frac{\partial}{\partial x}&\sum_{m=1}^{N+1}\left[a_m'(x)\cosh(\rho_m x)+b_m'(x)\sinh(\rho_m x)\right]=0 \label{eq:res12}\\
&\sum_{m=1}^{N+1}\left[\rho_m b_m'(x)\cosh(\rho_m x)+\rho_m a_m'(x)\sinh(\rho_m x)\right]=0 \label{eq:res13}
\end{align}
\end{subequations}

We see that in equation \eqref{eq:res12} the sum should be constant. Using equation \eqref{eq:res7} we see that this constant is zero.
\begin{equation}\label{eq:res14}
\sum_{m=1}^{N+1}\left[a_m'(x)\cosh(\rho_m x)+b_m'(x)\sinh(\rho_m x)\right]=0
\end{equation}

The remaining equations \eqref{eq:res5}, \eqref{eq:res13}, \eqref{eq:res14} form a system of differential equations with boundary conditions \eqref{eq:res11}.
\begin{equation}\label{eq:res15}
\begin{split}
&\sum_{m=1}^{N+1}\left[\frac{a_m'(x)}{k_j ^2-\rho_m^2}k_j \cosh(\rho_m x')+\frac{b_m'(x)}{k_j ^2-\rho_m^2}k_j \sinh(\rho_m x')\right]=0 \\
&\sum_{m=1}^{N+1}\left[\frac{a_m'(x)}{k_j ^2-\rho_m^2}\rho_m \sinh(\rho_m x')+\frac{b_m'(x)}{k_j ^2-\rho_m^2}\rho_m \cosh(\rho_m x')\right]=-R(z,B)y(x)\\
&\sum_{m=1}^{N+1}\left[\rho_m b_m'(x)\cosh(\rho_m x)+\rho_m a_m'(x)\sinh(\rho_m x)\right]=0\\
&\sum_{m=1}^{N+1}\left[a_m'(x)\cosh(\rho_m x)+b_m'(x)\sinh(\rho_m x)\right]=0
\end{split}
\end{equation}
We can rewrite these equations by introducing some matrices. We define the diagonal matrices $\hat{C}$, $\hat{S} \in C(\Omega,\mathbb{C}^{(N+1)\times (N+1)})$, the square matrices $\hat{K}$, $\hat{M}$, $\hat{Q} \in \mathbb{C}^{(N+1)\times (N+1)}$ and the operator $\hat{R} : Y\rightarrow Y^{N+1}$ as follows
\begin{equation}
\begin{split}
\hat{C}_{m,m}(x)&=\cosh(\rho_m x)\\
\hat{S}_{m,m}(x)&=\sinh(\rho_m x)\\
\hat{K}_{j,m}&=\rho_m \hat{Q}_{j,m}\\
\hat{M}_{j,m}&=k_j \hat{Q}_{j,m}\\
\hat{Q}_{j,m}&=\begin{cases}
\frac{1}{k_j ^2-\rho_m^2} &\quad \text{for } j\in \{1,\cdots,N\} \\
 1 &\quad \text{for } j=N+1\end{cases}\\
(\hat{R}y)_i&=\begin{cases}
R(z,B)y &\quad\text{for } j\in \{1,\cdots,N\} \\
0 &\quad\text{for } j=N+1\end{cases}
\end{split}
\label{eq:resdef}
\end{equation}
Here $j,m \in \{1,\cdots,N+1\}$ and we define $k_{N+1}:=1$. 

We seek functions $\mathbf{a}(x)$ and $\mathbf{b}(x)$ which solve the system of differential equations 
\begin{equation}
\begin{split}
\hat{M}(\hat{C}(x)\mathbf{a}'(x)+\hat{S}(x)\mathbf{b}'(x))&=\mathbf{0}\\
\hat{K}(\hat{S}(x)\mathbf{a}'(x)+\hat{C}(x)\mathbf{b}'(x))&=-\hat{R}y(x)
\end{split}
\label{eq:res16}
\end{equation}
with boundary conditions
\begin{equation}
\begin{split}
S^{z,even}\mathbf{a}(-1)-S^{z,odd}\mathbf{b}(-1)&=\mathbf{0}\\
S^{z,even}\mathbf{a}(1)+S^{z,odd}\mathbf{b}(1)&=\mathbf{0} 
\end{split}
\label{eq:res17}
\end{equation}

For $z\in \rho(A)$ we have that $S^{z,odd}$ and $S^{z,even}$ are invertible. Due to lemmas \ref{lem:nonzero} and \ref{lem:nonequal} when $z\notin \mathcal{S}$, $\hat{Q}$ satisfies the conditions of Lemma \ref{lem:Q_invert} and hence $\hat{Q}$ is invertible. We can write the determinant of $\hat{K}$ and $\hat{M}$ in terms of the determinant of $\hat{Q}$, $\det(\hat{M})=\det(\hat{Q})\prod_{j=1}^N k_j $, $|\hat{K}|=\det(\hat{Q})\prod_{m=1}^{N+1}\rho_m$ and so $\hat{K}$ and $\hat{M}$ are both invertible too.

Now we multiply the first line of \eqref{eq:res16} by $\hat{C}(x)\hat{M}^{-1}$ and second line by $\hat{S}(x)\hat{K}^{-1}$
\begin{equation}
\begin{split}
\hat{C}^2(x)\mathbf{a}'(x)+\hat{C}(x)\hat{S}(x)\mathbf{b}'(x)&=\mathbf{0} \\
\hat{S}^2(x)\mathbf{a}'(x)+\hat{C}(x)\hat{S}(x)\mathbf{b}'(x)&=-\hat{S}(x)\hat{K}^{-1}\hat{R}y(x)
\end{split}
\end{equation}

If we now subtract these equations and use the trigonometric identity $\hat{C}^2(x)-\hat{S}^2(x)=I$, we arrive at the following equation
\begin{equation}
\begin{split}
\mathbf{a}'(x)&=\hat{S}(x)\hat{K}^{-1}\hat{R}y(x)\\
\mathbf{b}'(x)&=-\hat{C}(x)\hat{K}^{-1}\hat{R}y(x)
\end{split}
\end{equation}

Here we get the second line by a similar procedure. We note that $\hat{R}y \in C^2(\Omega)$ and $A(x),B(x) \in C^\infty(\Omega)$, which implies that $\mathbf{a}(x), \mathbf{b}(x) \in C^3(\Omega)$. Hence we satisfy the regularity condition.

We can now find $\mathbf{a}(x)$ and $\mathbf{b}(x)$ by taking an anti-derivative plus some constants of integration, $\mathbf{a}^c$ and $\mathbf{b}^c$. To satisfy the boundary equations \eqref{eq:res17}, we take an anti-derivative such that $\mathbf{a}(-1)+\mathbf{a}(1)=2\mathbf{a}^c$ and $\mathbf{b}(-1)+\mathbf{b}(1)=2\mathbf{b}^c$.
\begin{equation}
\begin{split}
\mathbf{a}(x)&=\mathbf{a}^c+\frac{1}{2}\left(\int_{-1}^x \hat{S}(x')\hat{K}^{-1} \hat{R}y(x')\,dx'-\int_{x}^1 \hat{S}(x')\hat{K}^{-1} \hat{R}y(x')\,dx'\right)\\
\mathbf{b}(x)&=\mathbf{b}^c-\frac{1}{2}\left(\int_{-1}^x \hat{C}(x')\hat{K}^{-1} \hat{R}y(x')\,dx'-\int_{x}^1 \hat{C}(x')\hat{K}^{-1} \hat{R}y(x')\,dx'\right)
\end{split}
\label{eq:res_sol_part}
\end{equation}

By adding and subtracting the boundary equations \eqref{eq:res17} we find that the constants of integration equal
\begin{equation}
\begin{split}
\mathbf{a}^c&=\;\;\, \frac{1}{2}(S^{z,even})^{-1}S^{z,odd}\left(\int_{-1}^1 \hat{C}(x')\hat{K}^{-1} \hat{R}y(x')\,dx'\right)\\
\mathbf{b}^c&=-\frac{1}{2}(S^{z,odd})^{-1}S^{z,even}\left(\int_{-1}^1 \hat{S}(x')\hat{K}^{-1} \hat{R}y(x')\,dx'\right) 
\end{split}
\label{eq:resconst}
\end{equation}

We can simplify this as 
\begin{equation}
\begin{split}
\mathbf{a}(x)&=\;\;\,\frac{1}{2}\int_{-1}^1 \left(\hat{S}(x')\mathrm{sgn}(x-x')+(S^{z,even})^{-1}S^{z,odd}\hat{C}(x')\right)\hat{K}^{-1} \hat{R}y(x')\,dx'\\
\mathbf{b}(x)&=-\frac{1}{2}\int_{-1}^1 \left(\hat{C}(x')\mathrm{sgn}(x-x')+(S^{z,odd})^{-1}S^{z,even}\hat{S}(x')\right)\hat{K}^{-1} \hat{R}y(x')\,dx'
\end{split}
\label{eq:res_sol}
\end{equation}
\qed

For the computation of the first Lyapunov coefficient $l_1$ we need to evaluate the Dunford integral in \eqref{eq:normal_form_computation}. Similar to \citet{dijkstra_pitchforkhopf_2015} we can use residue calculus to find an expression for this integral.

\begin{theorem}\label{thm:residue}
Let $\lambda\in \sigma_p(A)$ be a simple eigenvalue and $\lambda \notin \mathcal{S}$. Let $C_\lambda$ be a sufficiently small closed disk such that $C_\lambda\cap \sigma(A)=\{\lambda\}$ and $ C_{\lambda} \cap \mathcal{S}=\emptyset$. 

If $\lambda$ is an 'even' eigenvalue with eigenvector
\begin{equation}
\psi(\theta)(x)= e^{\lambda \theta}\sum_{m=1}^{N+1} a_m \cosh(\rho_m(\lambda) x)
\end{equation}
where $\mathbf{a}$ is a non-trivial solution of $S^{\lambda,even}\mathbf{a}=0$. Then 
\begin{equation}\label{eq:complint}
\frac{1}{2\pi i} \oint_{\partial C_\lambda} e^{z \theta} \Delta^{-1}(z)y\, \,dz= \nu \psi(\theta)
\end{equation}
if and only if
\begin{equation}
\frac{\mathrm{adj}(S^{\lambda,even})}{2 \frac{d}{dz}(\det(S^{\lambda,even}))|_{z=\lambda}}S^{\lambda,odd} \int_{-1}^{1}\hat{C}(x')\hat{K}^{-1} \hat{R}y(x')\,dx'=\nu \mathbf{a}
\end{equation}
For all $y\in Y$, where $\mathrm{adj}(S^{\lambda,even})$ denotes the adjugate of $S^{\lambda,even}$ and using the definitions in \eqref{eq:resdef}. 

If $\lambda$ is an 'odd' eigenvalue with eigenvector
\begin{equation}
\psi(\theta)(x)= e^{\lambda \theta} \sum_{m=1}^{N+1} b_m \sinh(\rho_m(\lambda) x)
\end{equation}
where $\mathbf{b}$ is a non-trivial solution of $S^{\lambda,odd}\mathbf{b}=0$. Then 
\begin{equation}
\frac{1}{2\pi i} \oint_{\partial C_\lambda} e^{z \theta} \Delta^{-1}(z)y\, \,dz= \nu \psi(\theta)
\end{equation}
if and only if
\begin{equation}
\frac{-\mathrm{adj}(S^{\lambda,odd})}{2 \frac{d}{dz}(\det(S^{z,odd}))|_{z=\lambda}}S^{\lambda,even} \int_{-1}^{1}\hat{B}(x')\hat{K}^{-1} \hat{R}y(x')\,dx'=\nu \mathbf{b}
\end{equation}
For all $y\in Y$, where $\mathrm{adj}(S^{\lambda,odd})$ denotes the adjugate of $S^{\lambda,odd}$ and using the definitions in \eqref{eq:resdef}.
\label{eq:normal_form_computation2}
\end{theorem}
\textit{Proof.} As $\sigma_p(A)$ and $\sigma_p(B)$ contain only isolated eigenvalues and $\rho_m(z)$ and $\det(P^z(k_{i,j}(z)))$ are analytic in $z$, the set $\mathcal{S}$ contains only isolated values. Hence such a $C_{\lambda}$ exists. 
 
Suppose $\lambda$ is an even eigenvalue. As $\mathcal{S}\cap C_{\lambda}=\emptyset$ and $\sigma(A)\cap C_{\lambda}=\{\lambda\}$, we have that the $\Delta^{-1}(z)\mathbf{y}$ is given by Theorem \ref{thm:resolvent} for $z\in C_\lambda$. We observe that all components of the resolvent are analytic for all $z\in C_\lambda$ expect for the constants of integration $\mathbf{a}^c(z)$. This analyticity simplifies \eqref{eq:complint} to 
\[\frac{e^{\lambda \theta}}{2\pi i} \sum_{m=1}^{N(N+1)} \cosh(\rho_m(\lambda) x) \oint_{\partial C_\lambda} a_m^c(z) \,dz= \nu e^{\lambda \theta} \sum_{m=1}^{N(N+1)} a_m \cosh(\rho_m(\lambda) x)\]
for all $x\in \Omega, \theta \in [-h,0]$. We can substitute \eqref{eq:resconst} and use the residue formula
\[\frac{1}{2\pi i} \oint_{\partial C_\lambda} (S^{z,even})^{-1} \,dz = \mathrm{Res}\left(\frac{\mathrm{adj}(S^{z,even})}{\det(S^{z,even})},\lambda\right)=\frac{\mathrm{adj}(S^{\lambda,even})}{\frac{d}{dz}(\det(S^{z,even}))|_{z=\lambda}}\]
Due to linear independence of $\cosh(\rho_m(\lambda) x)$ for $m \in \{1,\cdots, N+1\}$, this results in the formula
\[\frac{\mathrm{adj}(S^{\lambda,even})}{2 \frac{d}{dz}(\det(S^{z,even}))|_{z=\lambda}}S^{\lambda,odd} \int_{-1}^{1}\hat{C}(x')\hat{K}^{-1} \hat{R}y(x')\,dx'=\nu \mathbf{a}\]
The reasoning for odd eigenvalues is similar. \qed

\section{Numerical Results}
\label{sec:num}
In this section we will examine a specific numerical example. We will compute eigenvalues and the first Lyapunov coefficient for a Hopf bifurcation and investigate the effect of varying the diffusion parameter $d$.

\begin{figure}[ht!]
\includegraphics[width=0.65\linewidth]{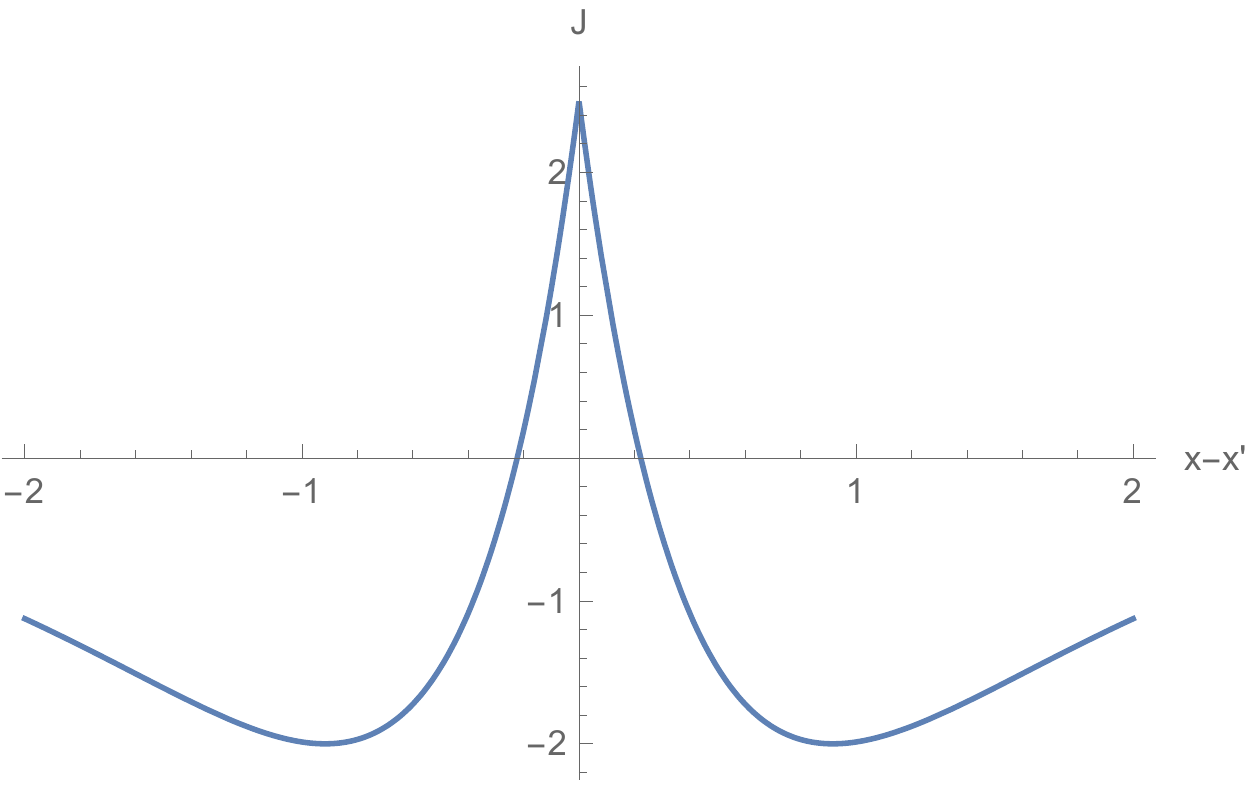}
\caption{The wizard-hat connectivity of \eqref{eq:connectivity}}\label{Figure2}
\end{figure}

For $J$ we choose the following difference of two exponentials, as  \citet{dijkstra_pitchforkhopf_2015}
\begin{equation}\label{eq:connectivity}
J(x,x')= \frac{25}{2}e^{-2|x-x'|}-10 e^{-|x-x'|}
\end{equation}
This connectivity is a model of a population of excitatory neurons acting on a short distance combined with a population of inhibitory neurons acting on a longer distance, see figure \ref{Figure2}

For the activation function $S$ we choose the sigmoidal function 
\begin{equation}
S(u)=\frac{1}{1+e^{-\gamma u}}-\frac{1}{2}
\end{equation}
As $S$ is an odd function, $S''(0)=0$ and hence $D^2G(0) \equiv 0$. This simplifies the computation of first Lyapunov coefficient $l_1$ of \eqref{eq:normal_form_computation} to
\begin{equation}
\frac{1}{4\pi i} \oint_{\partial C_\lambda} e^{z\theta} \Delta^{-1}(z)D^3G(0)(\psi,\psi,\bar{\psi})\,dz= c_1 \psi(\theta)
\end{equation}
We can compute this integral using Theorem \ref{thm:residue} with $y= \tfrac{1}{2}D^3G(0)(\psi,\psi,\bar{\psi})$.

We fix the following values for parameters $\alpha=1$ and $\tau^0=\frac{3}{4}$ and use $\gamma$ as the bifurcation parameter. We want to compare two cases: without diffusion, i.e. $d=0$, and with diffusion, i.e. $d>0$.
 
\subsection{Hopf bifurcation}
For $d=0$ we have an Hopf-bifurcation for $\gamma=3.3482$ at $\lambda=1.2403i$ with corresponding eigenvector 
\begin{equation}\label{eq:eigenvector1}
\begin{split}
\psi(\theta)(x)=e^{1.2403 i \theta}[&0.9998\cosh((0.2770-0.8878i)x)\\
&(-0.0178+0.0050i)\cosh((3.7185+3.2284i)x)]
\end{split}
\end{equation}
The normal form coefficient $c_1=-1.132-0.282i$ and the Lyapunov coefficient $\ell_1=-0.9123$ and hence the bifurcation is supercritical.

For $d=0.2$ we have an Hopf-bifurcation for $\gamma=3.3094$ at $\lambda=1.2379i$ with corresponding eigenvector 
\begin{equation}
\begin{split}\label{eq:eigenvector2}
\psi(\theta)(x)=e^{1.2379 i \theta}[&0.9972\cosh((0.2535-0.8490i)x)\\
&+(-0.0727-0.0177i)\cosh((1.7315+3.2475i)x)\\
&+(0.0029-0.0060i)\cosh((3.90746+0.3586i)x)]
\end{split}
\end{equation}
The normal form coefficient $c_1=-1.153-0.258i$ and the Lyapunov coefficient $\ell_1=-0.9314$ and hence the bifurcation is also supercritical. We have put these values for further reference in table \ref{tab:Hopfvalues}.

\begin{figure}
\centering
\includegraphics[scale=0.65]{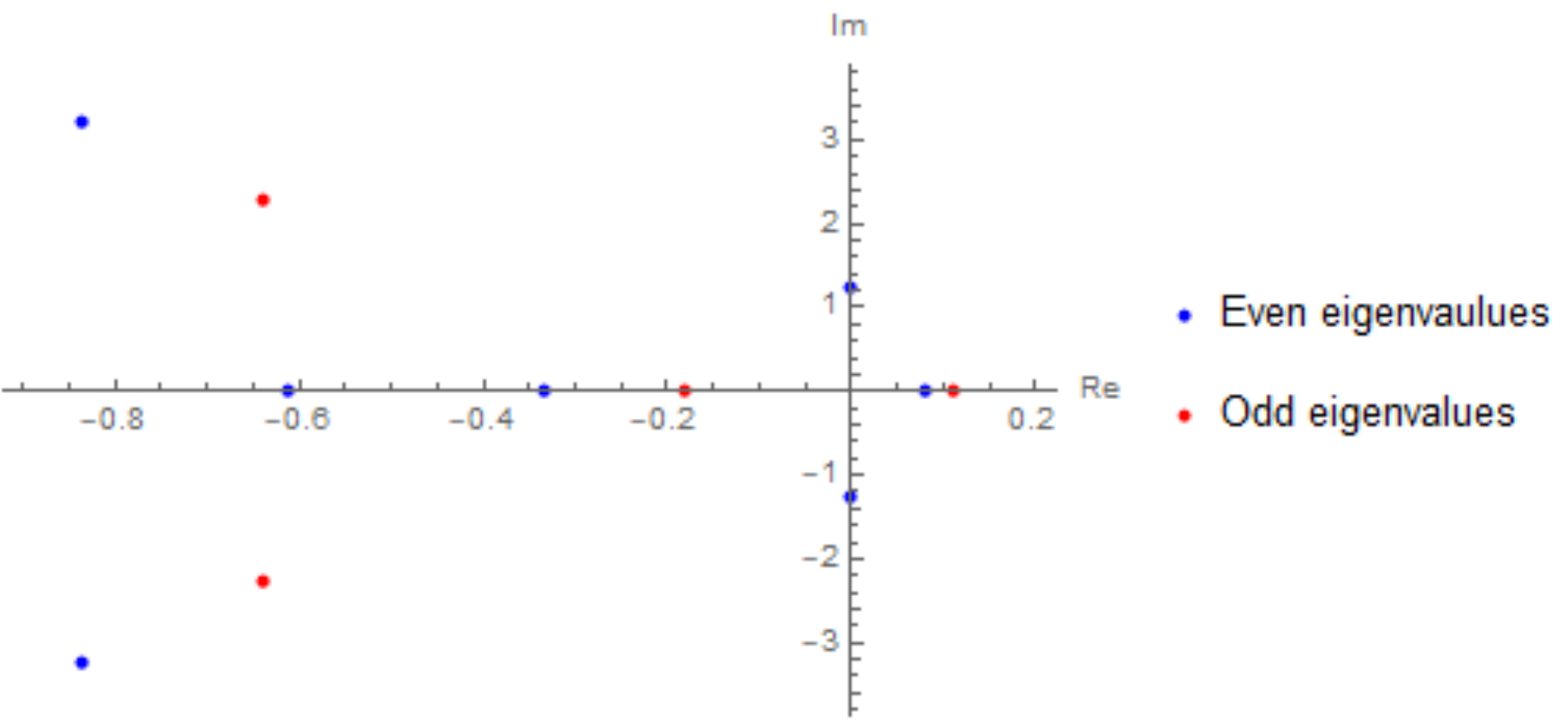}
\includegraphics[scale=0.65]{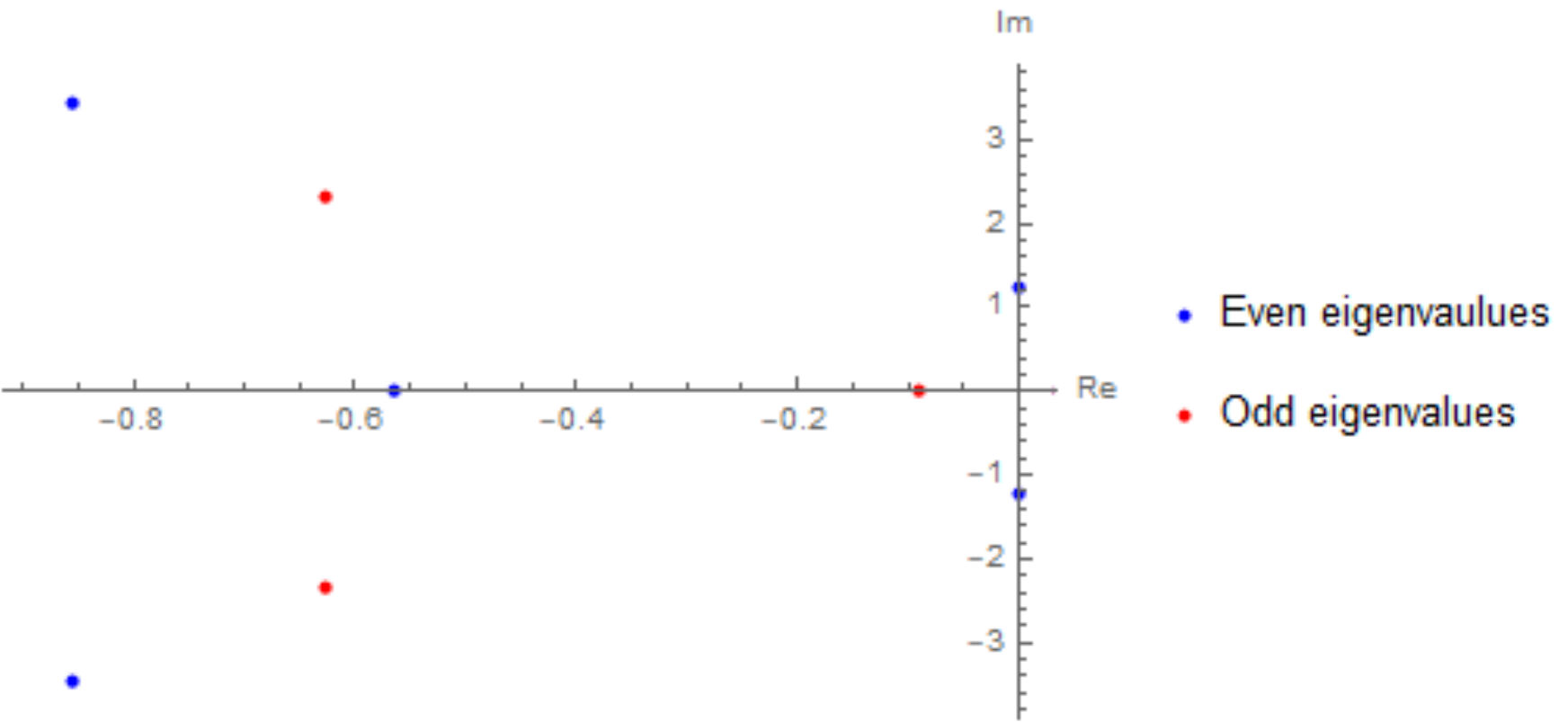}
\caption{The eigenvalues of $A$ at parameter values in table \ref{tab:Hopfvalues} of the Hopf bifurcation without and with diffusion respectively.}
\label{Figure3}
\end{figure}

\begin{figure}
\centering
\includegraphics[scale=0.6]{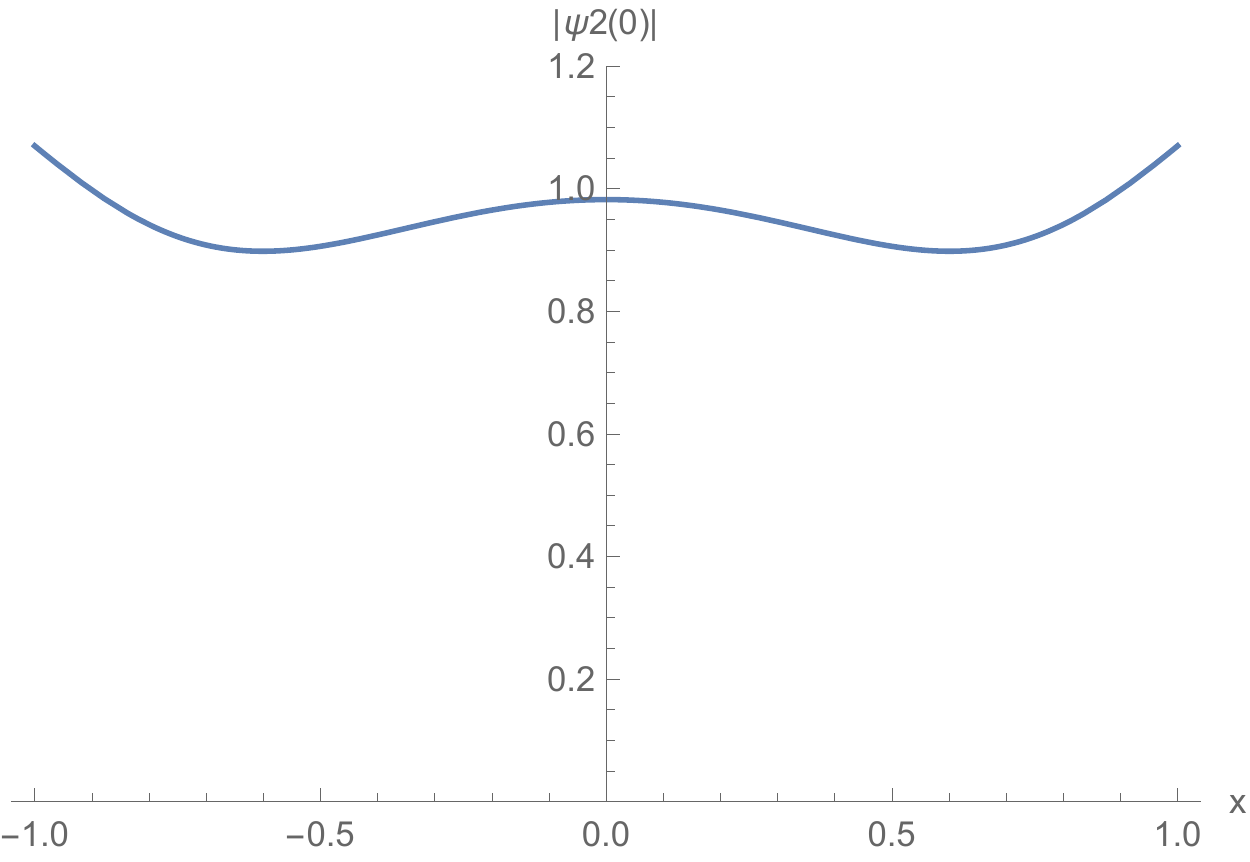}
\includegraphics[scale=0.6]{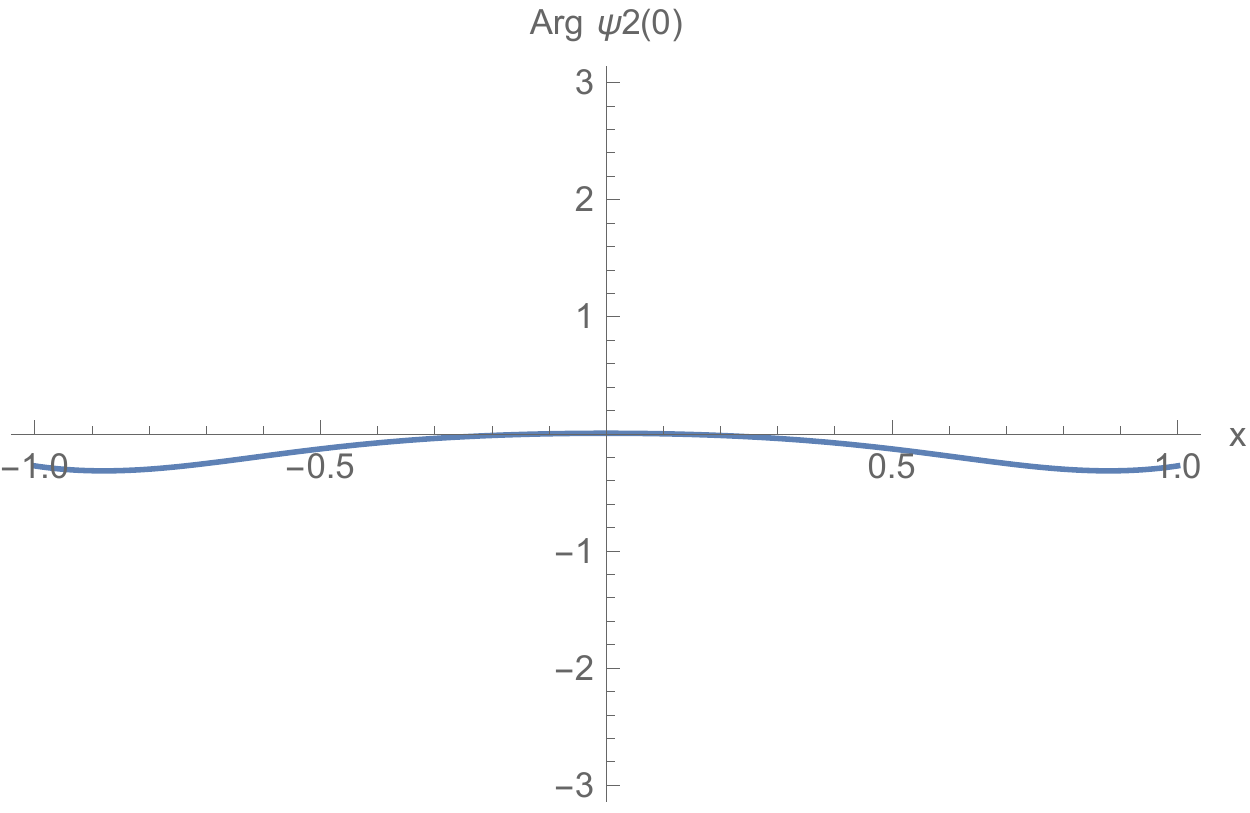}
\includegraphics[scale=0.6]{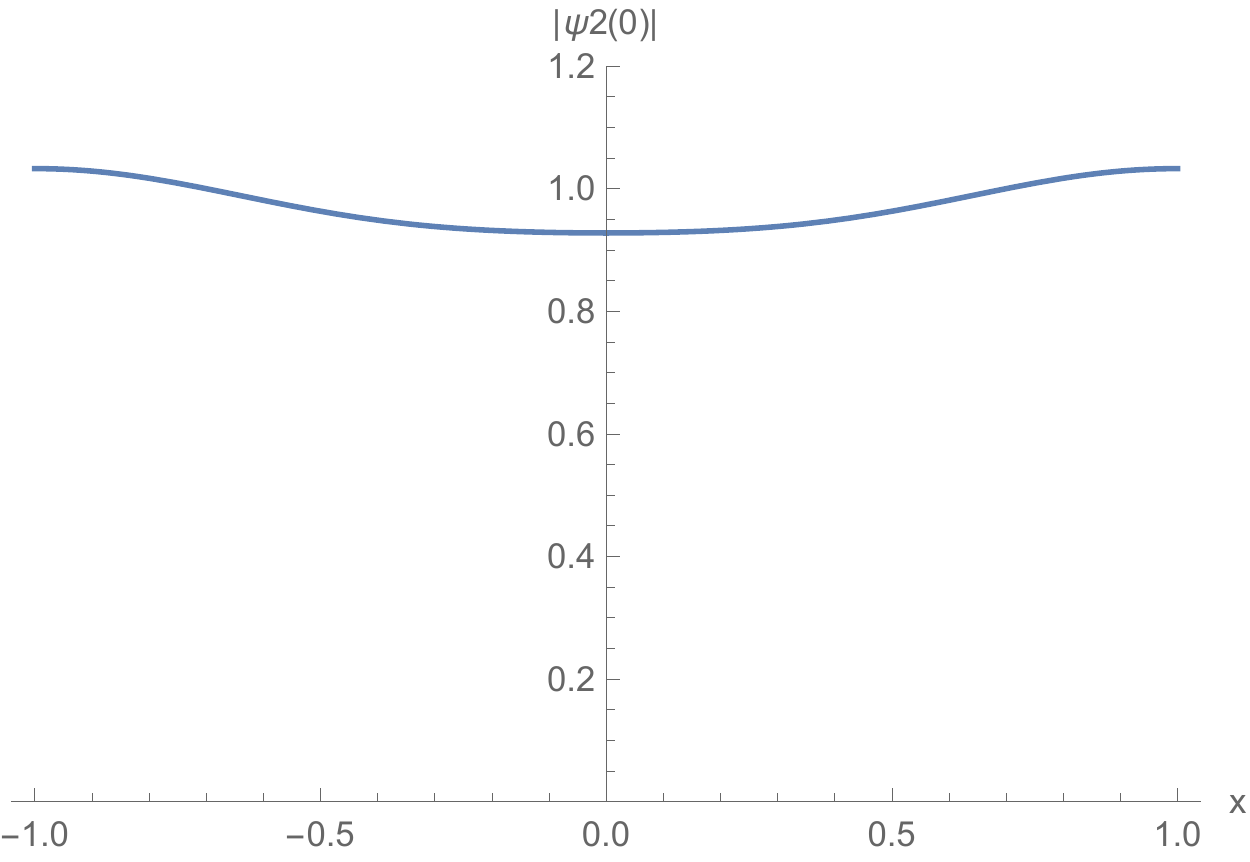}
\includegraphics[scale=0.6]{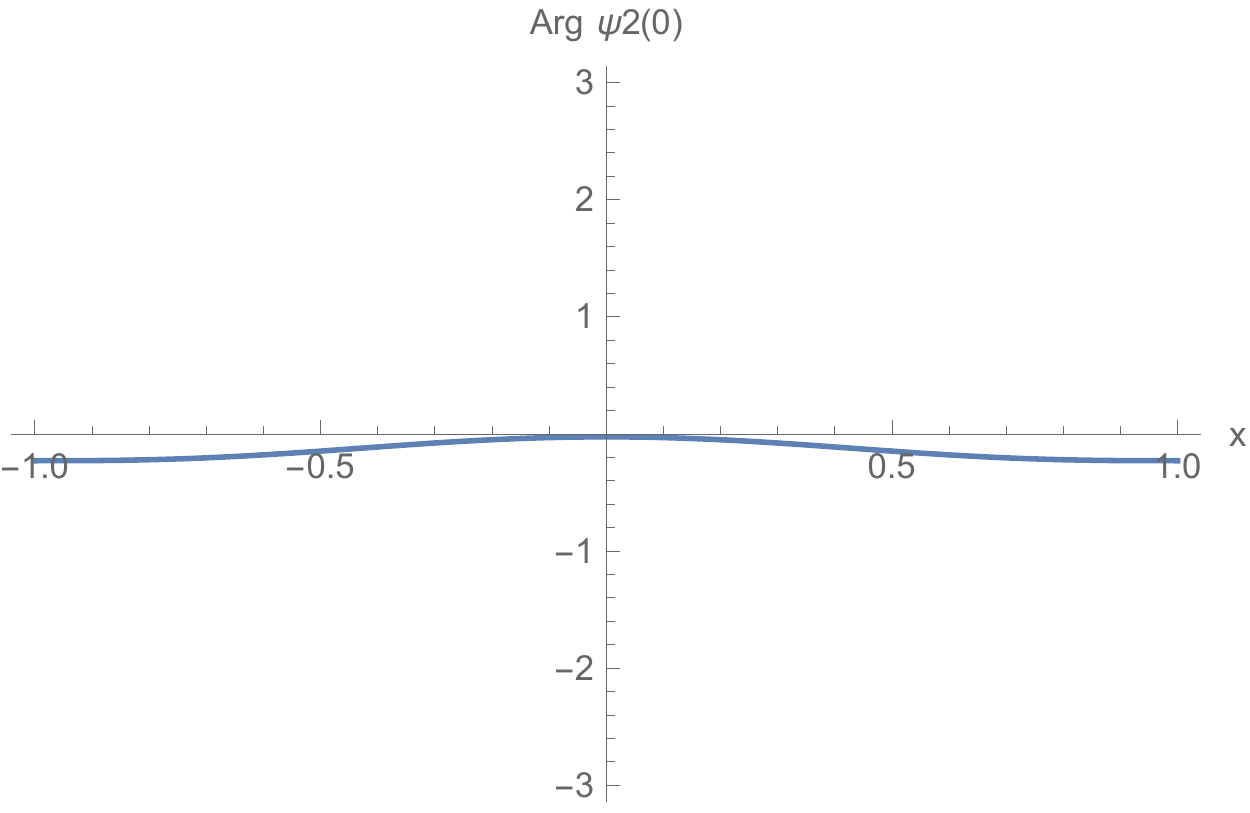}
\caption{The corresponding eigenvectors of the eigenvalue $\lambda = \omega i$ at parameter values in table \ref{tab:Hopfvalues} without and with diffusion respectively. Note that with diffusion the eigenvector satisfies the boundary conditions at $x=1$ and $x=-1$, while this is not the case without diffusion.}
\label{Figure4}
\end{figure}

As one might already have observed, the diffusion has little effect on the Hopf bifurcation. We observe more generally that the eigenvalues which are off the real axis are barely effected by the introduction of diffusion, while the eigenvalues on the real axis become more negative, see figure \ref{Figure3}.\footnote{Note that there is another positive $\lambda \in \mathbb{R}$, not shown in figure \ref{Figure3}, which solves $\det(S^{\lambda,odd})=0$ and $\det(S^{\lambda,even})=0$, however this is a degenerate case as $P^\lambda(\rho)$ has a double root. Simulations of the linearised system did not indicate the presence of an unstable mode, so we don't regard this point as an eigenvalue.} A possible explanation is that the eigenvector corresponding to the eigenvalue on the imaginary axis has very little spatial curvature, see figure \ref{Figure4}. As diffusion penalises curvature, its effect on this eigenvector would be small.

\begin{table}
\centering
\begin{tabular}{|c|cccccccccc|}
\hline
Bifurcation & $\alpha$ & $\tau^0$ & $\eta_1$ & $\eta_2$ & $\mu_1$ & $\mu_2$ & $d$ & $\gamma$ & $\lambda$ & $\ell_1$\\
\hline
Hopf 1 & 1 & 0.75 & 12.5 & -10 & 2 & 1 & 0 & 3.3482 & 1.2403i & -0.9123\\
Hopf 2 & 1 & 0.75 & 12.5 & -10 & 2 & 1 & 0.2 & 3.3094 & 1.2379i & -0.9314\\
\hline
\end{tabular}
\caption{Parameter values of the Hopf bifurcation without and with diffusion respectively.}
\label{tab:Hopfvalues}
\end{table}

\subsection{Discretisation}
To obtain an approximate solution of \eqref{ADDE} we discretise the spatial domain $\Omega$ into an equidistant grid of $n^x$ points, $x_1,\, \dots \, , x_{n^x}$, with a width of $\delta = \tfrac{2}{n^x-1}$. Like \citet{faye_theoretical_2010}, we discretise the integral operator $G$ using the trapezoidal rule and the diffusion operator $B$ using a central difference method and a reflection across the boundary for the boundary conditions. This results in a second order spatial discretisation. The discretisation of the \eqref{ADDE} for $n \in \{1,\cdots, n^x\}$ and $t\in \mathbb{R}^+$ becomes a set of delay equations \eqref{DDE} 
\begin{equation}
\begin{cases}
\frac{\partial u}{\partial t}(t,x_n)&= \frac{d}{2\delta^2} (u(t,x_{n-1}) - 2 u(t,x_n) + u(t,x_{n+1})) - \alpha u(t,x_n)\\
&+ \delta \sum_{m=1}^{n^x} \xi_m J(x_n,x_m)S(u(t-\tau(x_n,x_m),x_m)) \\
u(t,x_{0})&= u(t,x_{2})\\
u(t,x_{n^x+1})&= u(t,x_{n^x-1})\\
u(t,x_n)&=\varphi(t,x_n)
\end{cases}
\label{DDE} \tag{DDE}
\end{equation}
Here $\xi_m$ is defined as
\begin{equation}
\xi_m= \begin{cases}
1 &\qquad m \in \{2, \cdots, n^x-1\}\\
\frac{1}{2}&\qquad m=1 \text{ or } m=n^x
\end{cases}
\end{equation}
Now we are left with a set of $n^x$ ordinary delay differential equations which we solve with a standard DDE-solver. Note that the \eqref{DDE} is very similar to the discrete model \eqref{eq:model3} from which the \eqref{ADDE} is derived. Only the terms at the boundary are different due to the second order discretisation.

\subsection{Simulations}
\begin{figure}[htbp]
\centering
\includegraphics[scale=0.65]{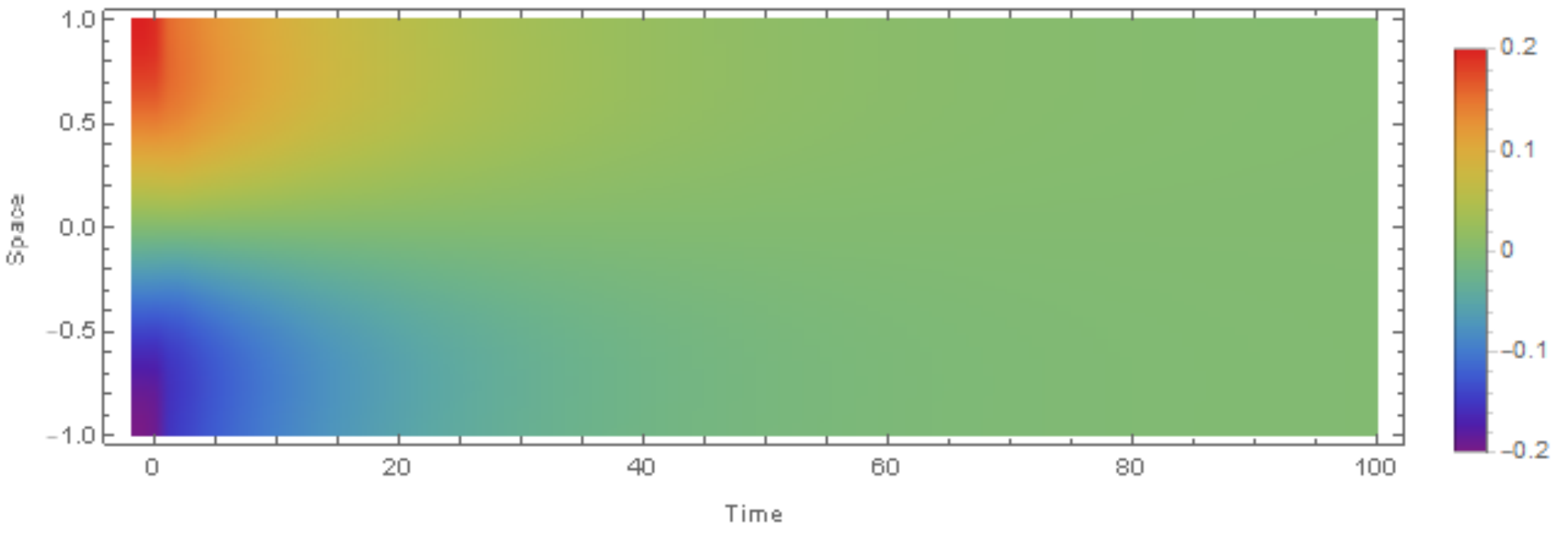}
\includegraphics[scale=0.65]{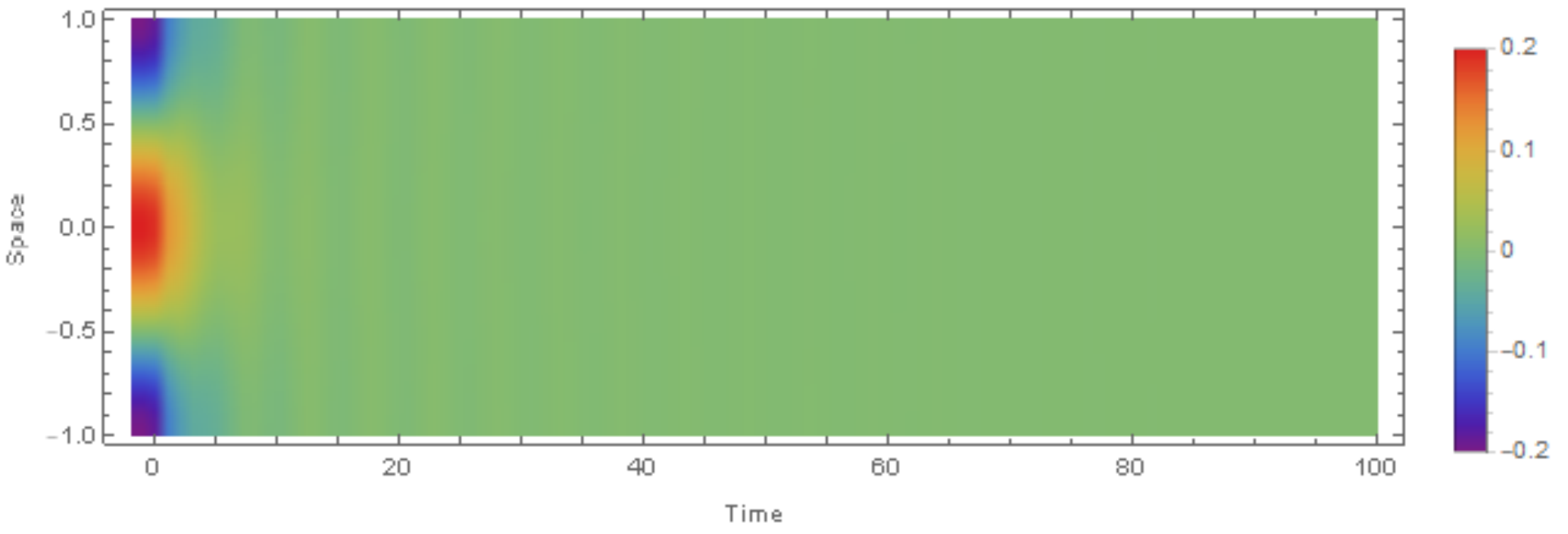}
\caption{Simulation of \eqref{DDE} with the initial conditions $\varphi_1, \varphi_2$ of \eqref{eq:initialcondition} and $\gamma = 3$ and $d=0.2$ }
\label{Figure5}
\end{figure}

\begin{figure}[htbp]
\centering
\includegraphics[scale=0.65]{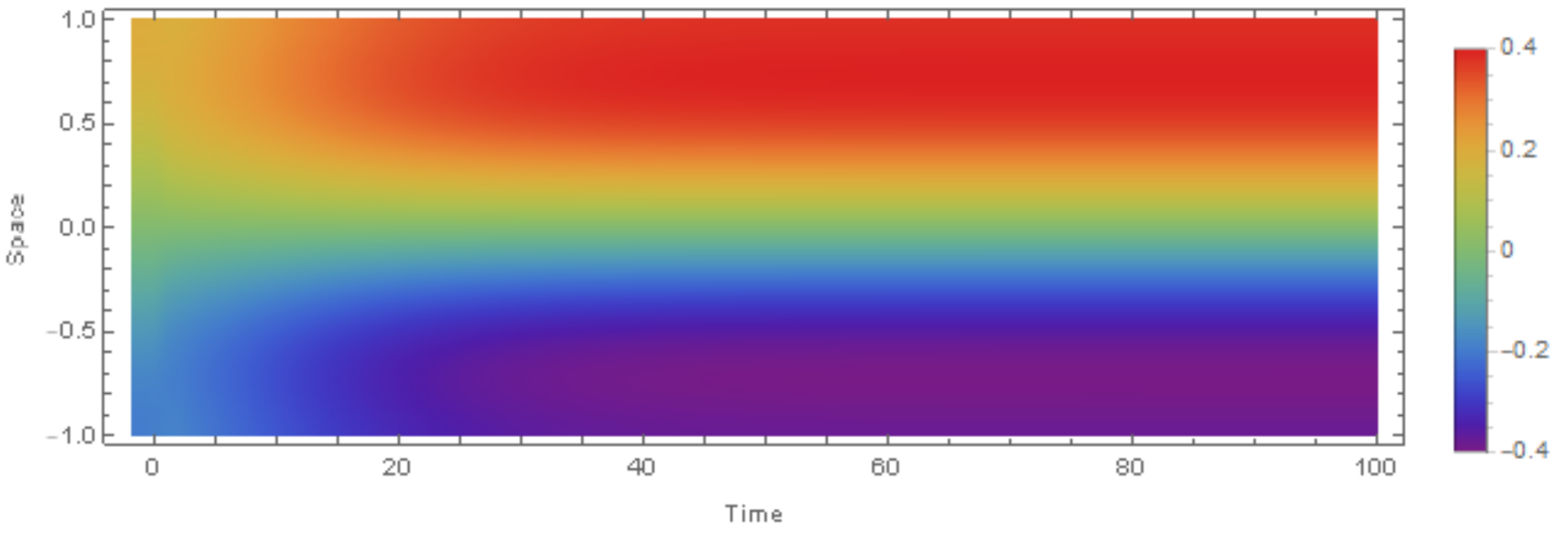}
\includegraphics[scale=0.65]{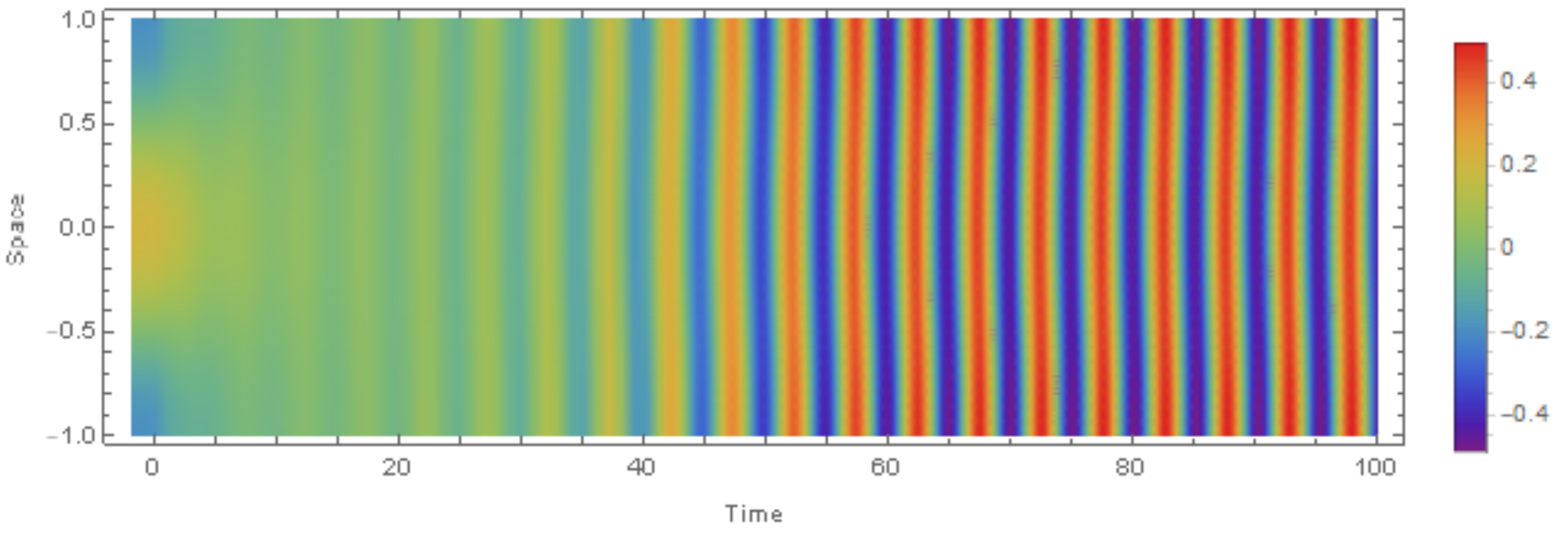}
\caption{Simulation of \eqref{DDE} with the initial conditions $\varphi_1, \varphi_2$ of \eqref{eq:initialcondition} and $\gamma = 4$ and $d=0.2$}
\label{Figure6}
\end{figure}

We will now perform some simulations around the Hopf bifurcation with diffusion. We set $n^x=50$ and take as initial conditions an odd function and an even function,
\begin{equation}\label{eq:initialcondition}
\begin{split}
\varphi_1(\theta)(x)&= \frac{1}{5}\sin{\frac{1}{2}\pi x}\\
\varphi_2(\theta)(x)&= \frac{1}{5}\cos{\pi x}\\
\end{split}
\end{equation}
For figure \ref{Figure5} we took $\gamma=3$ and for figure \ref{Figure6} $\gamma=4$. 

\begin{figure}[htbp]
\centering
\includegraphics[scale=0.75]{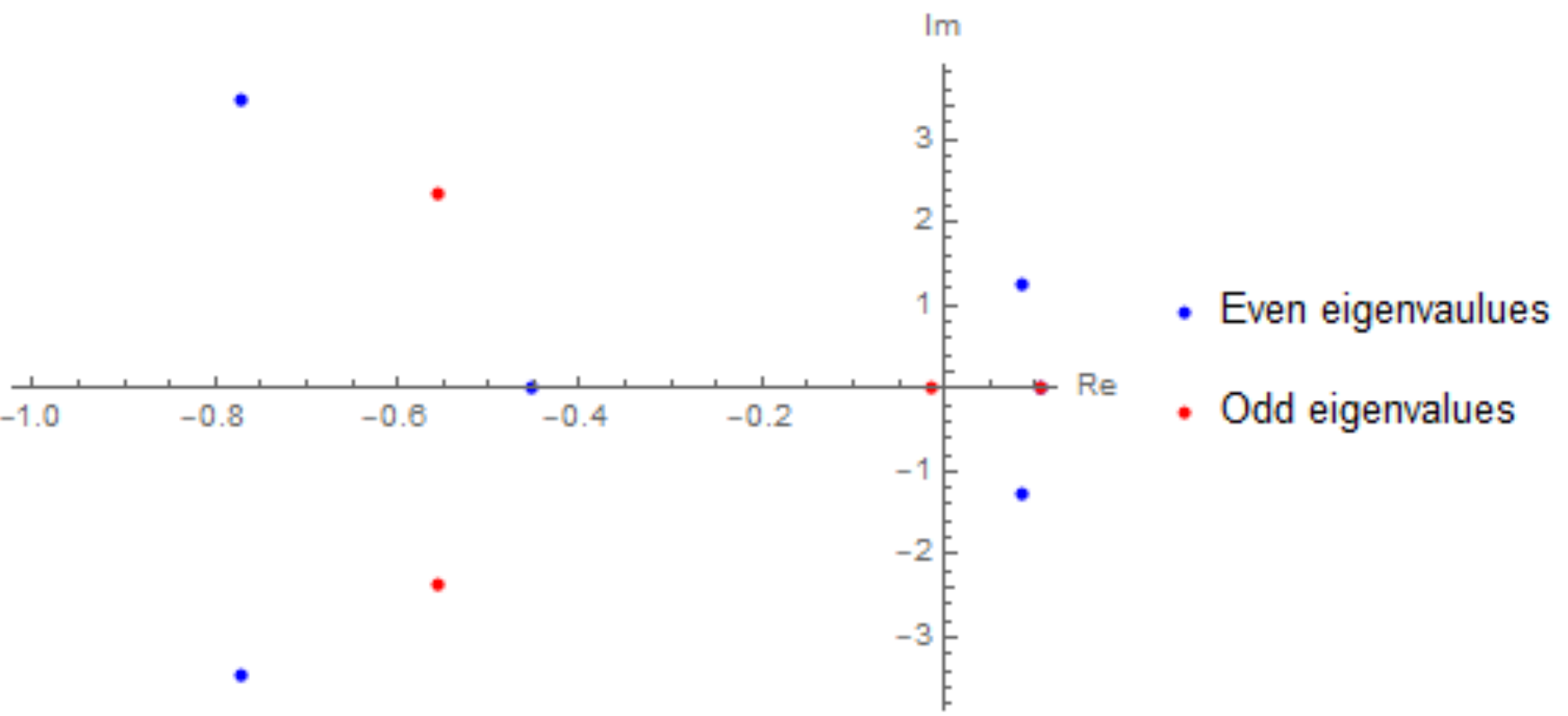}
\caption{The eigenvalues of $A$ for $\gamma = 4$ and $d=0.2$}
\label{Figure7}
\end{figure}

For $\gamma=3$, the solutions with both initial conditions \eqref{eq:initialcondition} converge to the trivial equilibrium. The one with the odd initial condition converges monotonously to the trivial equilibrium, while the one with the even initial condition converges to the trivial equilibrium in an oscillatory manner. For $\gamma=4$, there are (at least) two non-trivial stable states. The odd initial condition converges to some non-trivial equilibrium and the even initial condition converges to some limit cycle, which is due to the Hopf bifurcation. This is similar to the results of \citet{dijkstra_pitchforkhopf_2015}, where the non-trivial equilibrium arises from a pitchfork bifurcation. The bi-stability is also exemplified in the eigenvalues, see figure \ref{Figure7}, as we have a positive real eigenvalue and a pair of complex eigenvalues with a positive real component. 

\begin{figure}
\centering
\includegraphics[scale=0.65]{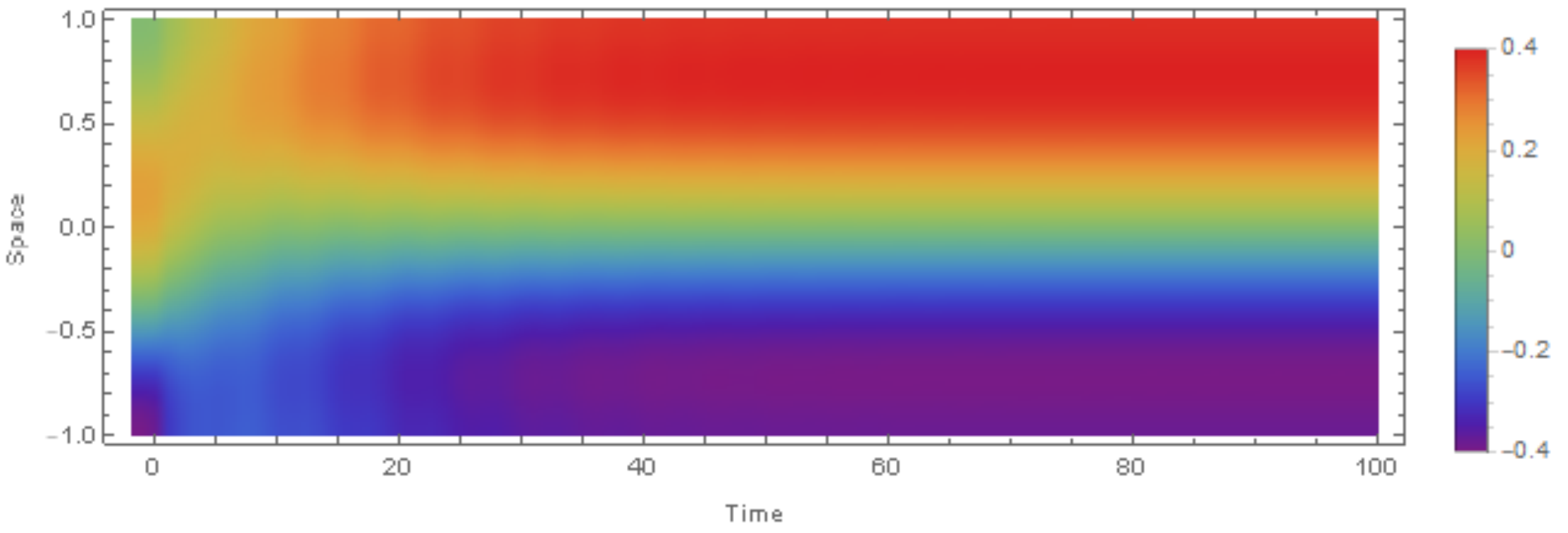}
\includegraphics[scale=0.65]{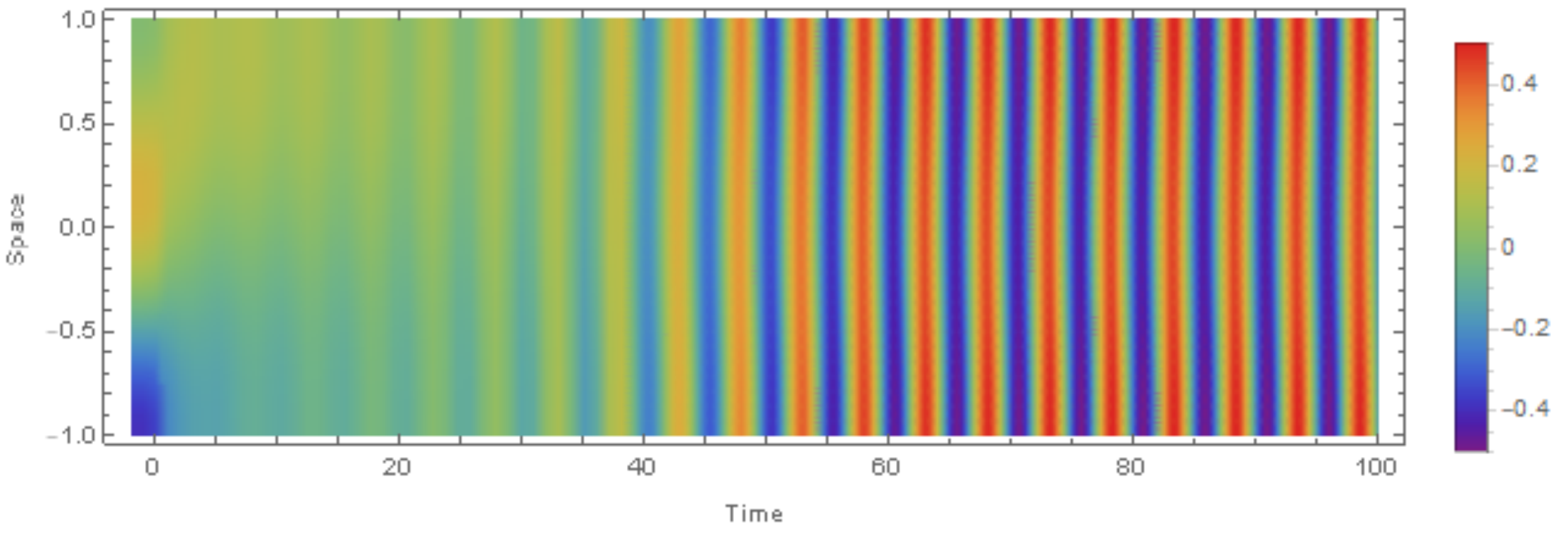}
\caption{Simulation of \eqref{DDE} with the same initial condition $\varphi_3$ \eqref{eq:initialcondition2} and $\gamma = 4$, $d=0.2$ and $\gamma = 4$, $d=0.5$ respectively.}
\label{Figure8}
\end{figure}

We have seen that increasing the value of $d$, decreases the eigenvalues on the real axis. This would imply that the non-trivial equilibrium becomes unstable or disappears, probably through a pitchfork bifurcation. Indeed when we use the initial condition
\begin{equation}\label{eq:initialcondition2}
\varphi_3 = \varphi_1 + \varphi_2\\
\end{equation}
and compare the dynamics for $d=0.2$ and $d=0.5$ in figure \ref{Figure8}. The initial condition converges to a non-trivial equilibrium when $d=0.2$, but it converges to a limit cycle when $d=0.5$.

\section{Discussion}\label{sec:discussion}
We have proved the necessary theorems to construct the sun-star calculus for abstract delay differential equations. In particular we proved a novel characterisation for sun-reflexivity in Theorem \ref{thm:sun-reflexivity}. The sun-star calculus provides a variation-of-constants formulation for the nonlinear problem and produces results on the spectral properties of the system, notably the essential spectrum. Using the results of  \citet{janssens_class_2020} on the center manifold reduction, we have derived a simple and explicit formula to compute the first Lyapunov coefficient for the Hopf bifurcation. This procedure can quite easily be extended to normal coefficients of other local bifurcations.

The neural field models, both with and without diffusion, can be cast as abstract delay differential equations to which the same theoretical results can be applied. In the sun-star calculus the relevant spaces, duality pairings and Fredholm alternative follow naturally by considering the strong continuity of adjoint operators. Hence there is no need to construct formal projectors. Moreover, for a specific example of the neural field we could calculate the first Lyapunov coefficient exactly and with arbitrary precision. Thus we conclude that the sun-star calculus for delay equations is a natural setting to study neural field models, with and without diffusion.

For certain specific connectivity functions we have derived analytical conditions for $\lambda$ to be an eigenvalue for a neural field with a connectivity function that is a sum of exponentials. We have also constructed the corresponding eigenvectors and the resolvent. Numerical results show that the diffusion term does not cause oscillations to arise due to a Hopf-bifurcation. However, stable equilibria which are not uniform disappear due to the smoothing effect of the diffusion. So increasing the diffusion in a bi-stable system with a non-uniform equilibrium and a synchronous oscillation leads to a system with only stable synchronous oscillations. We hypothesise that this is a more general feature of equations with diffusion and a delayed reaction.

Gap junctions, modelled by the diffusion term in our neural field, are thought to be linked to synchronisation in Parkinson's disease \citep{schwab_pallidal_2014}. Further research could be undertaken to see whether the effects can be observed in a Neural Field Model with physiological values for the parameters. 

We used a neural field model with a connectivity function which is a sum of exponentials. This connectivity function is commonly used to aggregate the effect of multiple different types of cells, e.g. excitatory and inhibitory neurons. However, introducing a diffusion term into this model leads to gap junctions between similar and different populations of neurons of the same strength. This may not be physiologically feasible. A way to circumvent this is to use a neural field model with multiple populations. In such a model, it is possible to introduce only gap junctions between neurons of the same population.

We have studied a neural field on a 1-dimensional closed domain. However, when model\-ling the neuronal activity in the cortex, it is common to use  2-dimensional domains \citep{coombes_tutorial_2014}. For a neural field with a rectangular domain, characterising the spectrum, as is done in this paper in section \ref{sec:single_pop}, is still an open problem. On a spherical domain, \citet{visser_standing_2017} have characterised the spectrum for a neural field with transmission delays and have computed normal form coefficients of Hopf and double Hopf bifurcations. It seems possible to extend the analysis of that paper to include a diffusion term into that neural field model. Due to the general nature of the theoretical results of section \ref{sec:duality}, these results, including the sun-star framework, the variation of constants formulation and the essential spectrum, also hold for neural field models on arbitrary domains. 

\begin{appendix}
\section{Properties of the Diffusion Operator}
\label{sec:diffusion}
In this appendix we investigate the properties of the diffusion operator $B$ in the context of the sun-star calculus. We consider the space of continuous functions $Y=C(\Omega)$, where we take our domain $\Omega$ to be the interval $[-1,1]$. We define $B: D(B) \rightarrow Y$, an unbounded, closed, linear operator as
\begin{equation}
\begin{split}
B q &:= d q'' - \alpha q\\
D(B) &:= \{ q\in Y | q\in C^2(\Omega), q'(\partial \Omega)=0 \}
\end{split}
\end{equation}

\subsection{Spectral properties}
We start our analysis with this result on the semigroup $S$ generated by $B$.
\begin{lemma}\label{lem:prop_b}
\citep[Proposition VI.6.19]{engel_one-parameter_1999}. The operator $(B,D(B))$ generates a strongly continuous, positive and immediately compact semigroup $(S(t))_{t\geq 0}$.
\end{lemma} 

Sturm–Liouville theory gives the following well-known results on the spectral properties of the diffusion operator. We can explicitly derive the eigenvalues and eigenvectors using separation of variables. This is entirely standard and therefore the calculation is omitted.
\begin{lemma}\label{lem:spectrumB}
For the spectrum of $B$ we have that $\sigma(B)=\sigma_p(B)$. All eigenvalues of $B$ are simple and given by $\lambda_n^{even}=-dn^2\pi^2-\alpha$ with even eigenvector $\cos(n\pi x)$ and $\lambda_n^{odd}=-d(n+\tfrac{1}{2})^2\pi^2-\alpha$ with odd eigenvector $\sin((n+\tfrac{1}{2})\pi x)$ for all $n\in \mathbb{N}_0$. Moreover these eigenvectors form a maximal set in $Y$, i.e. their span is dense in $Y$.
\end{lemma}

We can also explicitly find an explicit representation of the semigroup $S(t)$ and resolvent $R(z,B)$ in terms of the eigenvectors.
\begin{lemma}
The semigroup $S$ can be explicitly written as a convolution 
\begin{equation}
S(t)\varphi(x)=\int_{\Omega}\varphi(x')G(t,x,x') \,dx'   
\end{equation}
with Green's function 
\begin{equation}
\begin{split}
G(t,x,x')&:=\sum_{n=0}^{\infty} \left((1+\delta_{0n})^{-1}\cos(n\pi x)\cos(n\pi x')e^{\left(-dn^2\pi^2-\alpha\right)t} \right. \\
 & \qquad \quad \left. +\sin((n+\tfrac{1}{2})\pi x)\sin((n+\tfrac{1}{2})\pi x')e^{\left(-d\left(n+\tfrac{1}{2}\right)^2\pi^2-\alpha\right)t} \right) 
\end{split}
\end{equation}
The resolvent $R(z,B): Y \rightarrow D(B)$ for $z \in \rho(B)$ can be explicitly written as a convolution
\begin{equation}
R(z,B)y(x)=\int_{\Omega}y(x')G^{z}(x,x')\,dx' \label{eq:res_b}   
\end{equation}
with  Green's function 
\begin{equation}
\begin{split}
G^z(x,x')&:=\sum_{n=0}^{\infty}\left( (1+\delta_{0n})^{-1}\left(z+\alpha+dn^2\pi^2\right)^{-1}\cos(n\pi x)\cos(n\pi x') \right. \\
&\qquad\quad \left. +(z+\alpha+d\left(n+\tfrac{1}{2})^2\pi^2\right)^{-1}\sin((n+\tfrac{1}{2})\pi x)\sin((n+\tfrac{1}{2})\pi x') \right)
\end{split}
\end{equation}
Here $\delta_{mn}$ is the Kronecker delta. 
\end{lemma}

\subsection{Sun-star calculus}
We will now develop the sun-star calculus for the diffusion operator $B$. We can take $d=1$ and $\alpha=0$ for this section, without loss of generality, as the sun-star calculus is invariant with respect to bounded perturbations of the generator of the semi-group.

As a consequence of the Riesz representation theorem, $Y^*$ can be represented as $NBV(\Omega)$, the functions of bounded variation, normalised such that for $y\in Y^*, y(-1)=0$. The corresponding norm on $NBV(\Omega)$ is the total variation norm and the duality pairing is given by the Riemann-Stieltjes integral:
\begin{equation}
\langle y^*,y \rangle:=\int_{-1}^1 y\,\,dy^*
\end{equation}
We will now try to find a representation for $B^*$. 

\begin{theorem}
The dual space $Y^*$ can be represented as $NBV(\Omega)$. Furthermore, $y^*\in D(B^*)$ if and only if for $x\in (-1,1]$ 
\begin{equation}
y^*(x)= c_1 + \int_{-1}^x \left(c_2 + \int_{-1}^s z^*(x') \,dx'\right) \,ds
\label{eq:ystar}
\end{equation}
Where $c_1,c_2 \in \mathbb{R}$ and $z^* \in NBV(\Omega)$ with $z^*(1)=0$. For such $y^*$ we have that $B^* y^*= z^*$
\end{theorem}
\textit{Proof.} We start by proving the `only if' part of the Theorem. Let $y^* \in D(B^*)$, $y\in D(B)$ and $z^*=B^*y^*$. Furthermore, let 
\[w^*(s):=c_2 + \int_{-1}^s z^*(x') \,dx'\]
for some $c_2 \in \mathbb{R}$. As $y\in C^2(\Omega)$ and $y'(\pm 1)=0$ we get that using integration by parts for Riemann-Stieltjes integrals \citep[Proposition A.15, A.18, A.19]{janssens_class_2019}
\begin{align*}
\int_{-1}^1y''(x)\,dy^*(x) &= \langle y^*, B y\rangle = \langle z^*, y\rangle = \int_{-1}^1y \,dz^* \\
&= \left.z^*(x)y(x)\right|_{-1}^1 - \int_{-1}^1y'(x)z^*(x)\,dx\\
&= z^*(1)y(1)+\int_{-1}^1y''(x)w^*(x)\,dx
\end{align*}
If we take $y$ as a constant function we immediately see that $z(1)=0$ is a necessary condition. For any $-1<x'<x<1$ we can take a sequence of $y_n \in D(B)$ such that $y_n''(s)$ converges monotone to the characteristic function on the interval $[x',x]$. Then by the Lebesque monotone convergence theorem we get that
\[y^*(x)-y^*(x')=\int_{x'}^x \,dy^*(s) = \int_{x'}^x w^*(s) \,ds\]
Letting $x'\downarrow -1$ we get that
\[y^*(x)= \lim_{x' \downarrow -1} y^* (x') + \int_{-1}^x w^*(s) \,ds\]
So we can write this $y^*$ as 
\[y^*(x)= c_1 + \int_{-1}^x \left(c_2 + \int_{-1}^s z^*(x') \,dx'\right) \,ds\]

Next we prove the `if' part of the Theorem. Let $y^*$ have the form in equation \eqref{eq:ystar} with $z(1)=0$. Then for all $y\in D(B)$ we have again by using integration by parts that
\begin{align*}
\langle y^*, B y\rangle &= \int_{-1}^1y''(x)\,dy^*(x)\\
&= \int_{-1}^1y''(x)w^*(x)\,dx\\
&= -\int_{-1}^1y'(x)z^*(x)\,dx\\
&= \int_{-1}^1 y(x) \,dz^*(x) = \langle z^*, y\rangle
\end{align*}
Hence we can conclude that $y^*\in D(B^*)$ and $B^* y^* = z^*$. \qed

Now we are in a position to find $Y^\odot$, the sun-dual of $Y$ with respect to $S$, which is the closure of $D(B^*)$ with respect to the total variation Norm. 
\begin{theorem}
The sun dual $Y^\odot$ with respect to the semigroup $S$ can be represented as $\mathbb{R}\times L^1(\Omega)$. For the sun dual of $B$ we have that 
\begin{equation}
D(B^{\odot}):=\{(c,w^\odot)\in \mathbb{R}\times L^1(\Omega)| c \in \mathbb{R}, (w^\odot)'\in AC([-1,1]) ,(w^\odot)'(1)=0\}
\end{equation}
and $B^\odot (c,w^\odot) := ((w^\odot)'(-1),(w^\odot)'')$, where $(w^\odot)''$ is some $L^1$ function such that 
\begin{equation}
(w^\odot)'(x)= (w^\odot)'(-1) + \int_{-1}^x (w^\odot)''(s)\,ds
\end{equation}
\end{theorem}
\textit{Proof.} Let $y^* \in D(B^*)$. Again using the notation that for $x,s\in(-1,1]$, 
\begin{align*}
y^*(x)&= c_1 + \int_{-1}^x w^*(s) \,ds\\
w^*(s)&=c_2 + \int_{-1}^s z^*(x') \,dx'
\end{align*}
for some $c_1,c_2 \in \mathbb{R}$ and $z^* \in NBV(\Omega)$ with $z^*(1)=0$, we can rewrite the total variation norm as:
\[||y^*||_{Y^*}=|c_1| + ||w^*||_{L^1}\]
For the space 
\[W := \left\{c + \int_{-1}^s z^*(x') \,dx' \middle| c \in \mathbb{R}, z^* \in NBV(\Omega), z^*(1)=0\right\}\]
we have that $\{w^* \in C^2 | (w^*)'(-1)=0\}\subset W \subset L^1$. As this first space of $C^2$ functions is dense in $L^1$, we have that $W$ is dense in $L^1$. Hence, we can represent $Y^\odot$ as the space 
\small
\[\left\lbrace y^\odot\in NBV(\Omega)\middle| y^\odot(x) = c + \int_{-1}^x w^\odot(s) \,ds \text{ where } c\in \mathbb{R}, w^\odot \in L^1(\Omega) \text{ for } x \in(-1,1]\right\rbrace\]
\normalsize
which are the absolutely continuous functions on $(-1,1]$ with a jump from $0$ to $c$ at $x=-1$. 

We can equivalently express $Y^\odot$ as $\mathbb{R}\times L^1(\Omega)$ where $y^\odot=(c,w^\odot)$ with $c\in \mathbb{R}$ and $w^\odot 
\in L^1(\Omega)$ equipped with the norm 
\[||y^\odot||_{Y^\odot}:=|c_1| + ||w^\odot||_{L^1}\]
The domain of $B^{\odot}$ is defined as $D(B^{\odot})=\{y^\odot\in D(B^*) | B^* y^\odot \in Y^{\odot}\}$. Using equation \eqref{eq:ystar} we have $B^*y^* = z^*$. If $z^* \in Y^\odot$ then $z^*$ must be absolutely continuous on $(-1,1]$. So for $y^\odot=(c,w^\odot)$ we find that $(w^\odot)'=z^*$ is absolutely continuous on $(-1,1]$. As $(w^\odot)'$ is an $L^1$-function, we can redefine $(w^\odot)'(-1):= (w^\odot)'(-1+)$ to get a absolutely continuous function on $[-1,1]$. The boundary condition $z(1)=0$ is transformed into $(w^\odot)'(1)=0$ 

Thus we can write that $B^\odot (c,w^\odot)= ((w^\odot)'(-1), (w^\odot)'')$, where $(w^\odot)''$ is an $L^1$ function such that 
\[(w^\odot)'(x)= (w^\odot)'(-1) + \int_{-1}^x (w^\odot)''(s)\,ds\]
\qed

Note that the sun-dual $Y^\odot$ is almost the same as in the book by \citet[Theorem II.5.2]{diekmann_delay_1995}, where it is taken with respect to the first derivative with the condition $\dot{y}(0)=0$. However, in that case there was an extra condition in $Y^\odot$ that functions $g \in L^1$ could be extended be zero for $\theta \geq h$. In our case with diffusion we have a fixed domain on which the diffusion takes place, so this condition is not present.

Now we can take the dual again and end up at the dual space $Y^{\odot *}$.
\begin{theorem}
The dual space $Y^{\odot *}$ can be represented as $\mathbb{R}\times L^\infty(\Omega)$. For the the operator $B^{\odot *}$ we have that 
\small
\begin{equation}\label{eq:bsunstar}
D(B^{\odot *})=\{(\gamma,w^{\odot *})| (w^{\odot *})' \text{ is Lipschitz continuous, } w^{\odot *}(-1)=\gamma, (w^{\odot *})'(\pm 1)=0 \}
\end{equation}
\normalsize
and $B^{\odot*} (\gamma,w^{\odot *}) := (0, (w^{\odot *})'')$, where $(w^{\odot *})''$ is an $L^\infty(\Omega)$ function such that 
\begin{equation}
(w^{\odot *})'(x)=\int_{-1}^x (w^{\odot *})''(s)\,ds
\end{equation}
\end{theorem}
\textit{Proof.} The dual space of $\mathbb{R}\times L^1(\Omega)$ can be represented as $\mathbb{R}\times L^\infty(\Omega)$ with the duality pairing between $Y^{\odot *}$ and $Y^\odot$ being given by
\[\langle(\gamma, w^{\odot *}),(c,w^\odot)\rangle := \gamma c + \int_{-1}^1 w^{\odot *}(x)w^{\odot }(x)\,dx\]

First we prove the $\subseteq$ inclusion of \eqref{eq:bsunstar}. Let $(\gamma,w^{\odot *})\in D(B^{\odot *})$ and $B^{\odot *}(\gamma,w^{\odot *})=(\beta,z^{\odot *})$. Let 
\[v^{\odot *}(x) :=v^{\odot *}(-1)+\int_{-1}^x z^{\odot *}(s)\,ds\]
which is a Lipschitz continuous function as $z^{\odot *} \in L^\infty(\Omega)$. 
Then for all $(c,w^\odot)\in D(B^\odot)$ we get that
\small
\begin{align*}
\gamma (w^\odot)'(-1) + \int_\Omega w^{\odot *}(x)(w^\odot)''(x)\,dx&=\langle (\gamma,w^{\odot *}),B^\odot (c,w^\odot)\rangle =\langle(\beta,z^{\odot *}),(c,w^\odot)\rangle\\
&=\beta c + \int_\Omega z^{\odot *}(x)w^\odot(x)\,dx\\
&=\beta c + v^{\odot *}(x)w^\odot(x)|_{-1}^1 \\
&\qquad \;\; - \int_\Omega v^{\odot *}(x)(w^\odot)'(x)\,dx\\
&=\beta c + v^{\odot *}(-1)w^\odot(x)|_{-1}^1 + \gamma(w^\odot)'(-1) \\
&\qquad \;\; + \int_\Omega \left(\gamma +\int_{-1}^x v^{\odot *}(s)\,ds \right) (w^\odot)''(x)\,dx
\end{align*}
\normalsize
Here we used that $(w^\odot)' \in AC[-1,1]$ and $(w^\odot)'(1)=0$. As $c$ and $w^\odot(\pm 1)$ are arbitrary we see that necessarily $\beta=0, v^{\odot *}(\pm 1)=0$. Furthermore,
\[w^{\odot *}(x)=\gamma+\int_{-1}^x v^{\odot *}(s)\,ds\]
which implies that $(w^{\odot *})'=v^{\odot *}$ and $w^{\odot *}(-1)=\gamma$. 

Finally we prove the $\supseteq$ inclusion of \eqref{eq:bsunstar}. Let $(\gamma,w^{\odot *})$ be in the righthand side of \eqref{eq:bsunstar} and $(c,w^\odot)\in D(B^\odot)$. Then by the calculations above we get that 
\begin{align*}
\langle (\gamma,w^{\odot *}),B^\odot (c,w^\odot)\rangle &= \gamma (w^\odot)'(-1) + \int_\Omega w^{\odot *}(x)(w^\odot)''(x)\,dx\\
&= \int_\Omega (w^{\odot *})''(x)w^\odot(x)\,dx\\
&= \langle (0,(w^{\odot *})''), (c,w^\odot)\rangle
\end{align*}
From which we can conclude that $(\gamma,w^{\odot *}) \in D(B^{\odot *}$ and $B^{\odot*} (\gamma,w^{\odot *}) := (0, (w^{\odot *})'')$.
\qed

Finally we characterise the sun bi-dual $Y^{\odot \odot}$ which is the closure of $D(B^{\odot *})$ with respect to the $Y^{\odot *}$-norm, which is a supremum norm.

\begin{theorem}
The sun bi-dual $Y^{\odot \odot}$ can be represented as $\{(\gamma,w^{\odot \odot})|w^{\odot \odot}\in C(\Omega), w^{\odot \odot}(-1)=\gamma\}$. The canonical embedding $j_Y:Y\rightarrow Y^{\odot *}$ is given by $j_Y y=(y(-1),y)$. Moreover, $Y$ is sun-reflexive with respect to the semigroup $S$, i.e. $j_Y(Y)=Y^{\odot \odot}$. 
\end{theorem}
\textit{Proof.} Let $y^{\odot *} = (\gamma, w^{\odot *}) \in Y^{\odot *}$. As the supremum norm does not preserve derivatives, i.e. the $C^2$ functions are dense in $C^0$ with respect to the supremum norm, we have that only the continuity and the condition $w^{\odot *}(-1)=\gamma$ remain. For $j_Y y=(y(-1),y)$, it can be easily checked that for any $y^\odot \in Y^\odot$
\[\langle j_Y y, y^\odot \rangle = \langle y^\odot, y\rangle = \int_\Omega \left(\gamma +\int_{-1}^x v^{\odot *}(s)\,ds \right) (w^\odot)''(x)\,dx = \]
So the $j_Y$ is the canonical embedding between $Y$ and $Y^{\odot *}$ and it is an isomorphism between $Y$ and $Y^{\odot \odot}$. Hence $Y$ is sun-reflexive. \qed

\section{Proofs}
\label{sec:proofs}
\begin{lemma}\label{lem:int_by_parts}
Let $\Phi, \psi \in L^{\infty}([-h,0];Y^{**})$ and $g,\dot{g} \in L^1([0,h];Y^*)$ such that
\begin{align*}
\Phi(-t) &= \Phi(0) - \int_0^t \psi(-\theta) \,d\theta \\
g(t) &= g(0) + \int_0^t \dot{g}(\theta)\,d\theta
\end{align*}
for all $t\in [0,h]$, then it holds that 
\[\langle \Phi(-t),g(t) \rangle = \langle \Phi(0),g(0) \rangle + \int_0^t \langle \Phi(-\theta),\dot{g}(\theta) \rangle \,d\theta - \int_0^t \langle \psi(-\theta),g(\theta) \rangle \,d\theta\]
for all $t \in [0,h]$.
\end{lemma}
\textit{Proof.} Let $\Phi,\psi, g, \dot{g}$ as above and define the scalar function $\xi$
\[\xi(t) := \langle \Phi(-t), g(t) \rangle\]
for $t \in [0,h]$. As $\Phi \in L^{\infty}([-h,0];Y^{**})$ and $g \in L^1([0,h];Y^*)$, $\xi$ is integrable.

By definition $\xi$ is absolutely continuous on an interval $I$ if for every $\epsilon>0$, there is a $\delta>0$ such that whenever a finite sequence of pairwise disjoint sub-intervals $(s_k,t_k)$ of $I$ with $t_k,s_k\in I$ satisfies
\[\sum_k (t_k - s_k) < \delta\]
then
\[\sum_k \|\xi(t_k) - \xi(s_k)\| < \epsilon\]
Both $\Phi$ and $g$ are absolutely continuous and a.e. differentiable with derivative $\psi$ and $\dot{g}$ respectively \citep[Corollary 1.4.31]{cazenave_introduction_1998}.

For $t,s\in [0,h]$
\begin{align*}
|\xi(t)-\xi(s)| &= |\langle \Phi(-t), g(t) \rangle - \langle \Phi(-s), g(s) \rangle|\\
&= |\langle \Phi(-t) - \Phi(-s), g(t) \rangle + \langle \Phi(-s), g(t)-g(s) \rangle|\\
&\leq \|\Phi(-t) - \Phi(-s)\| \max_{t\in [0,h]}\|g(t)\| + \|g(t) - g(s)\| \max_{t\in [0,h]}\|\Phi(-t)\|
\end{align*}
Hence by the absolute continuity of $\Phi$ and $g$, $\xi$ is absolutely continuous and consequently has an a.e. derivative $\dot{\xi}$, which is integrable and for $t\in [0,h]$
\[\xi(t) = \xi(0) + \int_0^t \dot{\xi}(\theta) \,d\theta\]

Furthermore, we have that
\[\frac{\xi(t)-\xi(s)}{t-s} = \left\langle \Phi(-s), \frac{g(t)-g(s)}{t-s} \right\rangle - \left\langle \frac{\Phi(-t) - \Phi(-s)}{s-t}, g(t) \right\rangle\]
Taking the limit as $s \rightarrow t$ we can deduce that 
\[\dot{\xi}(t) = \langle \Phi(-t), \dot{g}(t) \rangle - \langle \psi(-t), g(t) \rangle\]
Hence we have that for $t\in [0,h]$
\[\langle \Phi(-t),g(t) \rangle = \langle \Phi(0),g(0) \rangle + \int_0^t \langle \Phi(-\theta),\dot{g}(\theta) \rangle \,d\theta - \int_0^t \langle \psi(-\theta),g(\theta) \rangle \,d\theta\] \qed

\begin{lemma}\label{lem:Q_invert}
Define the matrix $\hat{Q} \in C^{(N+1)\times (N+1)}$ as 
\[\hat{Q}_{j,m}=\begin{cases}
\frac{1}{n_j-p_m} &\quad\text{for } j\in \{1,\cdots,N\}, m\in \{1,\cdots,N+1\}\\
1 &\quad\text{for } j=N+1, m\in \{1,\cdots,N+1\}\end{cases}\]
When $n_i\neq n_j\neq p_m \neq p_l$ for $i,j\in \{1,\cdots,N\}$, $l,m\in \{1,\cdots,N+1\}$, $i\neq j$, $l\neq m$, then $\hat{Q}$ is invertible.
\end{lemma}
\textit{Proof.}
We subtract the last column from the other columns. We get the following matrix $\tilde{Q}$
\[\tilde{Q}_{j,m}=\begin{cases}
\frac{p_m-p_{N+1}}{(n_j-p_m)(n_j-p_{N+1})} &\quad\text{for } j,m\in \{1,\cdots,N\}\\
\frac{1}{n_j-p_{N+1}} &\quad\text{for } j\in \{1,\cdots,N\}, m=N+1\\
0 &\quad\text{for } j=N+1, m\in \{1,\cdots,N\}\\
1 &\quad\text{for } j=m=N+1 \end{cases}\]
Now row $j$ of matrix $\tilde{Q}$ contains the factor $\frac{1}{n_j-p_{N+1}}$ and column $m$ contains the factor $p_m-p_{N+1}$ for $j,m \in \{1,\cdots,N\}$. Hence we can rewrite the determinant of $\hat{Q}$ as 
\[\det(\hat{Q})=\det(\tilde{Q})=\det(Q)\prod_{i=1}^N\frac{p_i-p_{N+1}}{n_i-p_{N+1}}\]
Here matrix $Q \in C^{N\times N}$ is defined as
\[Q_{j,m}=\frac{1}{n_j-p_m} \quad\text{for } j,m\in \{1,\cdots,N\}\]
We observe that $Q$ is a Cauchy matrix when $n_i\neq n_j\neq p_m \neq p_l$ for $i,j,l,m\in \{1,\cdots,N\}$, $i\neq j$, $l\neq m$ and hence invertible. Furthermore the product $\prod_{i=1}^N\frac{p_i-p_{N+1}}{n_i-p_{N+1}}$ is non-zero, so we conclude that $\hat{Q}$ is invertible. \qed

\section*{Acknowledgements}
We want to thank Sebastiaan Janssens for useful discussions on the sun-star calculus and his inspiring preprint on ArXiv \cite{janssens_class_2019} on abstract delay differential equations, to which our model belongs. We acknowledge R. Bellingacci for numerical analysis of neural field model in his Master Thesis at Utrecht University.
\end{appendix}

\bibliographystyle{plainnat}
\bibliography{NeuralField} 

\begin{thebibliography}{67}
\providecommand{\natexlab}[1]{#1}
\providecommand{\url}[1]{\texttt{#1}}
\expandafter\ifx\csname urlstyle\endcsname\relax
  \providecommand{\doi}[1]{doi: #1}\else
  \providecommand{\doi}{doi: \begingroup \urlstyle{rm}\Url}\fi

\bibitem[Amari(1977)]{amari_dynamics_1977}
Shun-ichi Amari.
\newblock Dynamics of pattern formation in lateral-inhibition type neural
  fields.
\newblock \emph{Biological Cybernetics}, 27\penalty0 (2):\penalty0 77--87, June
  1977.
\newblock ISSN 1432-0770.
\newblock \doi{10.1007/BF00337259}.
\newblock URL \url{https://doi.org/10.1007/BF00337259}.

\bibitem[Amitai et~al.(2002)Amitai, Gibson, Beierlein, Patrick, Ho, Connors,
  and Golomb]{amitai_spatial_2002}
Yael Amitai, Jay~R. Gibson, Michael Beierlein, Saundra~L. Patrick, Alice~M. Ho,
  Barry~W. Connors, and David Golomb.
\newblock The {Spatial} {Dimensions} of {Electrically} {Coupled} {Networks} of
  {Interneurons} in the {Neocortex}.
\newblock \emph{Journal of Neuroscience}, 22\penalty0 (10):\penalty0
  4142--4152, May 2002.
\newblock ISSN 0270-6474, 1529-2401.
\newblock \doi{10.1523/JNEUROSCI.22-10-04142.2002}.
\newblock URL \url{http://www.jneurosci.org/content/22/10/4142}.

\bibitem[Bartle(1956)]{bartle_general_1956}
R.~Bartle.
\newblock A general bilinear vector integral.
\newblock \emph{Studia Mathematica}, 15\penalty0 (3):\penalty0 337--352, 1956.
\newblock ISSN 0039-3223.
\newblock URL \url{https://eudml.org/doc/216873}.

\bibitem[Bartle(2001)]{bartle_modern_2001}
Robert~Gardner Bartle.
\newblock \emph{A {Modern} {Theory} of {Integration}}.
\newblock American Mathematical Soc., 2001.
\newblock ISBN 978-0-8218-0845-0.

\bibitem[Batkai and Piazzera(2005)]{batkai_semigroups_2005}
Andras Batkai and Susanna Piazzera.
\newblock \emph{Semigroups for {Delay} {Equations}}.
\newblock CRC Press, September 2005.
\newblock ISBN 978-1-4398-6568-2.

\bibitem[Batkai and Piazzera(2001)]{batkai_semigroups_2001}
András Batkai and Susanna Piazzera.
\newblock Semigroups and {Linear} {Partial} {Differential} {Equations} with
  {Delay}.
\newblock \emph{Journal of Mathematical Analysis and Applications},
  264\penalty0 (1):\penalty0 1--20, December 2001.
\newblock ISSN 0022-247X.
\newblock \doi{10.1006/jmaa.2001.6705}.
\newblock URL
  \url{http://www.sciencedirect.com/science/article/pii/S0022247X01967058}.

\bibitem[Campbell(2007)]{jirsa_time_2007}
Sue~Ann Campbell.
\newblock Time {Delays} in {Neural} {Systems}.
\newblock In Viktor~K Jirsa and Ar~McIntosh, editors, \emph{Handbook of {Brain}
  {Connectivity}}, pages 65--90. Springer Berlin Heidelberg, Berlin,
  Heidelberg, 2007.
\newblock ISBN 978-3-540-71462-0 978-3-540-71512-2.
\newblock \doi{10.1007/978-3-540-71512-2_2}.
\newblock URL \url{http://link.springer.com/10.1007/978-3-540-71512-2_2}.

\bibitem[Cazenave et~al.(1998)Cazenave, Braides, Haraux, and
  Haraux]{cazenave_introduction_1998}
Thierry Cazenave, Andrea Braides, Alain Haraux, and Both Professors of
  Mathematics~Alain Haraux.
\newblock \emph{An {Introduction} to {Semilinear} {Evolution} {Equations}}.
\newblock Clarendon Press, 1998.
\newblock ISBN 978-0-19-850277-7.

\bibitem[Clément et~al.(1986)Clément, Diekmann, Gyllenberg, Heijmans, and
  Thieme]{clement_perturbation_1986}
Philippe Clément, Odo Diekmann, M.~Gyllenberg, Henk Heijmans, and H.~R.
  Thieme.
\newblock Perturbation theory for dual semigroups. {IV}. {The} intertwining
  formula and the canonical pairing.
\newblock \emph{Department of Applied Mathematics}, January 1986.
\newblock URL \url{https://ir.cwi.nl/pub/12515}.

\bibitem[Coombes(2005)]{coombes_waves_2005}
S.~Coombes.
\newblock Waves, bumps, and patterns in neural field theories.
\newblock \emph{Biological Cybernetics}, 93\penalty0 (2):\penalty0 91--108,
  August 2005.
\newblock ISSN 1432-0770.
\newblock \doi{10.1007/s00422-005-0574-y}.
\newblock URL \url{https://doi.org/10.1007/s00422-005-0574-y}.

\bibitem[Coombes(2010)]{coombes_large-scale_2010}
S.~Coombes.
\newblock Large-scale neural dynamics: {Simple} and complex.
\newblock \emph{NeuroImage}, 52\penalty0 (3):\penalty0 731--739, September
  2010.
\newblock ISSN 1053-8119.
\newblock \doi{10.1016/j.neuroimage.2010.01.045}.
\newblock URL
  \url{http://www.sciencedirect.com/science/article/pii/S1053811910000674}.

\bibitem[Coombes et~al.(2014)Coombes, beim Graben, and
  Potthast]{coombes_tutorial_2014}
Stephen Coombes, Peter beim Graben, and Roland Potthast.
\newblock \emph{Tutorial on neural field theory}.
\newblock Springer, 2014.

\bibitem[{Coombes Stephen} and {Laing
  Carlo}(2009)]{coombes_stephen_delays_2009}
{Coombes Stephen} and {Laing Carlo}.
\newblock Delays in activity-based neural networks.
\newblock \emph{Philosophical Transactions of the Royal Society A:
  Mathematical, Physical and Engineering Sciences}, 367\penalty0
  (1891):\penalty0 1117--1129, March 2009.
\newblock \doi{10.1098/rsta.2008.0256}.
\newblock URL
  \url{https://royalsocietypublishing.org/doi/full/10.1098/rsta.2008.0256}.

\bibitem[Coullet and Spiegel(1983)]{coullet_amplitude_1983}
P.~H. Coullet and E.~A. Spiegel.
\newblock Amplitude {Equations} for {Systems} with {Competing} {Instabilities}.
\newblock \emph{SIAM Journal on Applied Mathematics}, 43\penalty0 (4):\penalty0
  776--821, August 1983.
\newblock ISSN 0036-1399.
\newblock \doi{10.1137/0143052}.
\newblock URL \url{https://epubs.siam.org/doi/abs/10.1137/0143052}.
\newblock Publisher: Society for Industrial and Applied Mathematics.

\bibitem[Diekmann et~al.(1995)Diekmann, Gils, Lunel, and
  Walther]{diekmann_delay_1995}
Odo Diekmann, Stephan A.~van Gils, Sjoerd M.~V. Lunel, and Hans-Otto Walther.
\newblock \emph{Delay {Equations}: {Functional}-, {Complex}-, and {Nonlinear}
  {Analysis}}.
\newblock Springer Science \& Business Media, August 1995.
\newblock ISBN 978-1-4612-4206-2.

\bibitem[Dijkstra et~al.(2015)Dijkstra, Gils, Janssens, Kuznetsov, and
  Visser]{dijkstra_pitchforkhopf_2015}
K.~Dijkstra, S.A.~van Gils, S.G. Janssens, Yu.A. Kuznetsov, and S.~Visser.
\newblock Pitchfork–{Hopf} bifurcations in {1D} neural field models with
  transmission delays.
\newblock \emph{Physica D: Nonlinear Phenomena}, 297:\penalty0 88--101, March
  2015.
\newblock ISSN 01672789.
\newblock \doi{10.1016/j.physd.2015.01.004}.
\newblock URL
  \url{https://linkinghub.elsevier.com/retrieve/pii/S0167278915000111}.

\bibitem[Elphick et~al.(1987)Elphick, Tirapegui, Brachet, Coullet, and
  Iooss]{elphick_simple_1987}
C.~Elphick, E.~Tirapegui, M.~E. Brachet, P.~Coullet, and G.~Iooss.
\newblock A simple global characterization for normal forms of singular vector
  fields.
\newblock \emph{Physica D: Nonlinear Phenomena}, 29\penalty0 (1):\penalty0
  95--127, November 1987.
\newblock ISSN 0167-2789.
\newblock \doi{10.1016/0167-2789(87)90049-2}.
\newblock URL
  \url{http://www.sciencedirect.com/science/article/pii/0167278987900492}.

\bibitem[Engel and Nagel(1999)]{engel_one-parameter_1999}
Klaus-Jochen Engel and Rainer Nagel.
\newblock \emph{One-{Parameter} {Semigroups} for {Linear} {Evolution}
  {Equations}}, volume~63.
\newblock Springer, 1999.

\bibitem[Engelborghs et~al.(2002)Engelborghs, Luzyanina, and
  Roose]{engelborghs_numerical_2002}
K.~Engelborghs, T.~Luzyanina, and D.~Roose.
\newblock Numerical bifurcation analysis of delay differential equations using
  {DDE}-{BIFTOOL}.
\newblock \emph{ACM Transactions on Mathematical Software}, 28\penalty0
  (1):\penalty0 1--21, March 2002.
\newblock ISSN 0098-3500.
\newblock \doi{10.1145/513001.513002}.
\newblock URL \url{https://doi.org/10.1145/513001.513002}.

\bibitem[Ermentrout and Cowan(1980)]{ermentrout_large_1980}
G.~Ermentrout and J.~Cowan.
\newblock Large {Scale} {Spatially} {Organized} {Activity} in {Neural} {Nets}.
\newblock \emph{SIAM Journal on Applied Mathematics}, 38\penalty0 (1):\penalty0
  1--21, February 1980.
\newblock ISSN 0036-1399.
\newblock \doi{10.1137/0138001}.
\newblock URL \url{https://epubs.siam.org/doi/abs/10.1137/0138001}.

\bibitem[Ermentrout and Terman(2010)]{ermentrout_mathematical_2010}
G.~Bard Ermentrout and David~H. Terman.
\newblock \emph{Mathematical {Foundations} of {Neuroscience}}.
\newblock Springer Science \& Business Media, July 2010.
\newblock ISBN 978-0-387-87708-2.

\bibitem[Faria(2006)]{arino_normal_2006}
T.~Faria.
\newblock Normal {Forms} and {Bifurcations} {For} {Delay} {Differential}
  {Equations}.
\newblock In O.~Arino, M.L. Hbid, and E.~Ait Dads, editors, \emph{Delay
  {Differential} {Equations} and {Applications}}, volume 205, pages 227--282.
  Springer Netherlands, Dordrecht, 2006.
\newblock ISBN 978-1-4020-3645-3.
\newblock \doi{10.1007/1-4020-3647-7_7}.
\newblock URL \url{http://link.springer.com/10.1007/1-4020-3647-7_7}.

\bibitem[Faria and Magalhaes(1995{\natexlab{a}})]{faria_normal_1995-1}
T.~Faria and L.~T. Magalhaes.
\newblock Normal {Forms} for {Retarded} {Functional} {Differential} {Equations}
  with {Parameters} and {Applications} to {Hopf} {Bifurcation}.
\newblock \emph{Journal of Differential Equations}, 122\penalty0 (2):\penalty0
  181--200, November 1995{\natexlab{a}}.
\newblock ISSN 0022-0396.
\newblock \doi{10.1006/jdeq.1995.1144}.
\newblock URL
  \url{http://www.sciencedirect.com/science/article/pii/S0022039685711448}.

\bibitem[Faria and Magalhaes(1995{\natexlab{b}})]{faria_normal_1995}
T.~Faria and L.T. Magalhaes.
\newblock Normal {Forms} for {Retarded} {Functional} {Differential} {Equations}
  and {Applications} to {Bogdanov}-{Takens} {Singularity}.
\newblock \emph{Journal of Differential Equations}, 122\penalty0 (2):\penalty0
  201--224, November 1995{\natexlab{b}}.
\newblock ISSN 00220396.
\newblock \doi{10.1006/jdeq.1995.1145}.
\newblock URL
  \url{https://linkinghub.elsevier.com/retrieve/pii/S002203968571145X}.

\bibitem[Faye and Faugeras(2010)]{faye_theoretical_2010}
Grégory Faye and Olivier Faugeras.
\newblock Some theoretical and numerical results for delayed neural field
  equations.
\newblock \emph{Physica D: Nonlinear Phenomena}, 239\penalty0 (9):\penalty0
  561--578, May 2010.
\newblock ISSN 0167-2789.
\newblock \doi{10.1016/j.physd.2010.01.010}.
\newblock URL
  \url{http://www.sciencedirect.com/science/article/pii/S0167278910000229}.

\bibitem[Gowurin(1936)]{gowurin_uber_1936}
Mark Gowurin.
\newblock Über die {Stieltjessche} {Integration} abstrakter {Funktionen}.
\newblock \emph{Fundamenta Mathematicae}, 27:\penalty0 254--265, 1936.
\newblock ISSN 0016-2736, 1730-6329.
\newblock \doi{10.4064/fm-27-1-254-265}.
\newblock URL
  \url{https://www.impan.pl/en/publishing-house/journals-and-series/fundamenta-mathematicae/all/27/0/93340/uber-die-stieltjessche-integration-abstrakter-funktionen}.

\bibitem[Hale(1971)]{hale_theory_1971}
Jack~K. Hale.
\newblock \emph{Theory of {Functional} {Differential} {Equations}}.
\newblock Springer Science \& Business Media, December 1971.
\newblock ISBN 978-1-4612-9892-2.

\bibitem[Hodgkin and Huxley(1952)]{hodgkin_quantitative_1952}
A.~L. Hodgkin and A.~F. Huxley.
\newblock A quantitative description of membrane current and its application to
  conduction and excitation in nerve.
\newblock \emph{The Journal of Physiology}, 117\penalty0 (4):\penalty0
  500--544, 1952.
\newblock ISSN 1469-7793.
\newblock \doi{10.1113/jphysiol.1952.sp004764}.
\newblock URL
  \url{https://physoc.onlinelibrary.wiley.com/doi/abs/10.1113/jphysiol.1952.sp004764}.

\bibitem[Hutt(2008)]{hutt_local_2008}
A.~Hutt.
\newblock Local excitation-lateral inhibition interaction yields oscillatory
  instabilities in nonlocally interacting systems involving finite propagation
  delay.
\newblock \emph{Physics Letters A}, 372\penalty0 (5):\penalty0 541--546,
  January 2008.
\newblock ISSN 0375-9601.
\newblock \doi{10.1016/j.physleta.2007.08.018}.
\newblock URL
  \url{http://www.sciencedirect.com/science/article/pii/S0375960107011681}.

\bibitem[Hutt and Atay(2005)]{hutt_analysis_2005}
Axel Hutt and Fatihcan~M. Atay.
\newblock Analysis of nonlocal neural fields for both general and
  gamma-distributed connectivities.
\newblock \emph{Physica D: Nonlinear Phenomena}, 203\penalty0 (1):\penalty0
  30--54, April 2005.
\newblock ISSN 0167-2789.
\newblock \doi{10.1016/j.physd.2005.03.002}.
\newblock URL
  \url{http://www.sciencedirect.com/science/article/pii/S0167278905000989}.

\bibitem[Hutt and Atay(2007)]{hutt_spontaneous_2007}
Axel Hutt and Fatihcan~M. Atay.
\newblock Spontaneous and evoked activity in extended neural populations with
  gamma-distributed spatial interactions and transmission delay.
\newblock \emph{Chaos, Solitons \& Fractals}, 32\penalty0 (2):\penalty0
  547--560, April 2007.
\newblock ISSN 0960-0779.
\newblock \doi{10.1016/j.chaos.2005.10.091}.
\newblock URL
  \url{http://www.sciencedirect.com/science/article/pii/S0960077905010817}.

\bibitem[Hutt et~al.(2003)Hutt, Bestehorn, and Wennekers]{hutt_pattern_2003}
Axel Hutt, Michael Bestehorn, and Thomas Wennekers.
\newblock Pattern formation in intracortical neuronal fields.
\newblock \emph{Network: Computation in Neural Systems}, 14\penalty0
  (2):\penalty0 351--368, January 2003.
\newblock ISSN 0954-898X.
\newblock \doi{10.1088/0954-898X_14_2_310}.
\newblock URL \url{https://doi.org/10.1088/0954-898X_14_2_310}.
\newblock Publisher: Taylor \& Francis \_eprint:
  https://doi.org/10.1088/0954-898X\_14\_2\_310.

\bibitem[Janssens(2010)]{janssens_normalization_2010}
Sebastiaan~G Janssens.
\newblock On a {Normalization} {Technique} for {Codimension} {Two}
  {Bifurcations} of {Equilibria} of {Delay} {Differential} {Equations}.
\newblock \emph{Master Thesis at University of Utrecht}, November 2010.
\newblock URL \url{https://sebastiaanjanssens.nl/pdf/normalization.pdf}.

\bibitem[Janssens(2019)]{janssens_class_2019}
Sebastiaan~G. Janssens.
\newblock A class of abstract delay differential equations in the light of suns
  and stars.
\newblock \emph{arXiv:1901.11526 [math]}, January 2019.
\newblock URL \url{http://arxiv.org/abs/1901.11526}.
\newblock arXiv: 1901.11526.

\bibitem[Janssens(2020)]{janssens_class_2020}
Sebastiaan~G. Janssens.
\newblock A class of abstract delay differential equations in the light of suns
  and stars. {II}.
\newblock \emph{arXiv:2003.13341 [math]}, March 2020.
\newblock URL \url{http://arxiv.org/abs/2003.13341}.
\newblock arXiv: 2003.13341.

\bibitem[Jirsa and Haken(1996)]{jirsa_field_1996}
V.~K. Jirsa and H.~Haken.
\newblock Field {Theory} of {Electromagnetic} {Brain} {Activity}.
\newblock \emph{Physical Review Letters}, 77\penalty0 (5):\penalty0 960--963,
  July 1996.
\newblock \doi{10.1103/PhysRevLett.77.960}.
\newblock URL \url{https://link.aps.org/doi/10.1103/PhysRevLett.77.960}.
\newblock Publisher: American Physical Society.

\bibitem[Jirsa et~al.(2002)Jirsa, Jantzen, Fuchs, and
  Kelso]{jirsa_spatiotemporal_2002}
V.K. Jirsa, K.J. Jantzen, A.~Fuchs, and J.A.S. Kelso.
\newblock Spatiotemporal forward solution of the {EEG} and {MEG} using network
  modeling.
\newblock \emph{IEEE Transactions on Medical Imaging}, 21\penalty0
  (5):\penalty0 493--504, May 2002.
\newblock ISSN 1558-254X.
\newblock \doi{10.1109/TMI.2002.1009385}.
\newblock Conference Name: IEEE Transactions on Medical Imaging.

\bibitem[Katō(1995)]{kato_perturbation_1995}
Tosio Katō.
\newblock \emph{Perturbation theory for linear operators}.
\newblock Classics in mathematics. Springer, Berlin, 1995.
\newblock ISBN 978-3-540-58661-6.

\bibitem[Kuznetsov(2004)]{kuznetsov_elements_2004}
Yuri~A. Kuznetsov.
\newblock \emph{Elements of {Applied} {Bifurcation} {Theory}}.
\newblock Springer Science \& Business Media, 2004.
\newblock ISBN 978-1-4757-2421-9.

\bibitem[Laing(2015)]{laing_exact_2015}
C.~Laing.
\newblock Exact {Neural} {Fields} {Incorporating} {Gap} {Junctions}.
\newblock \emph{SIAM Journal on Applied Dynamical Systems}, 14\penalty0
  (4):\penalty0 1899--1929, January 2015.
\newblock \doi{10.1137/15M1011287}.
\newblock URL \url{https://epubs.siam.org/doi/abs/10.1137/15M1011287}.

\bibitem[Liley et~al.(2002)Liley, Cadusch, and Dafilis]{liley_spatially_2002}
D.~T.~J. Liley, P.~J. Cadusch, and M.~P. Dafilis.
\newblock A spatially continuous mean field theory of electrocortical activity.
\newblock \emph{Network: Computation in Neural Systems}, 13\penalty0
  (1):\penalty0 67--113, January 2002.
\newblock ISSN 0954-898X.
\newblock \doi{10.1080/net.13.1.67.113}.
\newblock URL \url{https://doi.org/10.1080/net.13.1.67.113}.
\newblock Publisher: Taylor \& Francis \_eprint:
  https://doi.org/10.1080/net.13.1.67.113.

\bibitem[Liu et~al.(2008)Liu, Magal, and Ruan]{liu_projectors_2008}
Zhihua Liu, Pierre Magal, and Shigui Ruan.
\newblock Projectors on the generalized eigenspaces for functional differential
  equations using integrated semigroups.
\newblock \emph{Journal of Differential Equations}, 244\penalty0 (7):\penalty0
  1784--1809, April 2008.
\newblock ISSN 0022-0396.
\newblock \doi{10.1016/j.jde.2008.01.007}.
\newblock URL
  \url{http://www.sciencedirect.com/science/article/pii/S0022039608000120}.

\bibitem[Liu et~al.(2014)Liu, Magal, and Ruan]{liu_normal_2014}
Zhihua Liu, Pierre Magal, and Shigui Ruan.
\newblock Normal forms for semilinear equations with non-dense domain with
  applications to age structured models.
\newblock \emph{Journal of Differential Equations}, 257\penalty0 (4):\penalty0
  921--1011, August 2014.
\newblock ISSN 0022-0396.
\newblock \doi{10.1016/j.jde.2014.04.018}.
\newblock URL
  \url{http://www.sciencedirect.com/science/article/pii/S0022039614001697}.

\bibitem[Magal and Ruan(2009{\natexlab{a}})]{magal_center_2009}
Pierre Magal and Shigui Ruan.
\newblock Center manifolds for semilinear equations with non-dense domain and
  applications to {Hopf} bifurcation in age structured models.
\newblock \emph{Memoirs of the American Mathematical Society}, 202\penalty0
  (951):\penalty0 0--0, 2009{\natexlab{a}}.
\newblock ISSN 0065-9266, 1947-6221.
\newblock \doi{10.1090/S0065-9266-09-00568-7}.
\newblock URL \url{http://www.ams.org/memo/0951}.

\bibitem[Magal and Ruan(2009{\natexlab{b}})]{magal_semilinear_2009}
Pierre Magal and Shigui Ruan.
\newblock On semilinear {Cauchy} problems with non-dense domain.
\newblock \emph{Advances in Differential Equations}, 14\penalty0
  (11/12):\penalty0 1041--1084, November 2009{\natexlab{b}}.
\newblock ISSN 1079-9389.
\newblock URL \url{https://projecteuclid.org/euclid.ade/1355854784}.
\newblock Publisher: Khayyam Publishing, Inc.

\bibitem[Nunez(1974)]{nunez_brain_1974}
Paul~L. Nunez.
\newblock The brain wave equation: a model for the {EEG}.
\newblock \emph{Mathematical Biosciences}, 21\penalty0 (3):\penalty0 279--297,
  December 1974.
\newblock ISSN 0025-5564.
\newblock \doi{10.1016/0025-5564(74)90020-0}.
\newblock URL
  \url{http://www.sciencedirect.com/science/article/pii/0025556474900200}.

\bibitem[Ostojic et~al.(2009)Ostojic, Brunel, and
  Hakim]{ostojic_synchronization_2009}
Srdjan Ostojic, Nicolas Brunel, and Vincent Hakim.
\newblock Synchronization properties of networks of electrically coupled
  neurons in the presence of noise and heterogeneities.
\newblock \emph{Journal of Computational Neuroscience}, 26\penalty0
  (3):\penalty0 369--392, June 2009.
\newblock ISSN 0929-5313, 1573-6873.
\newblock \doi{10.1007/s10827-008-0117-3}.
\newblock URL \url{http://link.springer.com/10.1007/s10827-008-0117-3}.

\bibitem[Roxin and Montbrió(2011)]{roxin_how_2011}
Alex Roxin and Ernest Montbrió.
\newblock How effective delays shape oscillatory dynamics in neuronal networks.
\newblock \emph{Physica D: Nonlinear Phenomena}, 240\penalty0 (3):\penalty0
  323--345, February 2011.
\newblock ISSN 0167-2789.
\newblock \doi{10.1016/j.physd.2010.09.009}.
\newblock URL
  \url{http://www.sciencedirect.com/science/article/pii/S0167278910002599}.

\bibitem[Roxin et~al.(2005)Roxin, Brunel, and Hansel]{roxin_role_2005}
Alex Roxin, Nicolas Brunel, and David Hansel.
\newblock Role of {Delays} in {Shaping} {Spatiotemporal} {Dynamics} of
  {Neuronal} {Activity} in {Large} {Networks}.
\newblock \emph{Physical Review Letters}, 94\penalty0 (23):\penalty0 238103,
  June 2005.
\newblock \doi{10.1103/PhysRevLett.94.238103}.
\newblock URL \url{https://link.aps.org/doi/10.1103/PhysRevLett.94.238103}.
\newblock Publisher: American Physical Society.

\bibitem[Roxin et~al.(2006)Roxin, Brunel, and Hansel]{roxin_rate_2006}
Alex Roxin, Nicolas Brunel, and David Hansel.
\newblock Rate {Models} with {Delays} and the {Dynamics} of {Large} {Networks}
  of {Spiking} {Neurons}.
\newblock \emph{Progress of Theoretical Physics Supplement}, 161:\penalty0
  68--85, January 2006.
\newblock ISSN 0375-9687.
\newblock \doi{10.1143/PTPS.161.68}.
\newblock URL
  \url{https://academic.oup.com/ptps/article/doi/10.1143/PTPS.161.68/1900315}.
\newblock Publisher: Oxford Academic.

\bibitem[Sanz~Leon et~al.(2013)Sanz~Leon, Knock, Woodman, Domide, Mersmann,
  McIntosh, and Jirsa]{sanz_leon_virtual_2013}
Paula Sanz~Leon, Stuart~A. Knock, M.~Marmaduke Woodman, Lia Domide, Jochen
  Mersmann, Anthony~R. McIntosh, and Viktor Jirsa.
\newblock The {Virtual} {Brain}: a simulator of primate brain network dynamics.
\newblock \emph{Frontiers in Neuroinformatics}, 7, 2013.
\newblock ISSN 1662-5196.
\newblock \doi{10.3389/fninf.2013.00010}.
\newblock URL
  \url{https://www.frontiersin.org/articles/10.3389/fninf.2013.00010/full#h1}.

\bibitem[Schwab et~al.(2014{\natexlab{a}})Schwab, Heida, Zhao, Gils, and
  Wezel]{schwab_pallidal_2014}
Bettina~C. Schwab, Tjitske Heida, Yan Zhao, Stephan A.~van Gils, and Richard J.
  A.~van Wezel.
\newblock Pallidal gap junctions-triggers of synchrony in {Parkinson}'s
  disease?
\newblock \emph{Movement Disorders}, 29\penalty0 (12):\penalty0 1486--1494,
  2014{\natexlab{a}}.
\newblock ISSN 1531-8257.
\newblock \doi{10.1002/mds.25987}.
\newblock URL \url{https://onlinelibrary.wiley.com/doi/abs/10.1002/mds.25987}.

\bibitem[Schwab et~al.(2014{\natexlab{b}})Schwab, Meijer, van Wezel, and van
  Gils]{schwab_synchronization_2014}
Bettina~C. Schwab, Hil~GE Meijer, Richard~JA van Wezel, and Stephan~A. van
  Gils.
\newblock Synchronization of the parkinsonian globus pallidus by gap junctions.
\newblock \emph{BMC Neuroscience}, 15\penalty0 (1):\penalty0 O17, July
  2014{\natexlab{b}}.
\newblock ISSN 1471-2202.
\newblock \doi{10.1186/1471-2202-15-S1-O17}.
\newblock URL \url{https://doi.org/10.1186/1471-2202-15-S1-O17}.

\bibitem[Singer(1957)]{singer_linear_1957}
I.~Singer.
\newblock Linear functionals on the space of continuous mappings of a compact
  {Hausdorff} space into a {Banach} spaces.
\newblock \emph{Rev. Math. Pures Appl.}, 2:\penalty0 301--315, 1957.
\newblock URL \url{https://ci.nii.ac.jp/naid/10009422054/}.

\bibitem[van Gils et~al.(2013)van Gils, Janssens, Kuznetsov, and
  Visser]{van_gils_local_2013}
S.~A. van Gils, S.~G. Janssens, Yu.~A. Kuznetsov, and S.~Visser.
\newblock On local bifurcations in neural field models with transmission
  delays.
\newblock \emph{Journal of Mathematical Biology}, 66\penalty0 (4):\penalty0
  837--887, March 2013.
\newblock ISSN 1432-1416.
\newblock \doi{10.1007/s00285-012-0598-6}.
\newblock URL \url{https://doi.org/10.1007/s00285-012-0598-6}.

\bibitem[van Gils et~al.(2012)van Gils, Janssens, Kuznetsov, and
  Visser]{van_gils_local_2012}
Stephan~A. van Gils, Sebastiaan~G. Janssens, Yuri~A. Kuznetsov, and Sid Visser.
\newblock On {Local} {Bifurcations} in {Neural} {Field} {Models} with
  {Transmission} {Delays}.
\newblock \emph{arXiv:1209.2849 [math]}, October 2012.
\newblock URL \url{http://arxiv.org/abs/1209.2849}.
\newblock arXiv: 1209.2849.

\bibitem[van Neerven(1990)]{van_neerven_reflexivity_1990}
J.~M. A.~M. van Neerven.
\newblock Reflexivity, the dual {Radon}-{Nikodym} property, and continuity of
  adjoint semigroups.
\newblock \emph{Indagationes Mathematicae}, 1\penalty0 (3):\penalty0 365--379,
  January 1990.
\newblock ISSN 0019-3577.
\newblock \doi{10.1016/0019-3577(90)90024-H}.
\newblock URL
  \url{http://www.sciencedirect.com/science/article/pii/001935779090024H}.

\bibitem[Veltz and Faugeras(2013)]{veltz_center_2013}
R.~Veltz and O.~Faugeras.
\newblock A {Center} {Manifold} {Result} for {Delayed} {Neural} {Fields}
  {Equations}.
\newblock \emph{SIAM Journal on Mathematical Analysis}, 45\penalty0
  (3):\penalty0 1527--1562, January 2013.
\newblock ISSN 0036-1410.
\newblock \doi{10.1137/110856162}.
\newblock URL \url{https://epubs.siam.org/doi/abs/10.1137/110856162}.

\bibitem[Veltz and Faugeras(2010)]{veltz_localglobal_2010}
Romain Veltz and Olivier Faugeras.
\newblock Local/{Global} {Analysis} of the {Stationary} {Solutions} of {Some}
  {Neural} {Field} {Equations}.
\newblock \emph{SIAM Journal on Applied Dynamical Systems}, 9\penalty0
  (3):\penalty0 954--998, January 2010.
\newblock \doi{10.1137/090773611}.
\newblock URL \url{https://epubs.siam.org/doi/abs/10.1137/090773611}.
\newblock Publisher: Society for Industrial and Applied Mathematics.

\bibitem[Veltz and Faugeras(2011)]{veltz_stability_2011}
Romain Veltz and Olivier Faugeras.
\newblock Stability of the stationary solutions of neural field equations with
  propagation delays.
\newblock \emph{The Journal of Mathematical Neuroscience}, 1\penalty0
  (1):\penalty0 1, May 2011.
\newblock ISSN 2190-8567.
\newblock \doi{10.1186/2190-8567-1-1}.
\newblock URL \url{https://doi.org/10.1186/2190-8567-1-1}.

\bibitem[Veltz and Faugeras(2015)]{veltz_erratum_2015}
Romain Veltz and Olivier Faugeras.
\newblock {ERRATUM}: {A} {Center} {Manifold} {Result} for {Delayed} {Neural}
  {Fields} {Equations}.
\newblock \emph{SIAM Journal on Mathematical Analysis}, 47\penalty0
  (2):\penalty0 1665--1670, January 2015.
\newblock ISSN 0036-1410.
\newblock \doi{10.1137/140962279}.
\newblock URL \url{https://epubs.siam.org/doi/abs/10.1137/140962279}.
\newblock Publisher: Society for Industrial and Applied Mathematics.

\bibitem[Venkov et~al.(2007)Venkov, Coombes, and Matthews]{venkov_dynamic_2007}
N.~A. Venkov, S.~Coombes, and P.~C. Matthews.
\newblock Dynamic instabilities in scalar neural field equations with
  space-dependent delays.
\newblock \emph{Physica D: Nonlinear Phenomena}, 232\penalty0 (1):\penalty0
  1--15, August 2007.
\newblock ISSN 0167-2789.
\newblock \doi{10.1016/j.physd.2007.04.011}.
\newblock URL
  \url{http://www.sciencedirect.com/science/article/pii/S0167278907001285}.

\bibitem[Visser et~al.(2017)Visser, Nicks, Faugeras, and
  Coombes]{visser_standing_2017}
S.~Visser, R.~Nicks, O.~Faugeras, and S.~Coombes.
\newblock Standing and travelling waves in a spherical brain model: {The}
  {Nunez} model revisited.
\newblock \emph{Physica D: Nonlinear Phenomena}, 349:\penalty0 27--45, June
  2017.
\newblock ISSN 0167-2789.
\newblock \doi{10.1016/j.physd.2017.02.017}.
\newblock URL
  \url{http://www.sciencedirect.com/science/article/pii/S0167278916306352}.

\bibitem[Webb(1976)]{webb_functional_1976}
G.~F Webb.
\newblock Functional differential equations and nonlinear semigroups in
  {Lp}-spaces.
\newblock \emph{Journal of Differential Equations}, 20\penalty0 (1):\penalty0
  71--89, January 1976.
\newblock ISSN 0022-0396.
\newblock \doi{10.1016/0022-0396(76)90097-8}.
\newblock URL
  \url{http://www.sciencedirect.com/science/article/pii/0022039676900978}.

\bibitem[Wilson and Cowan(1973)]{wilson_mathematical_1973}
H.~R. Wilson and J.~D. Cowan.
\newblock A mathematical theory of the functional dynamics of cortical and
  thalamic nervous tissue.
\newblock \emph{Kybernetik}, 13\penalty0 (2):\penalty0 55--80, September 1973.
\newblock ISSN 1432-0770.
\newblock \doi{10.1007/BF00288786}.
\newblock URL \url{https://doi.org/10.1007/BF00288786}.

\bibitem[Wilson and Cowan(1972)]{wilson_excitatory_1972}
Hugh~R. Wilson and Jack~D. Cowan.
\newblock Excitatory and {Inhibitory} {Interactions} in {Localized}
  {Populations} of {Model} {Neurons}.
\newblock \emph{Biophysical Journal}, 12\penalty0 (1):\penalty0 1--24, January
  1972.
\newblock ISSN 0006-3495.
\newblock \doi{10.1016/S0006-3495(72)86068-5}.
\newblock URL
  \url{http://www.sciencedirect.com/science/article/pii/S0006349572860685}.

\bibitem[Wu(2012)]{wu_theory_2012}
Jianhong Wu.
\newblock \emph{Theory and {Applications} of {Partial} {Functional}
  {Differential} {Equations}}.
\newblock Springer Science \& Business Media, December 2012.
\newblock ISBN 978-1-4612-4050-1.

\end{thebibliography}
\end{document}